\documentclass{article}
\usepackage{graphicx}
\usepackage{amsmath}
\usepackage{amssymb}
\usepackage{amsthm}
\usepackage{mathrsfs}
\usepackage{here}
\providecommand{\bysame}{\leavevmode\hbox to3em{\hrulefill}\thinspace}
\providecommand{\MR}{\relax\ifhmode\unskip\space\fi MR }

\numberwithin{equation}{section}

\newtheorem{dummy}{Dummy}[section]

\newtheorem{thm}[dummy]{Theorem}

\newtheorem{observation}[dummy]{Observation}
\theoremstyle{definition}

\newtheorem{rem}[dummy]{Remark}

\newtheorem{prob}[dummy]{Problem}

\let\emx=\relax

\begin{document}

\title{Numerical studies of the optimization of the first eigenvalue of the heat diffusion in inhomogeneous media}

\author{Kaname Matsue\footnote{The Institute of Statistical Mathematics, Tachikawa, Tokyo, 190-8562, Japan. {\ttfamily kmatsue@ism.ac.jp}}
  \,
  and 
  Hisashi Naito\footnote{Graduate school of Mathematics, Nagoya University, Nagoya, 464-8602, Japan. {\ttfamily naito@math.nagoya-u.ac.jp}}}
\date{}
\maketitle
\begin{abstract}
  In this paper, 
  we study optimization of the first eigenvalue of 
  $-\nabla \cdot (\rho(x) \nabla u) = \lambda u$
  in a bounded domain $\Omega \subset \mathbb{R}^n$
  under several constraints for the function $\rho$. 
  We consider this problem in various boundary conditions and various topologies of domains.
  As a result,
  we numerically observe several common criteria for $\rho$ for optimizing eigenvalues in terms of corresponding eigenfunctions,
  which are independent of topology of domains and boundary conditions.
  Geometric characterizations of optimizers are also numerically observed. 
\end{abstract}
\section{Introduction}
\label{section-intro}

In this paper,
we consider the following eigenvalue problem
\begin{equation}
  \label{ev-prob-rho}
  -\nabla \cdot (\rho(x) \nabla u) = \lambda u
  \quad x\in \Omega,\quad \rho \in \mathcal{K}
\end{equation}
in a bounded domain $\Omega$ in $\mathbb{R}^n$ with suitable boundary conditions,
where  
\begin{equation}
  \label{K}
  \begin{aligned}
    \mathcal{K} &:= \left\{ \rho\in L^\infty(\Omega)\mid \rho = 1\text{ or }c,
      \text{ a.e.\ on }\Omega,\ \int_\Omega \rho dx = (cm_0 + (1-m_0))|\Omega| \right\},\\
    &c > 0,
    \quad
    c \not= 1, 
    \quad
    m_0 \in (0,1)
  \end{aligned}
\end{equation}
and $|\Omega|$ denotes the Lebesgue measure of $\Omega$ in $\mathbb{R}^n$. 
Let $S := \{x\in \Omega \mid \rho(x) = c\}$, 
then (\ref{K}) immediately implies 
\begin{equation}
  \label{constraint}
  |S| / |\Omega| \equiv m_0.
\end{equation}
\par
Under these settings, we consider the following problem:
\begin{prob}
  \label{prob-rho-abst}
  Find $\rho_\ast \in \mathcal{K}$ 
  which attains the supremum or the infimum of the first eigenvalue 
  $\lambda_1(\rho)$ on $\mathcal{K}$ 
  under suitable boundary condition.
  If it exists, 
  characterize the shape of the domain $S_\ast = \{x\in \Omega \mid \rho_\ast(x) = c\}$.
\end{prob}
We shall call the smallest positive eigenvalue the first eigenvalue, 
since we consider this problem under Dirichlet, Neumann or mixed boundary condition.
\par
Our motivation for the above problem is as follows.
Assume that two different materials are in a given domain $\Omega$
with fixed volume ratio.
How do we arrange such materials to optimize the heat conductivity of $\Omega$? 
Since the long-time behaviour of heat transfer of the heat equation 
is controlled by the first eigenvalue of 
\begin{math}
  -\nabla \cdot (\rho(x)\nabla u)
\end{math}, 
we may consider Problem \ref{prob-rho-abst} as one of toy models of this problem.
\par
In this paper,
we numerically study the eigenvalue optimization for (\ref{ev-prob-rho}) 
and we observe the following results for Problem \ref{prob-rho-abst}:
\begin{enumerate}
\item 
  The element $\rho_\ast \in \mathcal{K}$ optimizing $\lambda_1(\rho)$ can be characterized by inequalities with respect to $\rho_\ast \nabla u_\ast$, 
  where $u_\ast$ is the eigenfunction associated with $\lambda_1(\rho_\ast)$ of (\ref{ev-prob-rho}).
  As a consequence, 
  the domain $S_\ast = \{x\in \Omega \mid \rho_\ast(x) = c\}$ is given by the super- or the sub-level set of $|\rho_\ast \nabla u_\ast|$.
  This characterization is independent of topology and geometry of $\Omega$ and boundary conditions on $\partial \Omega$. 
\item 
  If $\Omega$ is star-shaped and symmetric in a certain direction,
  then $S_\ast$ has the same symmetry.
\item 
  Optimized region $S_\ast$ depends continuously on a parameter of the boundary condition if the Robin boundary condition is imposed.
\end{enumerate}
The precise statements are described in Section \ref{section-numerical} 
(Observation \ref{observation-criteria-heat} - \ref{observation-dependence-heat}) with various numerical results.
\par
This paper is organized as follows.
In Section \ref{section-setting},
we provide more precise setting of our problems. 
A numerical method for finding optimizers we apply here, {\em the level set approach}, is also derived here. 
In Section \ref{section-prior},
numerical and mathematical known results for a well-considered problem are discussed. 
Section \ref{section-numerical} is where our main discussion is developed. 
We show several numerical observations about eigenvalue optimization criteria, 
geometry of the level set $S_\ast = \{x\in \Omega \mid \rho_\ast(x) = c\}$ 
for the optimizer $\rho_\ast$ and continuous dependence of $S_\ast$ on boundary conditions.

\section{Setting}
\label{section-setting}
\subsection{Setting of problems}
Here we provide the precise setting of Problem \ref{prob-rho-abst}.
\begin{prob}[Precise version of Problem \ref{prob-rho-abst}] 
  \label{prob-rho}
  For $\rho \in \mathcal{K}$ define
  \begin{equation}
    \label{ev-rho}
    \lambda_1(\rho) := \inf_{u \text{ : admissible}} \frac{\int_\Omega \rho(x) |\nabla u(x)|^2 dx}{\int_\Omega |u(x)|^2 dx}.
  \end{equation}
  Find an element $\rho_\ast \in \mathcal{K}$ which attains
  \begin{equation}
    \label{sup-rho}
    \sup_{\rho\in \mathcal{K} } \lambda_1(\rho).
  \end{equation}
  If it exists, characterize $\rho_\ast$ and 
  the shape of domain $S_\ast = \{x\in \Omega \mid \rho_\ast(x) = c\}$.
\end{prob}
We also consider the same question for $\inf_\rho \lambda_1(\rho)$.
\par
It is mathematically known that the first eigenvalue of the linear operator 
$A_\rho u = -\nabla \cdot (\rho(\cdot)\nabla u)$
is characterized by the Rayleigh quotient (\ref{ev-rho}) (see e.g.\ \cite{CH}). 
Here,
we also consider the minimization of $\lambda_1(\rho)$ as a comparison with the maximization of it.
\par
We consider the problem in the cases of various type of domains like ones with piecewise smooth boundary, 
non-convex or non-simply connected ones,
as well as various boundary conditions (Dirichlet, Neumann or mixed boundary condition).
\par
\begin{rem}
  A well-known mathematical theory (see e.g.\ \cite{Ban}) tells us that 
  the linear operator $A_\rho$
  associated with 
  (\ref{ev-prob-rho}) possesses discrete eigenvalues
  \begin{equation}
    \label{evs-prob-rho}
    (0 <)\ \lambda_1(\rho) \leq \lambda_2(\rho) \leq  \cdots \to \infty.
  \end{equation}
  In the case of the homogeneous Dirichlet boundary value problem, 
  $0$ is not an eigenvalue and hence the infimum in (\ref{ev-rho}) is attained in $u$ being not identically zero. 
  Thus, in the case of Dirichlet boundary value problems, 
  we call $u$ admissible if and only if 
  $u \not= 0$ in an appropriate function space.
  On the other hand, in the case of the homogeneous Neumann boundary value problem, 
  $0$ is admitted as an eigenvalue with a constant function as the eigenfunction.
  $0$ is thus the smallest eigenvalue.
  According to the general theory of eigenfunctions, 
  eigenfunctions associated with $\lambda_1(\rho)$ have to be orthogonal to constant functions in the inner product on $L^2(\Omega)$.
  We thus call $u$ 
  admissible, if and only if 
  $u$ satisfies $\int_\Omega u\,dx = 0$.
\end{rem}
\subsection{The level set method}
\label{section-levelset}
Our problem is one of typical problems called {\em topology optimization}.
One of well-known methods for determining topology of the optimal object is {\emx the level set approach}, 
originally developed in Osher and Sethian \cite{OS1}.
This method provides an efficient way of describing time-evolving curves and surfaces which may undergo topological changes.
Osher and Santosa \cite{OS} improve the method so that it can be applied to optimization problems with one or more constraints like volume constraint.
We follow their method as a numerical approach for determining optimizers in our problems. Here we briefly review implementations of the level set approach discussed in \cite{OS}.
\par
We describe a subset $S\subset \Omega$ as the super-level set of a function $\phi \colon \Omega\to \mathbb{R}$, i.e.\ 
\begin{equation*}
  S = \{x\in \Omega \mid \phi(x) > 0\} = \{x\in \Omega \mid \rho(x) = c\}.
\end{equation*}
Then, the function $\rho$ in Problem \ref{prob-rho} can be represented as a function of $\phi$,
and hence $\lambda_1$ in Problem \ref{prob-rho} can be also regarded as a function of $\phi$ since $\lambda_1$ originally depends on $\rho$.
We may thus rewrite our optimization problem in the minimization problem of the following energy functional $L$:

\begin{equation}
  \label{Lagrange}
  L(\phi) = F(\phi) + \nu G(\phi),
\end{equation}
where $F(\phi) = \lambda_1(\phi) := \lambda_1(\rho)$ for minimizing eigenvalues and  $F(\phi) = -\lambda_1(\phi)$ for maximizing eigenvalues,
$\nu$ is the Lagrange multiplier and 
$G(\phi) = \int_{\{\phi(x) > 0\}}  dx - m_0|\Omega|$.
Here, we remark that the equation $G(\phi) = 0$ corresponds to the volume constraint (\ref{constraint}).
\par
The first variation of $L$ and $\nu$ can be calculated by ordinary methods. 
Our optimization problem is then reduced to solving the gradient flow associated with $L$,
namely, the flow which decreases the energy $L$ along solutions.
More precisely, our problem is reduced to the problem to solve the following evolution equation: 
\begin{equation}
  \label{SDF}
  \frac{\partial \phi}{\partial t} = -(v_0 + \nu)|\nabla \phi|\quad \text{ on }\partial S 
\end{equation}
in a certain function space, where $v_0 = v_0(x)$ is given by
\begin{equation}
  \label{v0}
  v_0(x) = \frac{c-1}{\int_\Omega u_\phi^2dx }|\nabla u_\phi(x)|^2,
\end{equation}
and $u_\phi$ is the associated eigenfunction of $\lambda_1(\phi)$.
Here we omit the detail how to obtain this system,
since it is completely done by following the arguments in \cite{OS}.
\par
\begin{rem}
  Since the differentiability of the solution $\phi$ of (\ref{SDF}) may be lost during evolution of (\ref{SDF}),
  we may not construct numerical solution of (\ref{SDF}). To avoid this difficulty,
  instead of (\ref{SDF}), 
  we solve
  \begin{equation}
    \label{viscos}
    \frac{\partial \phi}{\partial t} = \epsilon \Delta \phi -(v_0 + \nu)|\nabla \phi|
  \end{equation}
  with sufficiently small $\epsilon > 0$.
  Then (\ref{viscos}) is a semi-linear parabolic evolution equation,
  and hence the well-known explicit or implicit scheme enables us to solve this equation numerically keeping smoothness of $\phi$ during evolutions.
  This technique is well-known as {\em viscosity vanishing method}.
  As for Hamilton-Jacobi type equations like (\ref{SDF}), 
  it is known that solutions of (\ref{viscos}) approach to those of (\ref{SDF}) as $\epsilon \to 0$ in a suitable sense (see e.g.\ \cite{CL}). 
\end{rem}

\section{Known results}
\label{section-prior}
As for eigenvalue optimization problems, 
there are a lot of earlier works for
\begin{equation}
  \label{ev-prob-sigma}
  -\Delta u = \mu \sigma(x) u,\quad x\in \Omega,\quad \sigma \in \mathcal{K}.
\end{equation}
This is a generalized eigenvalue problem of the Laplacian, 
which is, for example,
well-considered for studying frequency of drums with spatially inhomogeneous density on $\Omega$.
Osher-Sethian \cite{OS1} and Osher-Santosa \cite{OS}
considered
(\ref{ev-prob-sigma}) with the level set approach (see Section \ref{section-levelset}), 
which has been 
well treated during a variety of improvements.
In \cite{OS1, OS} and many related works, 
the following problem is considered, 
which is an alternative one of Problem \ref{prob-rho}.
\begin{prob}
  \label{prob-sigma}
  For $\sigma \in \mathcal{K}$ define
  \begin{equation}
    \label{ev-sigma}
    \mu_1(\sigma) := \inf_{u \text{ : admissible}} \frac{\int_\Omega |\nabla u(x)|^2 dx}{\int_\Omega \sigma(x) |u(x)|^2 dx}.
  \end{equation}
  Find an element $\sigma_\ast \in \mathcal{K}$ which attains
  \begin{equation}
    \label{inf-sigma}
    \inf_{\sigma \in \mathcal{K}}  \mu_1(\sigma).
  \end{equation}
  If it exists, characterize $\sigma_\ast$ and
  the shape of domain $S_\ast = \{x\in \Omega \mid \sigma_\ast(x) = c\}$.
\end{prob}
\par
For Problem \ref{prob-sigma}, 
there are several mathematical results which determine optimizers of $\mu_1(\sigma)$ as well as geometries of domain $S_\ast = \{x\in \Omega \mid \sigma_\ast(x) = c\}$.
\par
In \cite{Kre}, 
Krein mathematically considered the largest and the smallest $k$-th natural frequency of 
strings with fixed endpoints on the interval $(0,1)$.
In words of Problem \ref{prob-sigma}, 
he considered $\mu_k$, $k\in \mathbb{N}$, 
on $\Omega = (0,1)$ with the homogeneous Dirichlet boundary condition. 
He completely solved this problem by detecting the optimizer $\sigma_\ast$.
A couple of decades later, 
Cox and McLaughlin extended Krein's arguments in arbitrary dimensions \cite{CM, CM2}.
They obtained the following optimization criteria for eigenvalues in terms of corresponding eigenfunctions. 
\par
\begin{thm}[Cox and McLaughlin \cite{CM2}]
  \label{preceding-sigma}
  Consider Problem \ref{prob-sigma} with the homogeneous Dirichlet boundary condition.
  If $\sigma_{\min}\in \mathcal{K}$ is the minimizer of $\mu_1$ in $\mathcal{K}$,
  then there is a positive constant $\alpha > 0$ such that,
  for the corresponding eigenfunction $u_1$, the following inequalities hold:
  \begin{align*}
    \sigma_{\min}(x) = c \quad &\Rightarrow \quad u_1(x) \geq \alpha,\\
    \sigma_{\min}(x) = 1 \quad &\Rightarrow \quad u_1(x) \leq \alpha.
  \end{align*}
  Similarly if $\sigma_{\max}\in \mathcal{K}$ is the maximizer of $\mu_1$ in $\mathcal{K}$,
  then there is a positive constant $\alpha' > 0$ such that the following inequalities hold:
  \begin{align*}
    \sigma_{\max}(x) = c \quad &\Rightarrow \quad u_1(x) \leq \alpha',\\
    \sigma_{\max}(x) = 1 \quad &\Rightarrow \quad u_1(x) \geq \alpha'.
  \end{align*}
\end{thm}
Moreover, 
they obtained results for geometric properties of $S_\ast$ using symmetry arguments derived from the maximum principle.
\begin{thm}[Cox and McLaughlin \cite{CM2}]
  \label{preceding-geom}
  Consider Problem \ref{prob-sigma} with the homogeneous Dirichlet boundary condition.
  Assume that $\sigma_{\min} \in \mathcal{K}$ is the minimizer of $\mu_1$ in $\mathcal{K}$ and $S_{\min}  = \{x\in \Omega \mid \sigma_{\min}(x) = c\}$.
  If $\Omega$ is convex and symmetric in $N$ orthogonal directions,
  then 
  $S_{\min}$ is also convex and symmetric in $N$ directions.
  Moreover, $S_{\min}$ is star-shaped with respect to the center of symmetry.
  Similarly, if $\sigma_{\max} \in \mathcal{K}$ is the maximizer of $\mu_1$ in $\mathcal{K}$ and $S_{\max}  = \{x\in \Omega \mid \sigma_{\max}(x) = c\}$.
  Then, under the same assumption as minimizers,
  the same statements hold for $S_{\max}^c$,
  where $A^c := \Omega\setminus  A$ is the complement of $A\subset \Omega$.
\end{thm}
Theorem \ref{preceding-geom} claims that geometric properties of $S_{\min}$ and $S_{\max}$ are also characterized by corresponding eigenfunctions,
since $S_{\min}$ and $S_{\max}$ are determined by the super- or the sub-level sets of corresponding eigenfunctions.
On the other hand, 
these researches focus on Problem \ref{prob-sigma} only with the homogeneous Dirichlet boundary condition.
Furthermore \cite{CM}, \cite{CM2} and \cite{Kre}
 are not concerned with Problem \ref{prob-sigma} with any other boundary condition, nor are they concerned with Problem \ref{prob-rho}. 
And we know no paper about them (even numerically) in the last few decades.

\begin{rem}
  Lou and Yanagida \cite{LY} discuss an indefinite linear eigenvalue problem related to biological invasion of species. 
  The formation of their problem is similar to Problem \ref{prob-sigma} with the homogeneous Neumann condition,
  but the assumption on the weight $\sigma(x)$ is different. 
  In \cite{LY}, 
  $\sigma$ is assumed to be bounded with a fixed negative total weight.
  They consider only one-dimensional problems and prove that the global minimizer $\sigma_\ast$ of the principal eigenvalue $\mu_1$ 
  (i.e.\ the first eigenvalue in our arguments) should be the specific two-valued function,
  which is often called \lq\lq bang-bang" type.
  Such bang-bang type weight $\sigma_\ast$ is characterized by the super-level set 
  of the eigenfunction associated with $\mu_1(\sigma_\ast)$, which is similar to Theorem \ref{preceding-sigma}.
\end{rem}
\par
In \cite{CGIKO}, Chamillo, Greiser, Imai, Kurata and Onishi 
discuss geometries of the minimizing configuration $S_\ast$ for the first eigenvalue problem 
$-\Delta u + (c-1) \chi_{S} u = \lambda u$ under the volume constraint (\ref{constraint}).
They prove that $S_\ast$ does not inherit symmetry of $\Omega$ for appropriate values of $c$ and $m_0$ when $\Omega$ is either an annulus or a dumbbell.
\par
Problem \ref{prob-rho} is considered by Cox and Lipton \cite{Cox-Lipton}, 
Conca, Laurain and Mahadevan \cite{Conca-Laurain-Mahadevan}, 
Conca, Mahadevan and Sanz \cite{Conca-Mahadevan-Sanz} and some other related works. 
In \cite{Cox-Lipton},
Cox and Lipton study an shape optimization problem for the heat equation with the Dirichlet boundary condition, which corresponds to Problem \ref{prob-rho}.
They discuss the existence of optimal eigenvalues and optimizers. Moreover, they obtain the exact optimal configuration in the one-dimensional problem.
On the other hand, in Conca, Laurain and Mahadevan \cite{Conca-Laurain-Mahadevan} and Conca, Mahadevan and Sanz \cite{Conca-Mahadevan-Sanz}, 
they discuss the asymptotic expansion of the first eigenvalue $\lambda$ and the asymptotic behavior of the minimized eigenvalue.
They also discuss a conjecture of the optimal configuration when $\Omega$ is a ball.
The conjecture claims that the optimal configuration $S_\ast$ for the minimization problem of the first Dirichlet eigenvalue is a concentric ball.
However, the conjecture is disproved in \cite{Conca-Dambrine-Mahadevan-Quintero}.
At least in two and three dimensional cases, the optimal configuration is a union of concentric annuli and a disk under a certain condition of $c$ and $m_0$.
On the other hand, 
since all the discussions therein are based on the asymptotic expansion of eigenpairs, hence the characterization of optimal configurations remains as an open problem for large $c$. 
\par
We also note that the attainability of the supremum or the infimum of the first positive eigenvalue in Problem \ref{prob-rho} still remains open.
In fact, the subset $\mathcal{K}$ is {\emx not} weak${}^\ast$-closed in $L^\infty(\Omega)$ (e.g.\ \cite{Cox-Lipton, CM}).
The supremum and the infimum of the first eigenvalue are thus attained in the weak${}^\ast$-closure $\mathcal{K}^\ast$ of $\mathcal{K}$.
In \cite{Cox-Lipton}, 
they prove that, in Problem \ref{prob-rho}, 
$\inf_{\rho \in \mathcal{K}}\lambda(\rho) = \inf_{\rho \in \mathcal{K}^\ast}\lambda(\rho)$ and the similar result for supremum also holds.
But there is still a possibility that $\rho_\ast \in \mathcal{K}^\ast \setminus \mathcal{K}$.
We finally remark that the optimizer $\sigma_\ast$ in Problem \ref{prob-sigma} can be chosen as an element of $\mathcal{K}$ under 
sufficient smoothness assumptions on $\partial \Omega$ and eigenfunctions.
This fact is the consequence of \cite{CM} and \cite{CM2}.

\section{Numerical Study}
\label{section-numerical}
\subsection{Optimization criteria}
\label{section-prob-rho}
So far Problem \ref{prob-rho} is well-considered neither mathematically nor numerically, but the same level set approach as in the case of Problem \ref{prob-sigma} can be applied.
\par
Recall that optimizers of eigenvalues are determined by stationary solutions of (\ref{SDF}).
 If $\phi$ is a stationary solution of (\ref{SDF}), then $\phi$ should satisfy $-(v_0 + \nu)|\nabla \phi| = 0$ on $\partial S = \phi^{-1}(0)$. 
Moreover if $\phi$ satisfies $|\nabla \phi| \not = 0$ on $\partial S$,
then it should associate $\rho_\ast \in \mathcal{K}$, 
which is given by $\rho_\ast(x) = c$ 
if and only if $\phi(x) > 0$, and the eigenfunction $u_\ast$ of (\ref{ev-prob-rho}) 
with $\rho = \rho_\ast$ such that $|\nabla u_\ast| = \text{ constant}$ on $\partial S$, 
since $v_0 = v_0(x) = C|\nabla u_\ast(x)|^2$ is equal to the constant $\nu$ on $\partial S$, 
where $C\not = 0$ is a constant. 
Hence one can guess that optimizers are characterized by $|\nabla u_\ast|$.
\par
\begin{figure}[H]
  \begin{center}
    \begin{tabular}{llll}
      {\bfseries (a)} 
      &
      {\bfseries (b)} 
      & 
      {\bfseries (c)}
      & 
      {\bfseries (d)}
      \cr
      \includegraphics[width=3.0cm]{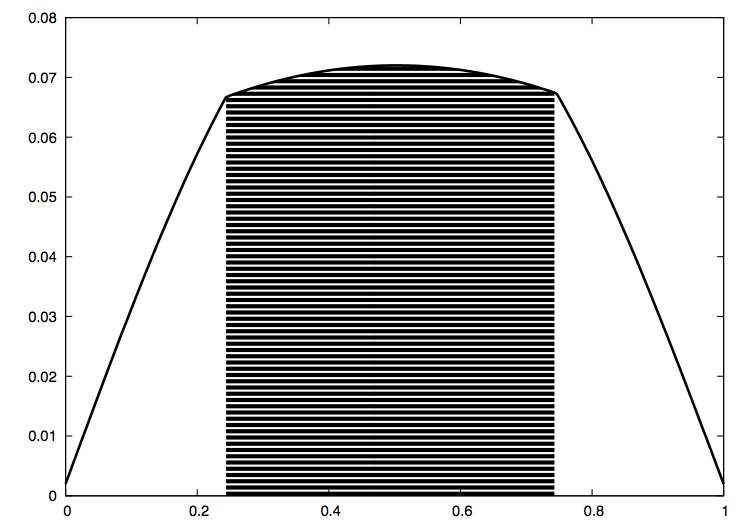}
      &
      \includegraphics[width=3.0cm]{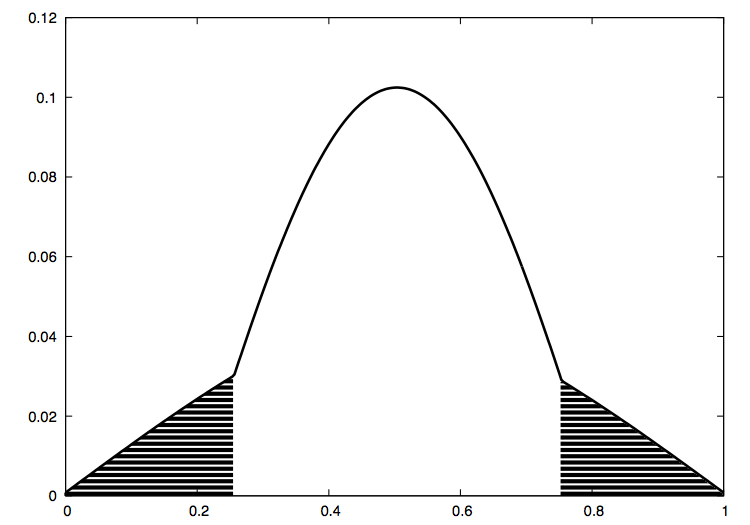}
      &
      \includegraphics[width=3.0cm]{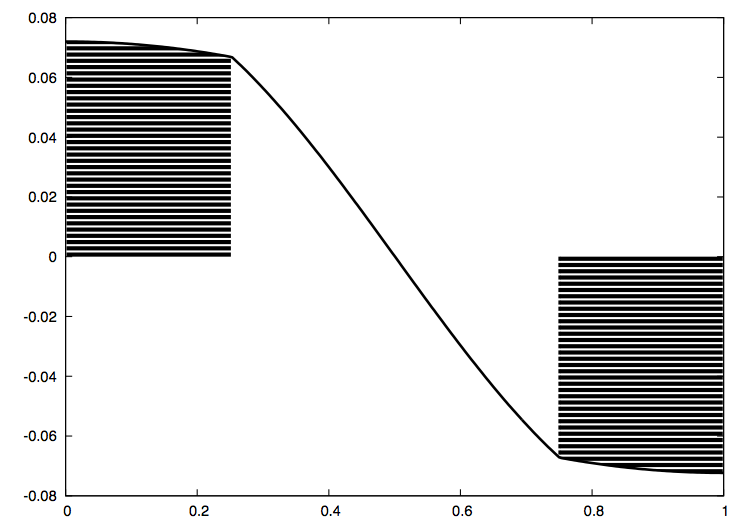}
      &
      \includegraphics[width=3.0cm]{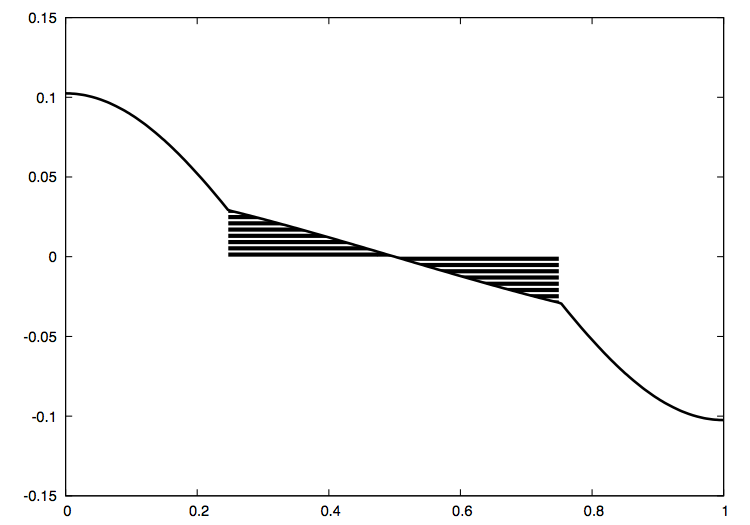}
      \cr
    \end{tabular}
  \end{center}
  \caption{Optimal regions and eigenfunctions for Problem \ref{prob-rho}}
  (a) : Minimization of $\lambda_1(\rho)$ with the homogeneous Dirichlet boundary condition. 
  (b) : Maximization of $\lambda_1(\rho)$ with the homogeneous Dirichlet boundary condition. 
  (c) : Minimization of $\lambda_1(\rho)$ with the homogeneous Neumann boundary condition. 
  (d) : Maximization of $\lambda_1(\rho)$ with the homogeneous Neumann boundary condition. 
  The region in $\Omega_0$ where impulses are hung on is $S_\ast$ in each figure. 
  Computed eigenvalues with $c=5$ are $(a) : 11.158517, (b) : 26.563359, (c) : 11.487066$ and $(d) : 27.073145$.
  One can see that $\partial S$ corresponds to the discontinuity of the differential $u'$ of $u$. 
  \label{fig-rho1dim-1}
\end{figure}
\begin{figure}[H]
  \begin{center}
    \begin{tabular}{llll}
      {\bfseries (a)} 
      &
      {\bfseries (b)} 
      & 
      {\bfseries (c)}
      & 
      {\bfseries (d)}
      \cr
      \includegraphics[width=3.0cm]{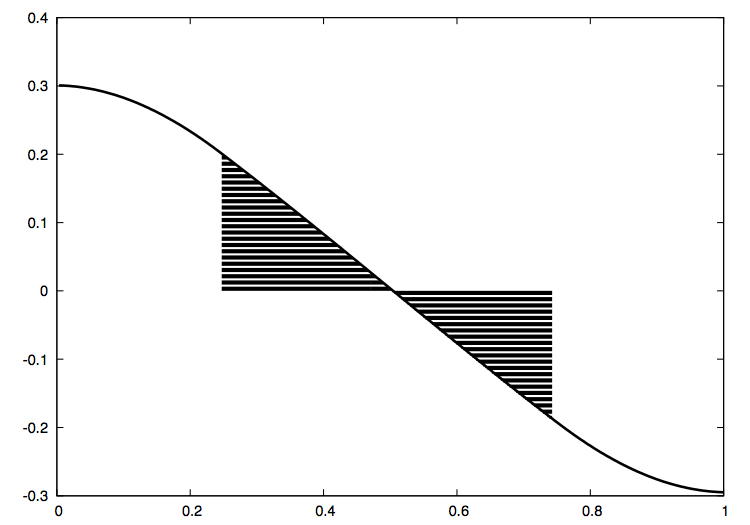}
      &
      \includegraphics[width=3.0cm]{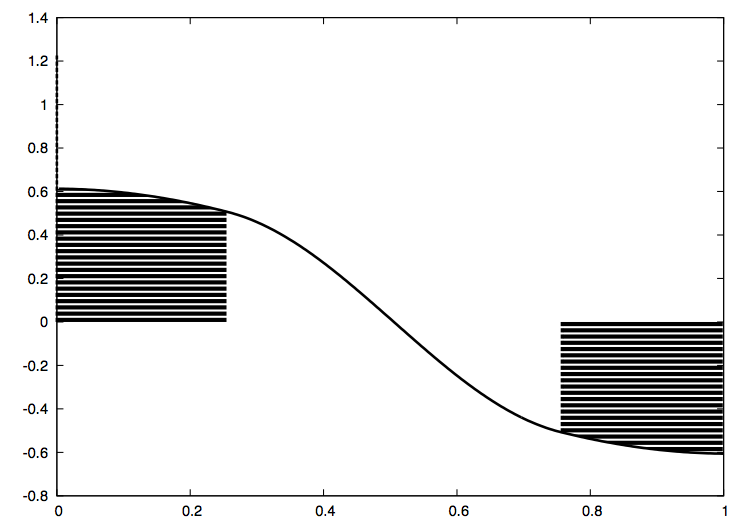}
      &
      \includegraphics[width=3.0cm]{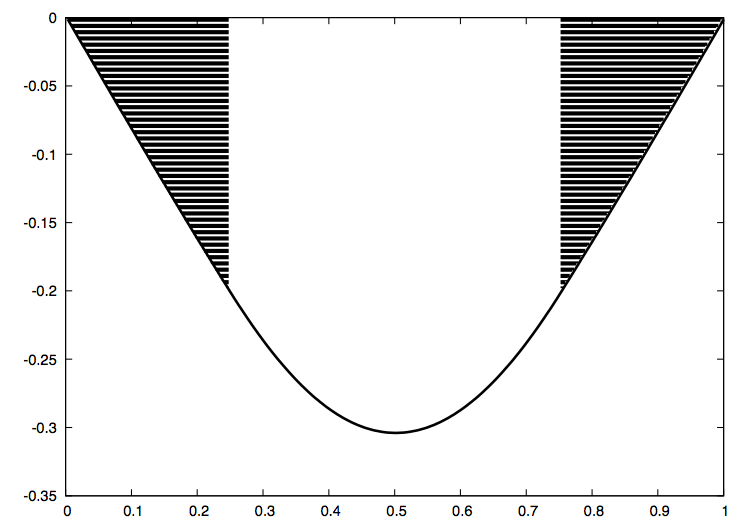}
      &
      \includegraphics[width=3.0cm]{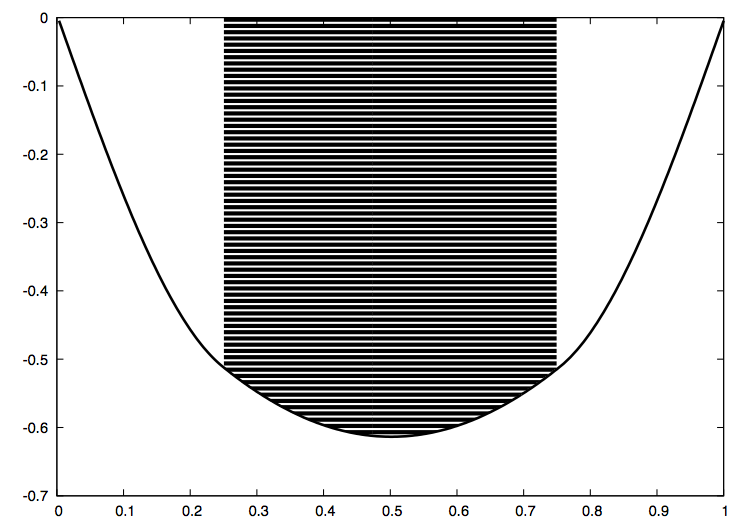}
      \cr
    \end{tabular}
  \end{center}
  \caption{The graph of $\rho u'$ of associated eigenfunction $u$ of the corresponding symbol in Figure \ref{fig-rho1dim-1}.}
  (a) is the graph of $\rho u'$ of $u$ in Figure \ref{fig-rho1dim-1}-(a).
  The region in $\Omega_0$ where impulses are hung on is $S_\ast$.
  The rest of figures are drawn in the same manner.
  One can expect that there is a certain correspondence between $S_\ast$ and the super- or the sub-level set of $\rho u'$. 
  \label{fig-rho1dim-2}
\end{figure}

First we study the eigenvalue problem in 
$\Omega_0 := (0,1)\subset \mathbb{R}$ shown in Figures \ref{fig-rho1dim-1} and \ref{fig-rho1dim-2}. 
As shown in Figure \ref{fig-rho1dim-1}, 
the differential $u' (=\nabla u = du/dx)$ of the eigenfunction $u$ has discontinuities on $\partial S$ 
and there is little hope of the correspondence between $S$ for optimizers and $u'$.
Next we consider the correspondence between $S$ for optimizers and $\rho u'$ instead of $u'$.
The graph of $\rho u'$ and $S$ for optimizers are shown in Figure \ref{fig-rho1dim-2}.
Both in the case of Dirichlet and Neumann boundary conditions, 
we may guess a correspondence between $S$ for optimizers and level sets of $\rho u'$ as similar to Theorem \ref{preceding-sigma}.
\par
\begin{figure}[H]
  \begin{center}
    \begin{tabular}{llll}
      {\bfseries (a)} 
      &
      {\bfseries (b)} 
      & 
      {\bfseries (c)}
      & 
      {\bfseries (d)}
      \cr
      \includegraphics[width=3.0cm]{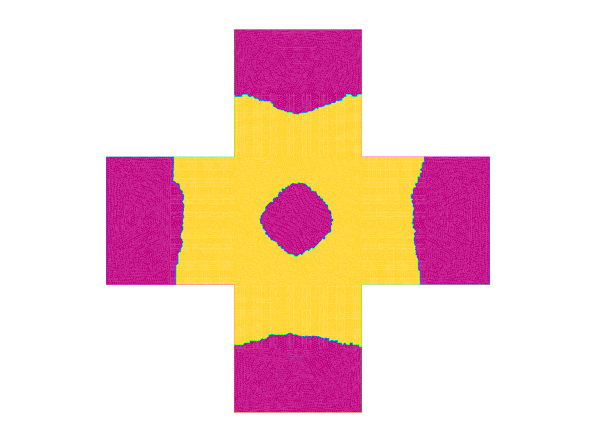}
      &
      \includegraphics[width=3.0cm]{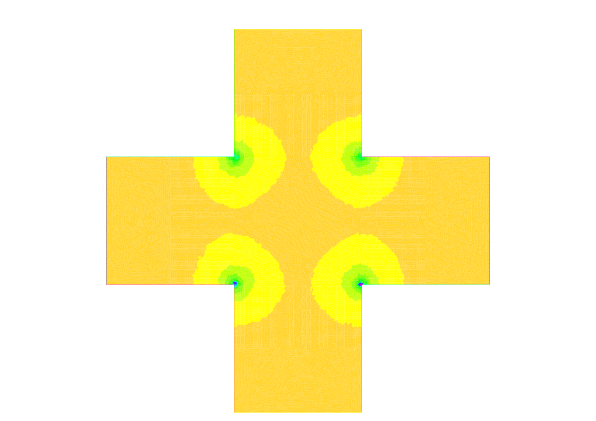}
      &
      \includegraphics[width=3.0cm]{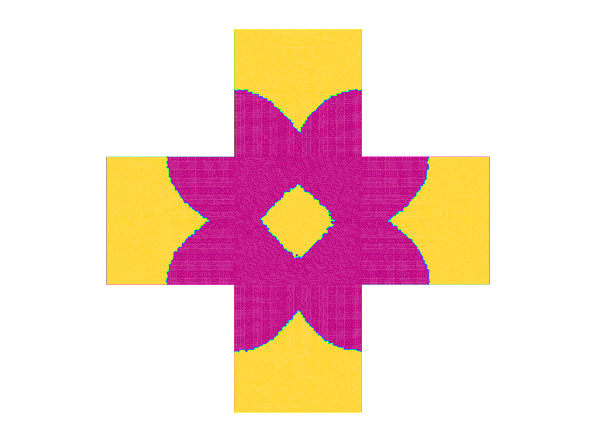}
      &
      \includegraphics[width=3.0cm]{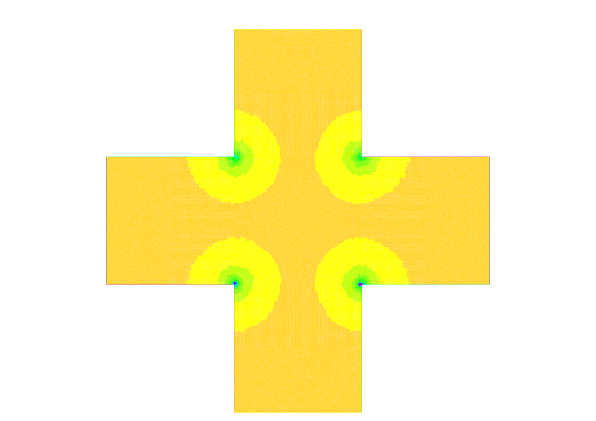}
      \cr
    \end{tabular}
  \end{center}
  \caption{Optimizers and eigenfunctions for Problem \ref{prob-rho} on $\Omega_2$}
  (a) minimizer $\rho_{\min}$, 
  (b) $|\rho \nabla u|^2$ of the associated eigenfunction of $\lambda_1(\rho_{\min})$.
  (c) maximizer $\rho_{\max}$, 
  (d) $|\rho \nabla u|^2$ of the associated eigenfunction of $\lambda_1(\rho_{\max})$.
 \label{fig-rhoD-cross}
\end{figure}
\begin{figure}[H]
  \begin{center}
    \begin{tabular}{llll}
      {\bfseries (a)} 
      &
      {\bfseries (b)} 
      & 
      {\bfseries (c)}
      & 
      {\bfseries (d)}
      \cr
      \includegraphics[width=3.0cm]{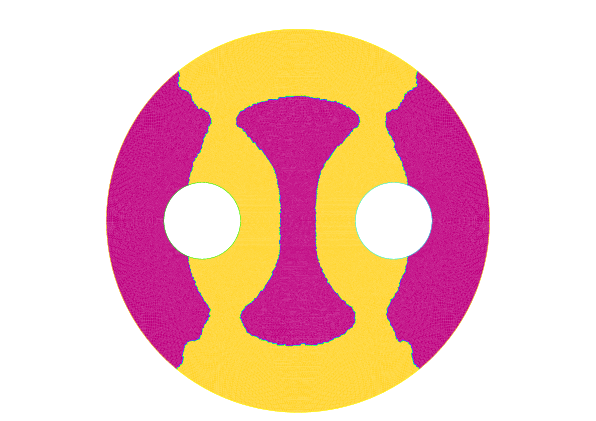}
      &
      \includegraphics[width=3.0cm]{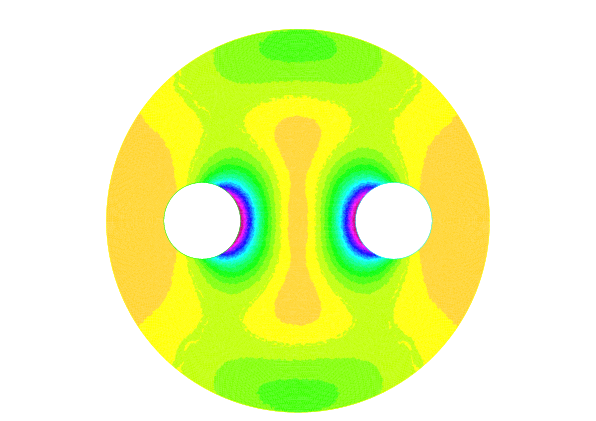}
      &
      \includegraphics[width=3.0cm]{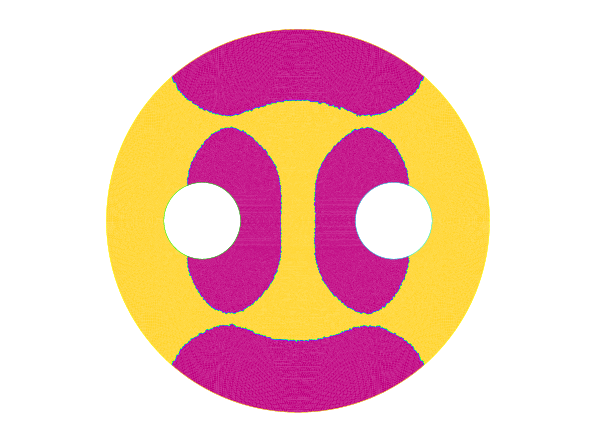}
      &
      \includegraphics[width=3.0cm]{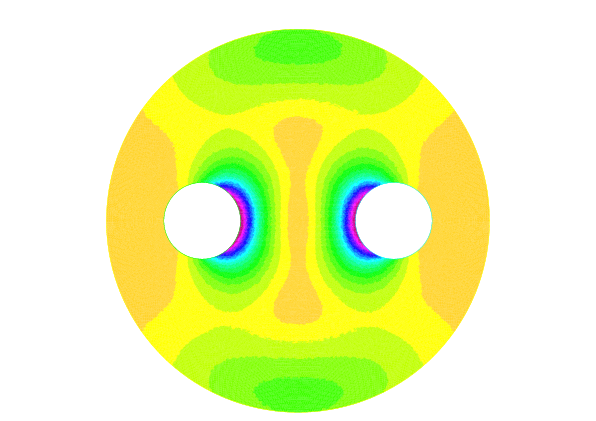}
      \cr
    \end{tabular}
  \end{center}
  \caption{Optimizers and eigenfunctions for Problem \ref{prob-rho} on $\Omega_4$}
  (a) minimizer $\rho_{\min}$, 
  (b) $|\rho \nabla u|^2$ of the associated eigenfunction of $\lambda_1(\rho_{\min})$.
  (c) maximizer $\rho_{\max}$, 
  (d) $|\rho \nabla u|^2$ of the associated eigenfunction of $\lambda_1(\rho_{\max})$.
  \label{fig-rhoD-nonstar}
\end{figure}
\begin{figure}[H]
  \begin{center}
    \begin{tabular}{llll}
      {\bfseries (a)} 
      &
      {\bfseries (b)} 
      & 
      {\bfseries (c)}
      & 
      {\bfseries (d)}
      \cr
      \includegraphics[width=3.0cm]{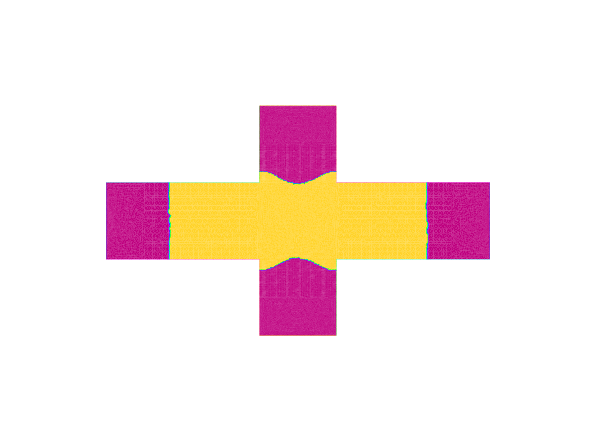}
      &
      \includegraphics[width=3.0cm]{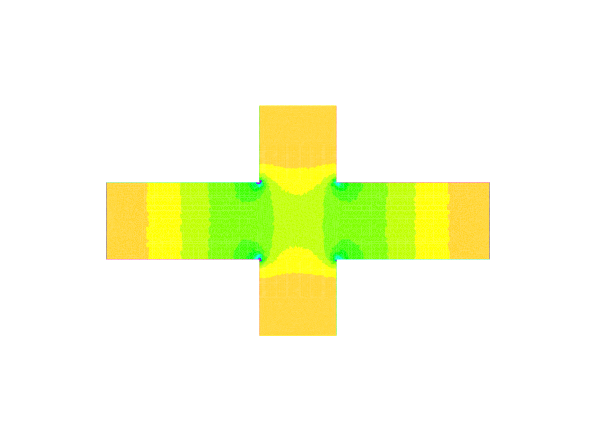}
      &
      \includegraphics[width=3.0cm]{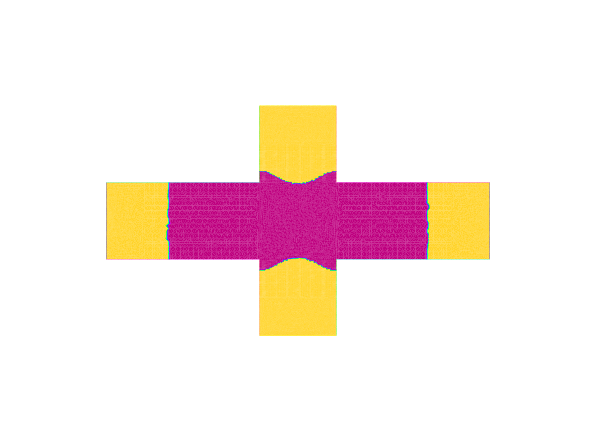}
      &
      \includegraphics[width=3.0cm]{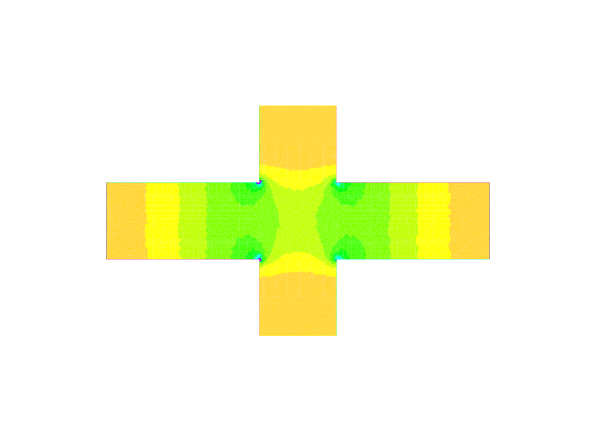}
      \cr
    \end{tabular}
  \end{center}
  \caption{Optimizers and eigenfunctions for Problem \ref{prob-rho} on $\Omega_6$}
  (a) minimizer $\rho_{\min}$, 
  (b) $|\rho \nabla u|^2$ of the associated eigenfunction of $\lambda_1(\rho_{\min})$.
  (c) maximizer $\rho_{\max}$, 
  (d) $|\rho \nabla u|^2$ of the associated eigenfunction of $\lambda_1(\rho_{\max})$.
  \label{fig-rhoN-cross}
\end{figure}

Next, we consider two-dimensional problems.
In the following numerical experiments, we fix $c=1.1$
and the constant of volume constraint $m_0= 0.5$ in (\ref{constraint}) unless otherwise noted and consider the following $\Omega$: 
\begin{align*}
  &\Omega_1 := (-\pi,\pi)\times (-\pi,\pi)\subset \mathbb{R}^2,\quad \text{(square)}\\
  &\Omega_2 := 
  \scalebox{0.8}[1.0]{$
    (0,3)\times (0,3)\setminus 
    \{((0,1]\times (0,1]) \cup ((0,1]\times [2,3)) \cup ([2,3)\times (0,1]) \cup ([2,3)\times [2,3))\}
    $}
  \subset \mathbb{R}^2,\\
  &\quad \quad \text{(cross, e.g.\ Figure \ref{fig-rhoD-cross})}\\
  &\Omega_3 := \{ x^2+y^2 < 1\}\subset \mathbb{R}^2,\quad \text{(disk)}\\
  &\Omega_4 := \{x^2+y^2 < 1\} \setminus  [\{(x-0.5)^2+y^2 < 0.04\}\cup \{(x+0.5)^2+y^2 < 0.04\}]\subset \mathbb{R}^2,\\
  &\quad \quad \quad \text{(disk with two holes, e.g.\ Figure \ref{fig-rhoD-nonstar})}\\
  &\Omega_5 := (-\pi, \pi)\times (-2\pi, 2\pi)\subset \mathbb{R}^2,\quad \text{(rectangle)}\\
  &\Omega_6 := 
  \scalebox{0.8}[1.0]{$
    (0,5)\times (0,3)\setminus \{((0,2]\times (0,1]) \cup ((0,2]\times [2,3)) \cup ([3,5)\times (0,1]) \cup ([3,5)\times [2,3))\}
    $}
  \subset \mathbb{R}^2,\\
  &\quad \quad \quad \text{(rectangular cross, e.g.\ Figure \ref{fig-rhoN-cross})}\\
  &\Omega_7 := \{x^2+y^2/4 < 1\}\subset \mathbb{R}^2\quad \text{(ellipse)}\quad \text{ and }\\
  &\Omega_8 := 
  \scalebox{0.9}[1.0]{$
    \{ x^2+y^2/4 < 1\} \setminus  [\{(x-0.5)^2+y^2/16 < 0.04\}\cup \{(x+0.5)^2+y^2/16 < 0.04\}]
    $}
  \subset \mathbb{R}^2.\\
  &\quad \quad \quad \text{(ellipse with two holes, e.g.\ Figure \ref{fig-rhoN-nonstar})}
\end{align*}
In two-dimensional problems,
we apply the FreeFEM++ library \cite{FreeFEM++} to all the computations.
\par
\begin{figure}[H]
  \begin{center}
    \begin{tabular}{llll}
      {\bfseries (a)} 
      &
      {\bfseries (b)} 
      & 
      {\bfseries (c)}
      & 
      {\bfseries (d)}
      \cr
      \includegraphics[width=3.0cm]{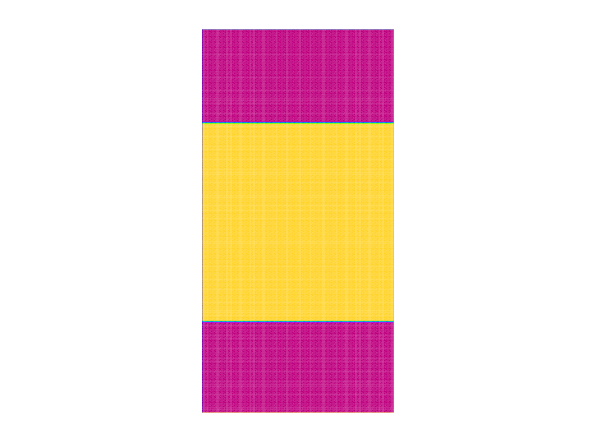}
      &
      \includegraphics[width=3.0cm]{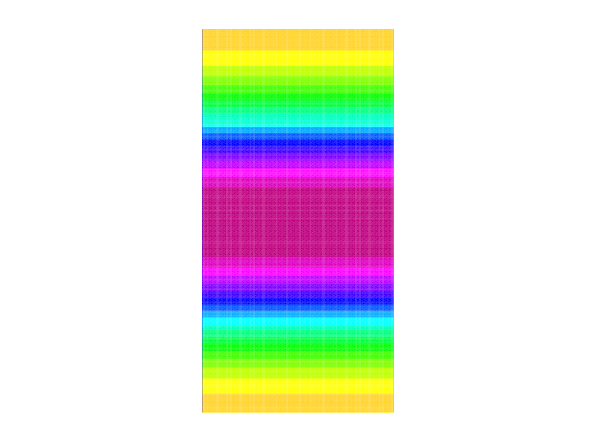}
      &
      \includegraphics[width=3.0cm]{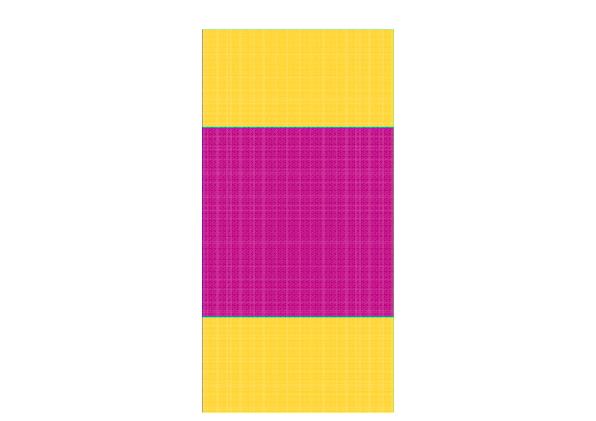}
      &
      \includegraphics[width=3.0cm]{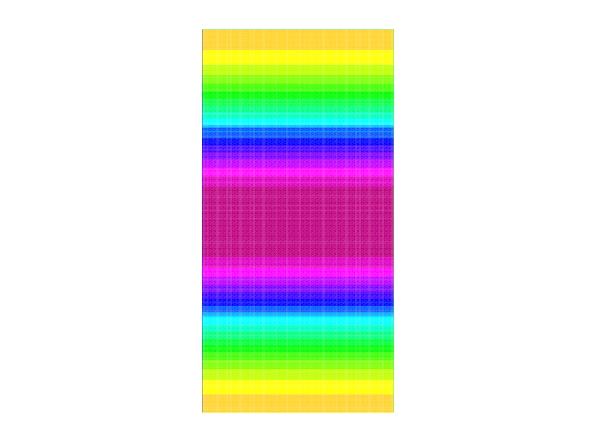}
      \cr
    \end{tabular}
  \end{center}
  \caption{Optimizers and eigenfunctions for Problem \ref{prob-rho} on $\Omega_5$}
  (a) minimizer $\rho_{\min}$, 
  (b) $|\rho \nabla u|^2$ of the associated eigenfunction of $\lambda_1(\rho_{\min})$.
  (c) maximizer $\rho_{\max}$, 
  (d) $|\rho \nabla u|^2$ of the associated eigenfunction of $\lambda_1(\rho_{\max})$.
  \label{fig-rhoN-rect}
\end{figure}
\begin{figure}[H]
  \begin{center}
    \begin{tabular}{llll}
      {\bfseries (a)} 
      &
      {\bfseries (b)} 
      & 
      {\bfseries (c)}
      & 
      {\bfseries (d)}
      \cr
      \includegraphics[width=3.0cm]{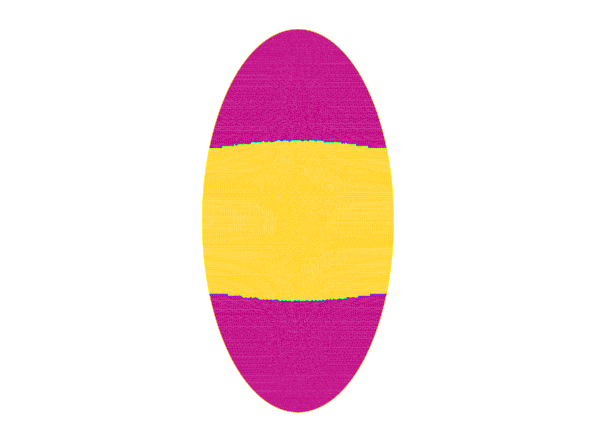}
      &
      \includegraphics[width=3.0cm]{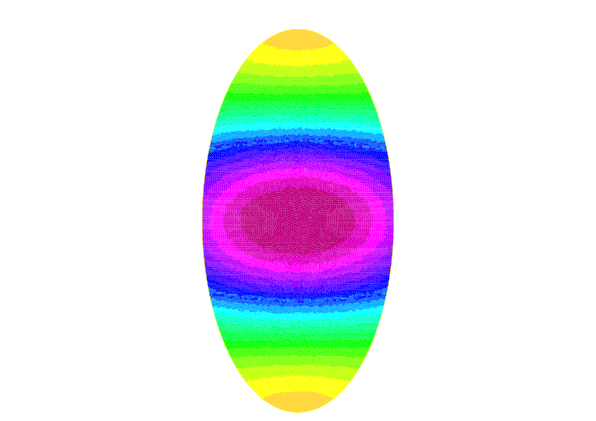}
      &
      \includegraphics[width=3.0cm]{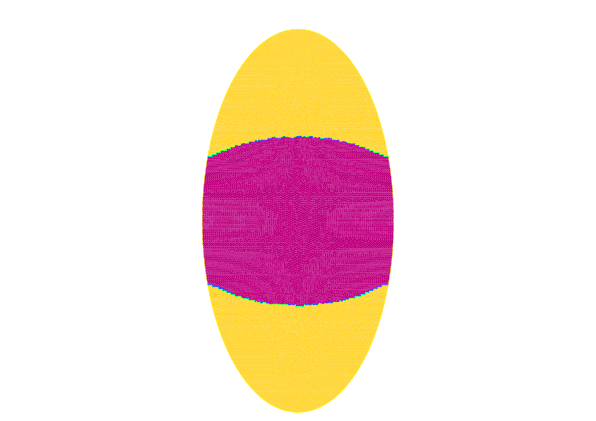}
      &
      \includegraphics[width=3.0cm]{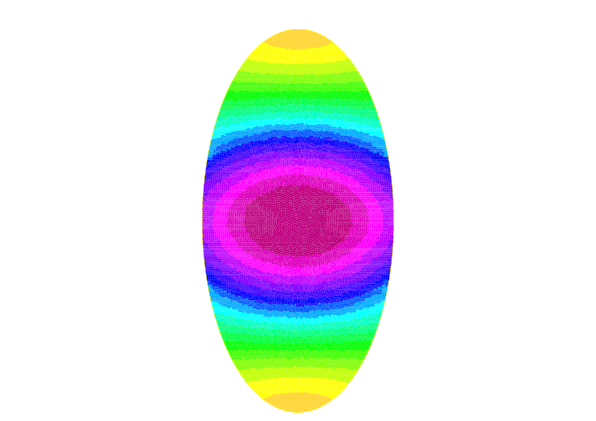}
      \cr
    \end{tabular}
  \end{center}
  \caption{Optimizers and eigenfunctions for Problem \ref{prob-rho} on $\Omega_7$}
  (a) minimizer $\rho_{\min}$, 
  (b) $|\rho \nabla u|^2$ of the associated eigenfunction of $\lambda_1(\rho_{\min})$.
  (c) maximizer $\rho_{\max}$, 
  (d) $|\rho \nabla u|^2$ of the associated eigenfunction of $\lambda_1(\rho_{\max})$.
  \label{fig-rhoN-ellipse}
\end{figure}
\begin{figure}[H]
  \begin{center}
    \begin{tabular}{llll}
      {\bfseries (a)} 
      &
      {\bfseries (b)} 
      & 
      {\bfseries (c)}
      & 
      {\bfseries (d)}
      \cr
      \includegraphics[width=3.0cm]{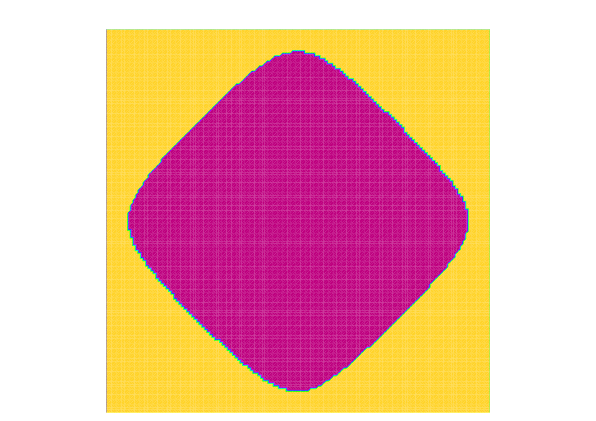}
      &
      \includegraphics[width=3.0cm]{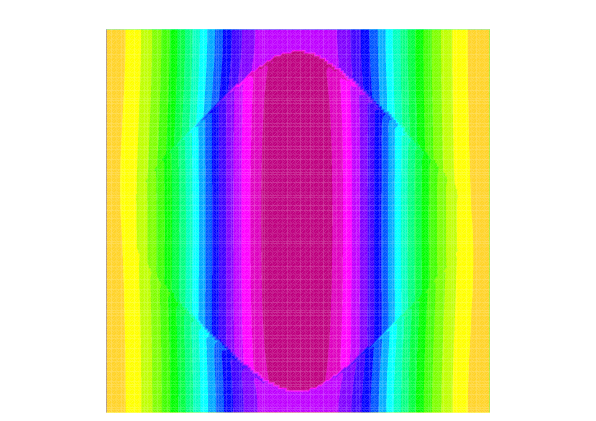}
      &
      \includegraphics[width=3.0cm]{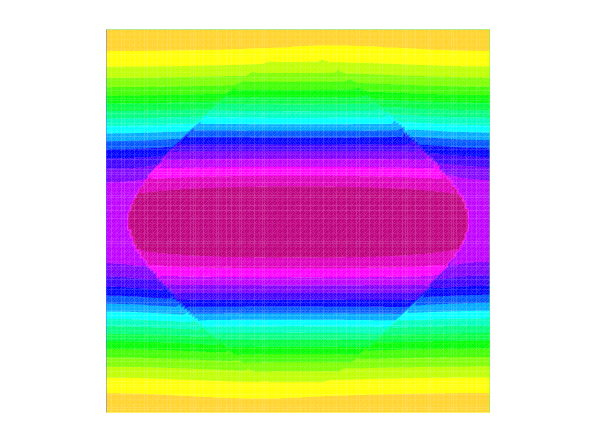}
      &
      \includegraphics[width=3.0cm]{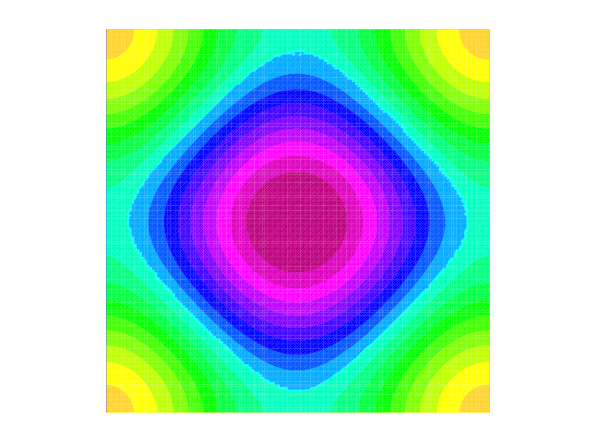}
      \cr
    \end{tabular}
  \end{center}
  \caption{Maximizers and eigenfunctions for Problem \ref{prob-rho} on $\Omega_1$}
  Figure (a) shows the maximizer $\rho_{\max}$ of $\lambda_1(\rho)$. 
  Figure (b) and (c) show corresponding $|\rho_{max} \nabla u_{1,\max}|^2$ and $|\rho_{\max} \nabla u_{2,\max}|^2$, 
  where $u_{1,\max}$ and $u_{2,\max}$ are associated eigenfunctions of $\lambda_1(\rho_{\max})$ and $\lambda_2(\rho_{\max})$ (actually equal to $\lambda_1(\rho_{\max})$), respectively, 
  after the normalization so that $\int_\Omega |u_{i,\max}|^2 = 1$ holds.  
  Figure (d) shows $|\rho_{\max} \nabla u_{1,\max}|^2+|\rho_{\max} \nabla u_{2,\max}|^2$ after normalizations. 
  \label{fig-rhoN-square}
\end{figure}

Focus on Figure \ref{fig-rhoD-square}, \ref{fig-rhoD-cross}, 
\ref{fig-rhoD-disk}, \ref{fig-rhoD-nonstar}, \ref{fig-rhoN-rect}, 
\ref{fig-rhoN-cross}, \ref{fig-rhoN-ellipse} and \ref{fig-rhoN-nonstar}.
In these cases, 
$S$ for optimizers corresponds to the super- or the sub-level set of 
$|\rho_\ast \nabla u_\ast|$ regardless of topologies of $\Omega$ and boundary conditions as in one-dimensional problems, 
where $u_\ast$ is the associated eigenfunction of $\lambda_1(\rho_\ast)$. 
More precisely, 
if $\rho_\ast\in \mathcal{K}$ minimizes or maximizes $\lambda_1$ and if $\lambda_1$ is {\emx simple} (see Table \ref{table-simple-rho}), 
then the corresponding eigenfunction $u_1$ satisfies
\begin{align}
  \label{inequality-min-rho}
  \text{Minimization : }\quad &|\rho_\ast(x) \nabla u_1(x)| \leq |\rho_\ast(y) \nabla u_1(y)|\quad \text{a.e.\ }x\in S\text{ and }y\in S^c.\\
  \label{inequality-max-rho}
  \text{Maximization : }\quad &|\rho_\ast(y) \nabla u_1(y)| \leq |\rho_\ast(x) \nabla u_1(x)|\quad \text{a.e.\ }x\in S\text{ and }y\in S^c.
\end{align}
Inequalities (\ref{inequality-min-rho}) and (\ref{inequality-max-rho}) can be considered as analogue of those in Theorem \ref{preceding-sigma}. 
If (\ref{inequality-min-rho}) and (\ref{inequality-max-rho}) hold, 
then the simpleness of $\lambda_1(\rho_\ast)$ holds and vice versa (Table \ref{table-simple-rho}). 
Indeed, 
{\emx
  if $\lambda_1$ is not simple,  
  then the inequality (\ref{inequality-max-rho}) does not hold%
}, 
which will be seen from Figure \ref{fig-rhoN-square} and Table \ref{table-simple-rho}.
\begin{rem}
  The maximum principle guarantees the simpleness of the smallest eigenvalue for (uniformly) elliptic operators.
  In the case of (\ref{ev-prob-rho}) with 
  the homogeneous Dirichlet boundary condition, 
  the smallest eigenvalue $\lambda_1$ is positive.
  The smallest eigenvalue $\lambda_1(\rho)$ is then always simple for any 
  $\rho \in \mathcal{K}$ in the case of the homogenous Dirichlet boundary condition.
  However, under the homogeneous Neumann boundary condition, 
  $\lambda_1(\rho)$ is the second smallest eigenvalue, 
  since the smallest eigenvalue is $0$, 
  and the simpleness of $\lambda_1(\rho)$ is not guaranteed. 
\end{rem}
\par
Assume that $\rho_{\max}$ attains the maximum of $\lambda_1$ on $\mathcal{K}$ and that $\lambda_1(\rho_{\max})$ has multiplicity two as shown in Figure \ref{fig-rhoN-square}.
Then the energy functional in (\ref{Lagrange}) with multiplicity two eigenvalue $\lambda_1(\rho_{\max})=\lambda_2(\rho_{\max})$ can be also written by
\begin{align*}
  L(\phi) &= F(\phi) + \nu_1 G_1(\phi) + \nu_2 G_2(\phi),\\
  F(\phi) &= -\lambda_1(\phi) = -(a_1\lambda_1(\phi) + a_2\lambda_2(\phi)),\quad a_1, a_2 \geq 0\text{ with }a_1+a_2 = 1,\\
  G_1(\phi) &= \int_{\{\phi(x) > 0\}}dx - m_0|\Omega|,\quad G_2(\phi) =\lambda_2(\phi)-\lambda_1(\phi) ,\quad \nu_1, \nu_2\in \mathbb{R}.
\end{align*}
Corresponding $v_0(x)$ in (\ref{v0}) 
for obtaining the steepest descent flow is
\begin{equation*}
  v_0(x) 
  = 
  -\left(\frac{a_1 (c-1)}{\int_\Omega u_{1,\phi}^2dx } 
    |\nabla u_{1,\phi}(x)|^2 + \frac{a_2 (c-1)}{\int_\Omega u_{2,\phi}^2dx } |\nabla u_{2,\phi}(x)|^2\right).
\end{equation*}
\par
Since $\lambda_1(\phi) = \lambda_2(\phi)$ holds for optimizers,
then one can see that the optimizer should satisfy 
$a_1' |\nabla u_1(x)|^2 + a_2' |\nabla u_2(x)|^2 \equiv \text{constant}$ on $\partial S$,
where $a_i' = a_i(\int_\Omega u_{i,\phi}^2dx)^{-1}$.
By the definition of $S$ by $\phi$, 
$\rho = 1$ holds on $\partial S$ and 
hence the above equality on $\partial S$ is equivalent to 
\begin{equation}
  \label{inequality-max-rho-s}
  \rho(x)^s(a_1' |\nabla u_1(x)|^2 + a_2' |\nabla u_2(x)|^2) 
  \equiv 
  \text{constant}\quad \text{a.e.\ on }\partial S \quad \text{ for some }s>0.
\end{equation}
Moreover, 
the constraints $G_1(\phi) = 0$ and $G_2(\phi)=0$ must be kept during evolution of $\phi$ by (\ref{SDF}).
Optimizers $\phi$ then have to satisfy
\begin{equation*}
  D_\phi G_i(\phi) \delta \phi = 0,\quad i=1,2
\end{equation*}
for the first variation $\delta \phi$ of $\phi$.
In particular, it follows that
\begin{equation*}
  \int_{\partial S} \rho(x)^s
  \left( 
    \frac{a_1|\nabla u_{1,\phi}(x)|^2}{\int_\Omega u_{1,\phi}^2dx }  - \frac{a_2|\nabla u_{2,\phi}(x)|^2}{\int_\Omega \sigma(\phi) u_{2,\phi}^2dx }
  \right)
  ds(x) = 0
\end{equation*}
should be satisfied by calculations discussed in \cite{OS}.
In particular, positive constants $a_1'$ and $a_2'$ must be identical 
and hence $a_1$, $a_2 > 0$. 
\par
Since $a_1$ and $a_2$ are arbitrary, 
we may choose $a_1 = a_2 = 1/2$. 
If we normalize eigenfunctions $u_1$ and $u_2$ 
so that $\int_\Omega |u_i|^2 = 1$ ($i=1,2$), 
then the following inequality will be the maximization criterion 
for $\lambda_1(\rho)$ with multiplicity two:
\begin{equation}
  \label{inequality-max-rho-2}
  \rho(y)^2 (|\nabla u_1(y)|^2+|\nabla u_2(y)|^2) 
  \leq 
  \rho(x)^2 (|\nabla u_1(x)|^2 + |\nabla u_2(x)|^2),
  \ \text{ a.e.\ }x\in S\text{ and }y\in S^c.
\end{equation}
Here we chose $s=2$ in (\ref{inequality-max-rho-s}),
which is natural because $\rho$ and $\nabla u$ 
have the same order in (\ref{inequality-min-rho}) and (\ref{inequality-max-rho}).
We can see that (\ref{inequality-max-rho-2}) is actually satisfied 
(see Figure \ref{fig-rhoN-square}).
\par
We conclude our numerical observations for optimization criteria in Problem \ref{prob-rho}.
\begin{observation}
\label{observation-criteria-heat}
  Consider Problem \ref{prob-rho} with bounded a domain $\Omega \subset \mathbb{R}^n$.
  If $\rho_{\min}\in \mathcal{K}$ is the minimizer of $\lambda_1$ in $\mathcal{K}$,
  then the eigenfunction $u_1$ associated with $\lambda_1$ satisfies (\ref{inequality-min-rho}).
  Similarly, 
  if $\rho_{\max}\in \mathcal{K}$ is the maximizer of $\mu_1$ in $\mathcal{K}$ and if $\lambda_1(\rho_{\max})$ is simple,
  then the eigenfunction $u_1$ associated with $\lambda_1(\rho_{\max})$ satisfies (\ref{inequality-max-rho}).
  If $\rho_{\max}\in \mathcal{K}$ is the maximizer of $\lambda_1$ in $\mathcal{K}$ and if $\lambda_1(\rho_{\max})$ has multiplicity two,
  then the eigenfunction $u_1$ associated with $\lambda_1(\rho_{\max})$ satisfies (\ref{inequality-max-rho-2}).
\end{observation}
These optimization criteria will be generalized to higher dimensional problems, in particular, the case that $\lambda_1(\rho)$ has multiplicity $k \geq 3$ in the similar manner.

\subsection{Geometry of optimizers}
\label{section-symmetry}
Next, we consider geometry of the domain 
$S_\ast = \{x\in \Omega \mid \rho_\ast (x) = c\}$ 
defined by the optimizer $\rho_\ast$. 
Inequalities (\ref{inequality-min-rho}), (\ref{inequality-max-rho})  
and (\ref{inequality-max-rho-2}) imply that 
various properties of $S_\ast$ for optimizers come from corresponding eigenfunctions.
It is natural to expect that
{\emx
  some geometric properties inherit from $\Omega$.
}
Here, we focus on connectivity, convexity, star-shapedness and symmetry.
\par
\begin{figure}[H]
  \begin{center}
    \begin{tabular}{llll}
      {\bfseries (a)} 
      &
      {\bfseries (b)} 
      & 
      {\bfseries (c)}
      & 
      {\bfseries (d)}
      \cr
      \includegraphics[width=3.0cm]{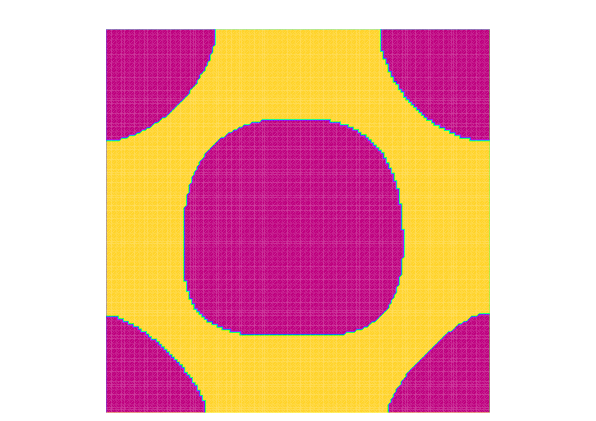}
      &
      \includegraphics[width=3.0cm]{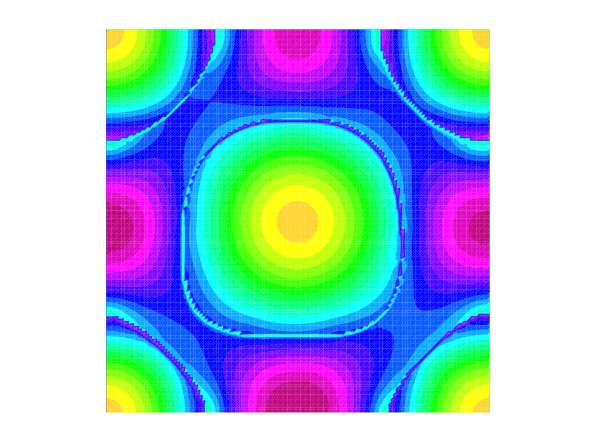}
      &
      \includegraphics[width=3.0cm]{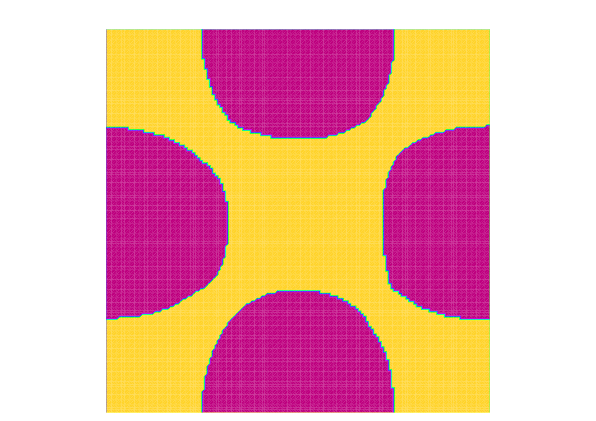}
      &
      \includegraphics[width=3.0cm]{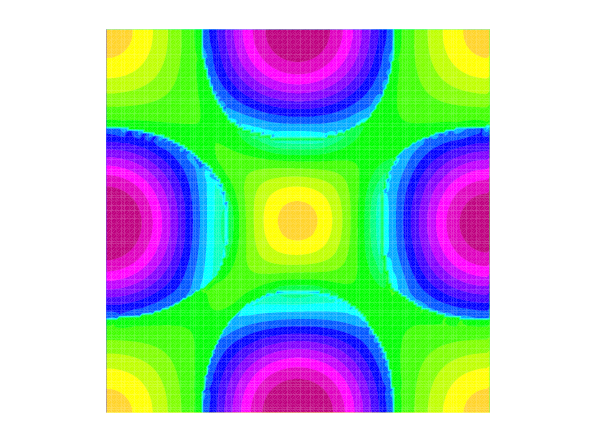}
      \cr
    \end{tabular}
  \end{center}
  \caption{Optimizers and eigenfunctions for Problem \ref{prob-rho} on $\Omega_1$.}
  (a) minimizer $\rho_{\min}$, 
  (b) $|\rho \nabla u|^2$ of the associated eigenfunction of $\lambda_1(\rho_{\min})$.
  (c) maximizer $\rho_{\max}$, 
  (d) $|\rho \nabla u|^2$ of the associated eigenfunction of $\lambda_1(\rho_{\max})$.
  \label{fig-rhoD-square}
\end{figure}
\begin{figure}[H]
  \begin{center}
    \begin{tabular}{llll}
      {\bfseries (a)} 
      &
      {\bfseries (b)} 
      & 
      {\bfseries (c)}
      & 
      {\bfseries (d)}
      \cr
      \includegraphics[width=3.0cm]{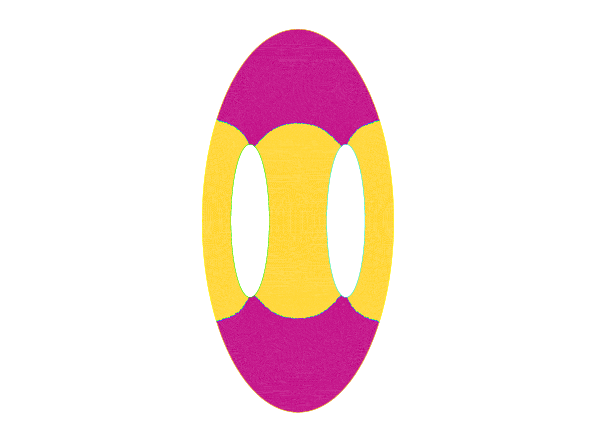}
      &
      \includegraphics[width=3.0cm]{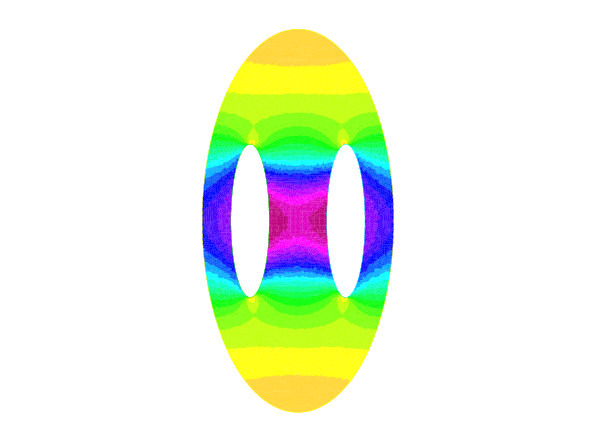}
      &
      \includegraphics[width=3.0cm]{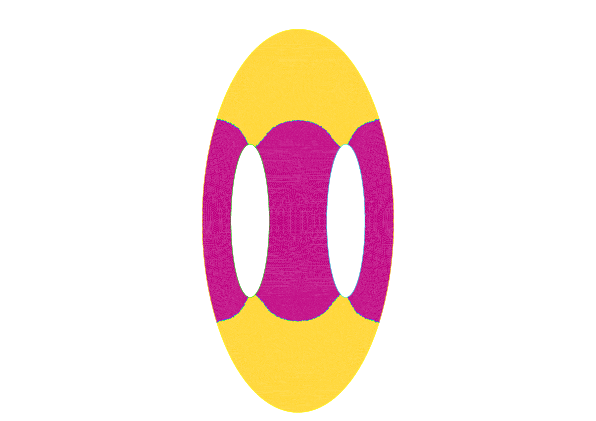}
      &
      \includegraphics[width=3.0cm]{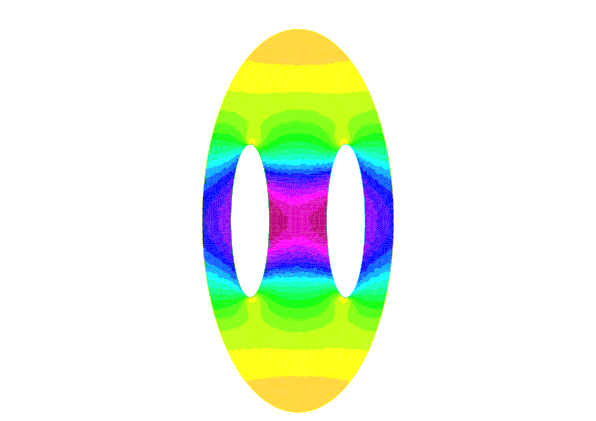}
      \cr
    \end{tabular}
  \end{center}
  \caption{Optimizers and eigenfunctions for Problem \ref{prob-rho} on $\Omega_8$}
  (a) minimizer $\rho_{\min}$, 
  (b) $|\rho \nabla u|^2$ of the associated eigenfunction of $\lambda_1(\rho_{\min})$.
  (c) maximizer $\rho_{\max}$, 
  (d) $|\rho \nabla u|^2$ of the associated eigenfunction of $\lambda_1(\rho_{\max})$.
  \label{fig-rhoN-nonstar}
\end{figure}

In Figure \ref{fig-rhoD-square} 
(square with the homogeneous Dirichlet boundary),
neither $S_\ast$ nor $S_\ast^c$ are even connected, 
even if $\Omega$ is convex.
In the case of the homogeneous Neumann boundary value problem,
either $S_\ast$ or $S_\ast^c$ is convex if $\Omega$ is convex,
according to our computation results. On the other hand,
in the case of the non-simply connected domain $\Omega_8$ (Figure \ref{fig-rhoN-nonstar}),
$S_\ast$ maximizing $\lambda_1(\rho)$ is not connected,
although $\Omega_8$ is connected.
However, by the inequality (\ref{inequality-max-rho}),
one can easily confirm that we can choose $m_0 \in (0.5,1)$ so that $S_\ast$ maximizing $\lambda_1(\rho)$ is connected. 
As a consequence, 
there is generally no topological correspondence
between $\Omega$ and $S_\ast$ (or $S_\ast^c$)
which are independent of boundary conditions on $\partial \Omega$ or $m_0$.
\par
\begin{figure}[H]
  \begin{center}
    \begin{tabular}{llll}
      {\bfseries (a)} 
      &
      {\bfseries (b)} 
      & 
      {\bfseries (c)}
      & 
      {\bfseries (d)}
      \cr
      \includegraphics[width=3.0cm]{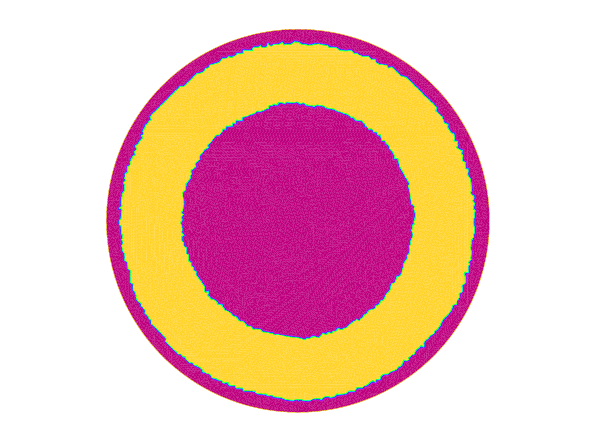}
      &
      \includegraphics[width=3.0cm]{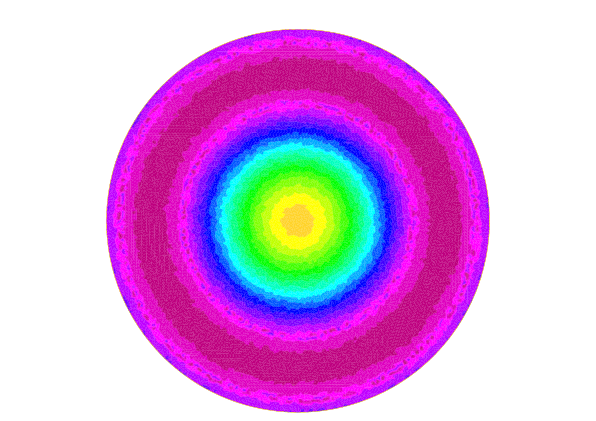}
      &
      \includegraphics[width=3.0cm]{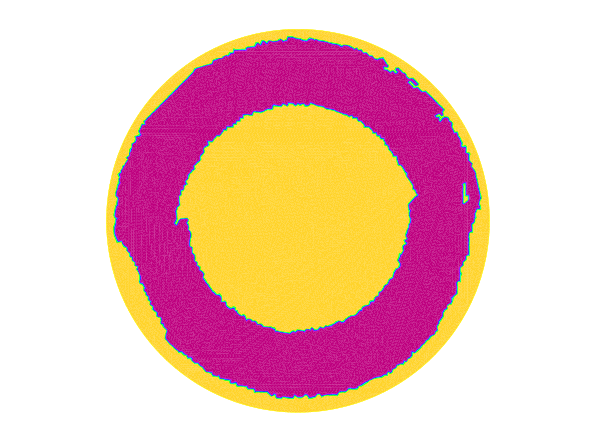}
      &
      \includegraphics[width=3.0cm]{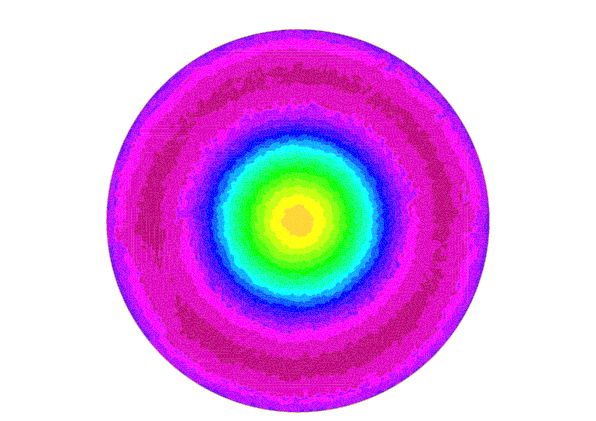}
      \cr
    \end{tabular}
  \end{center}
  \caption{Optimizers and eigenfunctions for Problem \ref{prob-rho} on $\Omega_3$}
  (a) minimizer $\rho_{\min}$, 
  (b) $|\rho \nabla u|^2$ of the associated eigenfunction of $\lambda_1(\rho_{\min})$.
  (c) maximizer $\rho_{\max}$, 
  (d) $|\rho \nabla u|^2$ of the associated eigenfunction of $\lambda_1(\rho_{\max})$.
  \label{fig-rhoD-disk}
\end{figure}

Next we focus on symmetry of $S_\ast$.
Inequalities (\ref{inequality-min-rho}) and (\ref{inequality-max-rho}) 
imply that symmetry of $S_\ast$ comes from that of $\rho \nabla u$ given by corresponding eigenfunction.
If $\rho \nabla u$ is symmetric in a certain axis direction or in rotation,
then, thanks to the original equation 
$-\nabla \cdot (\rho \nabla u) = \lambda u$,
$u$ and $\nabla u$ will be also symmetric.
Finally symmetry of $\rho$ will hold from symmetry of $\rho \nabla u$.
The key consideration is that {\emx whether the symmetry of $\rho_\ast \nabla u_\ast$ associated with the optimizer $\rho_\ast$ inherits from $\Omega$}.
Our numerical simulations argue that symmetry of 
$\rho_\ast \nabla u_\ast$ inherits from $\Omega$
(Figure \ref{fig-rhoD-square}, 
\ref{fig-rhoD-cross},
\ref{fig-rhoD-disk},
\ref{fig-rhoD-nonstar},
\ref{fig-rhoN-rect},
\ref{fig-rhoN-cross},
\ref{fig-rhoN-ellipse} and \ref{fig-rhoN-nonstar}).
Since $u_\ast$ is the solution of (\ref{ev-prob-rho}), 
symmetry of $\rho_\ast \nabla u_\ast$ leads to that of $u_\ast$ and $\rho_\ast$.
One then observes the following.
\begin{observation}
\label{observation-symmetry-heat}
  For Problem \ref{prob-rho}, 
  let $\rho_\ast$ be the optimizer of $\lambda_1(\rho)$,
  $u_\ast$ be the associated eigenfunctions of $\lambda_1(\rho_\ast)$ and $S_\ast = \{x\in \Omega \mid \rho_\ast(x) = c\}$.
  If $\Omega$ is star-shaped and symmetric in a certain direction or in rotation,
  then 
  $S_\ast$ and $S_\ast^c$ are also symmetric in the direction
  regardless of boundary conditions on $\partial \Omega$.
\end{observation}
Our examples also show that there is a possibility that 
both $S_\ast$ and $S_\ast^c$ are symmetric 
even if $\Omega$ is not star-shaped (see Figures \ref{fig-rhoD-nonstar} and \ref{fig-rhoN-nonstar}).
However, in the Dirichlet boundary case, there is also a possibility that symmetry breaking of $S_\ast$ occurs in the case of $\Omega$ being either an annulus or a dumbbell,
which are not star-shaped. This is indeed the case of the eigenvalue problem $-\Delta u + (c-1) \chi_{S} u = \lambda u$ as mentioned in Section \ref{section-prior} and 
the same phenomenon may occur in Problem \ref{prob-rho}.
\par
On the minimization problem for the disk with the homogeneous Dirichlet boundary condition,
we also mention that the optimal configuration $S_\ast$ is a union of concentric annuli and a disk (See Figure \ref{fig-rhoD-disk}).
This observation implies that the result of \cite{Conca-Dambrine-Mahadevan-Quintero}
mentioned in Section \ref{section-prior} may also hold if c is not close to 1.

\subsection{Continuous dependence of optimizers on boundary conditions}
\label{section-cont}
\begin{figure}[H]
  \begin{center}
    \begin{tabular}{cccccc}
      (a) 
      &
      \includegraphics[width=2.25cm]{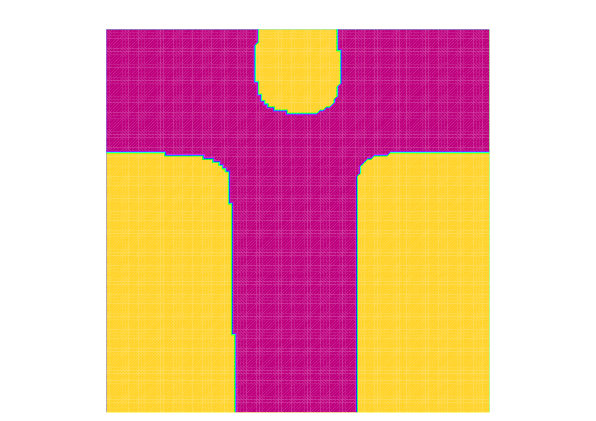}
      &
      \includegraphics[width=2.25cm]{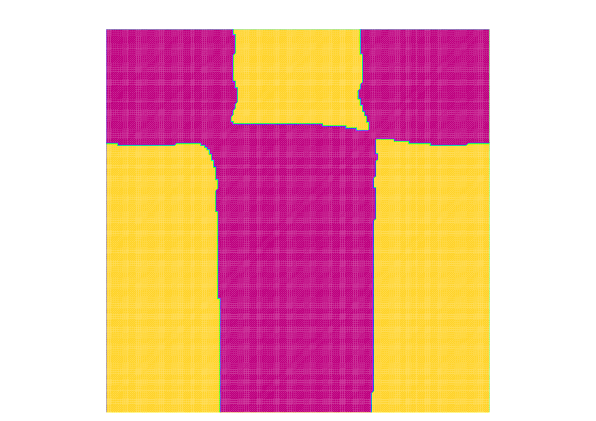}
      &
      \includegraphics[width=2.25cm]{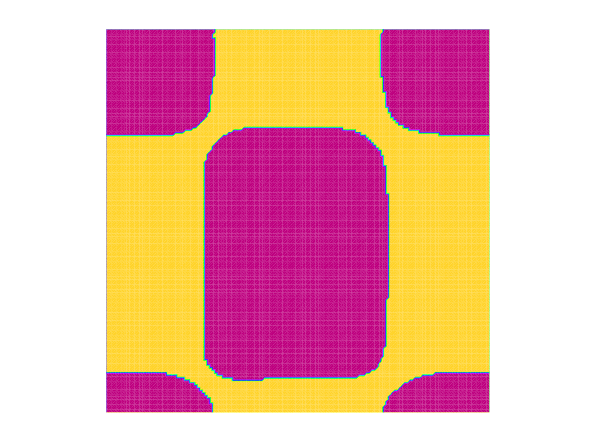}
      &
      \includegraphics[width=2.25cm]{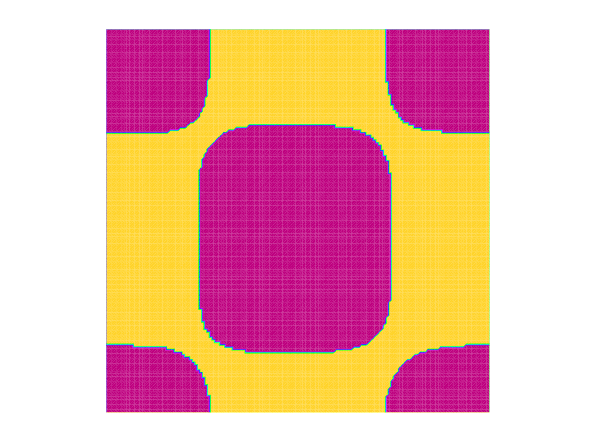}
      &
      \includegraphics[width=2.25cm]{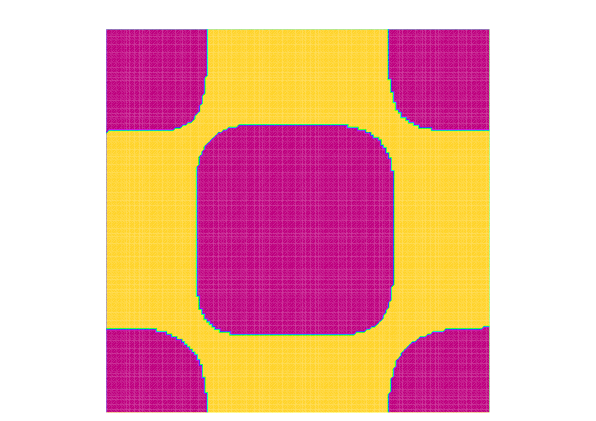}
      \cr
      (b) 
      &
      \includegraphics[width=2.25cm]{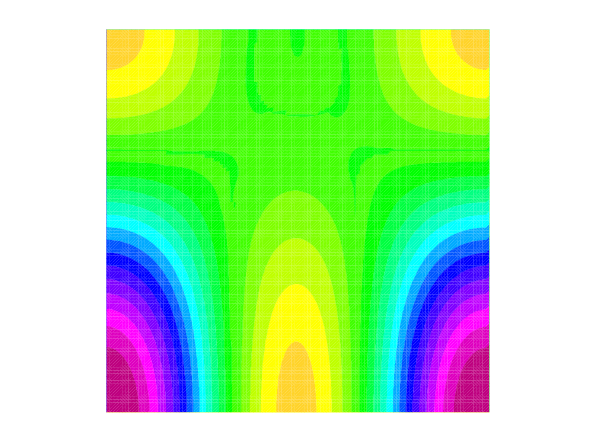}
      &
      \includegraphics[width=2.25cm]{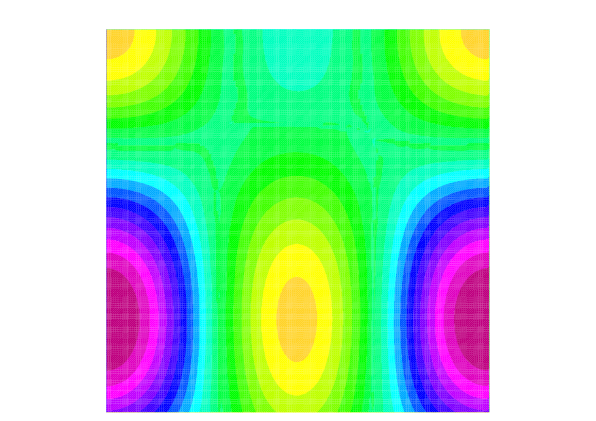}
      &
      \includegraphics[width=2.25cm]{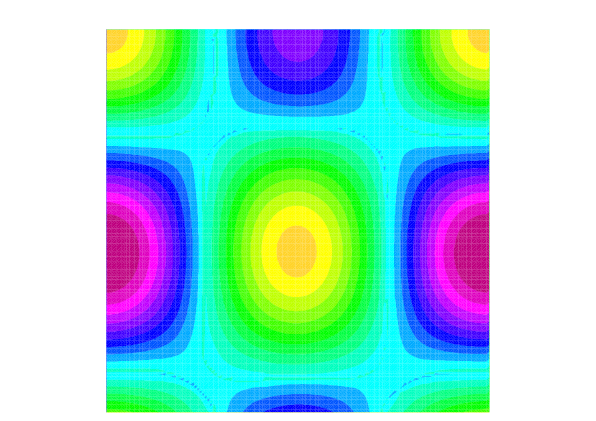}
      &
      \includegraphics[width=2.25cm]{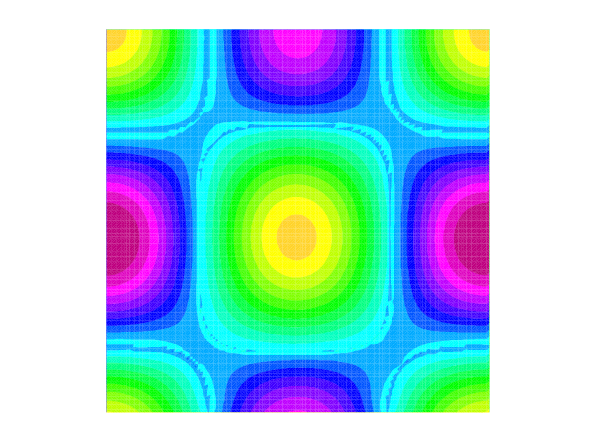}
      &
      \includegraphics[width=2.25cm]{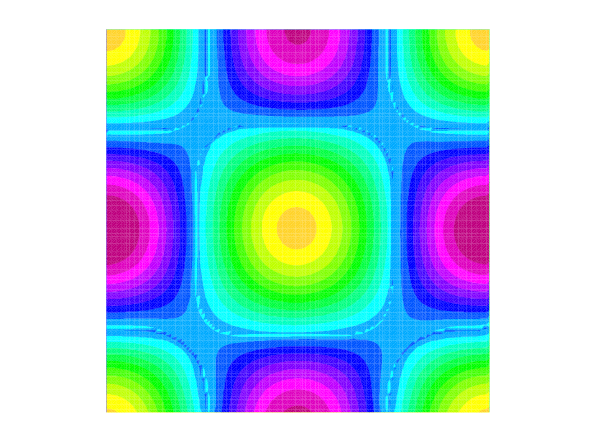}
      \cr
      &
      $\eta=0$
      &
      $\eta=0.2$
      &
      $\eta=0.5$
      &
      $\eta=1.0$
      &
      $\eta=5.0$
      \cr
    \end{tabular}
  \end{center}
  \caption{Minimizer for Problem \ref{prob-rho} on $\Omega_1$}
  (a) shows the minimizer $\rho_{\min}$ and 
  (b) shows $|\rho_{\min} \nabla u_{\min}|^2$ for the eigenfunction $u_{\min}$ of $\rho_1(\lambda_{\min})$ with various $\eta$
  The minimization criterion (\ref{inequality-min-rho}) is actually satisfied in each case. 
  The larger $\eta$ becomes, the closer $\rho_{\min}$ is to the minimizer with homogeneous Dirichlet boundary condition (Figure \ref{fig-rhoD-square}, a).
  \label{fig-rho-cont-min}
\end{figure}

\begin{figure}[H]
  \begin{center}
    \begin{tabular}{cccccc}
      (a) 
      &
      \includegraphics[width=2.25cm]{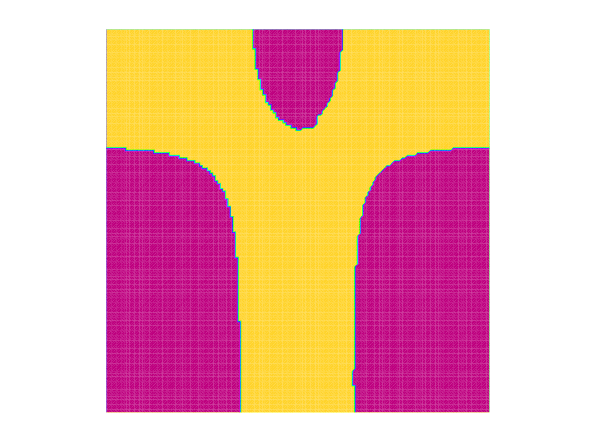}
      &
      \includegraphics[width=2.25cm]{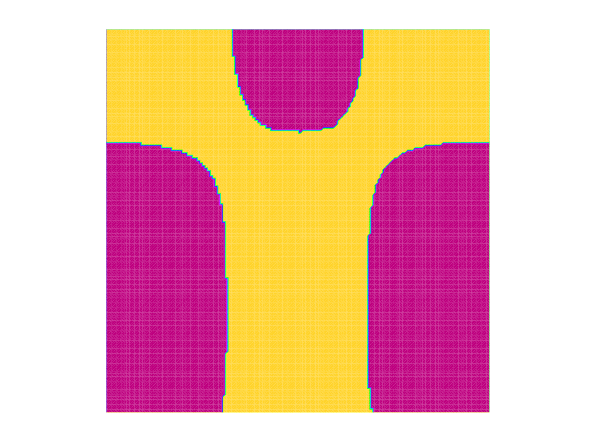}
      &
      \includegraphics[width=2.25cm]{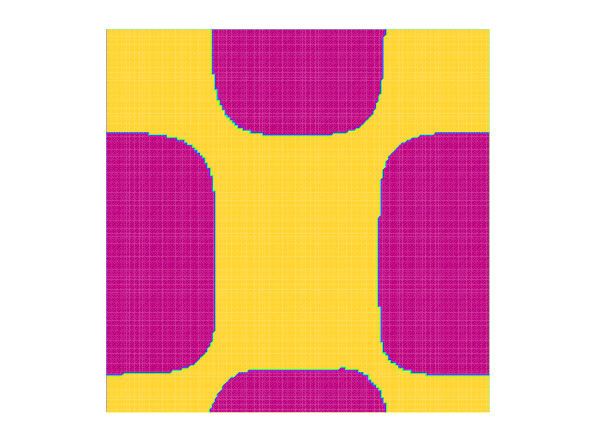}
      &
      \includegraphics[width=2.25cm]{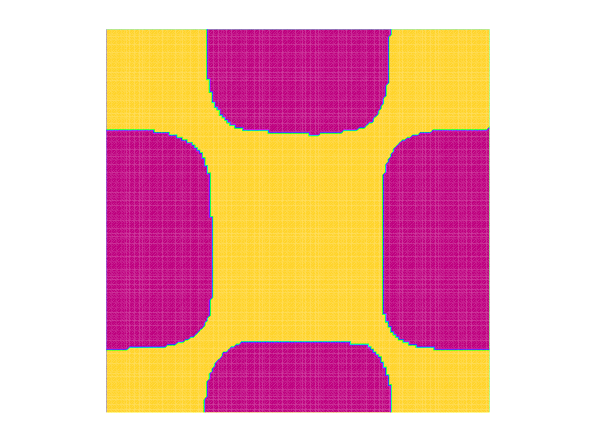}
      &
      \includegraphics[width=2.25cm]{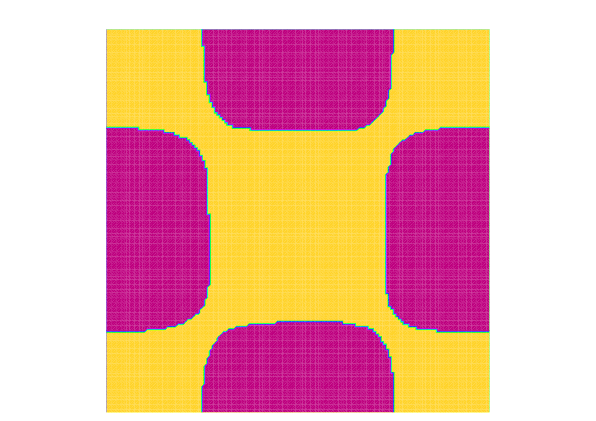}
      \cr
      (b) 
      &
      \includegraphics[width=2.25cm]{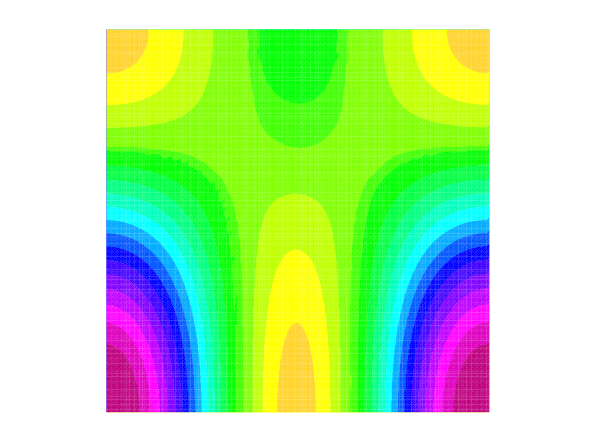}
      &
      \includegraphics[width=2.25cm]{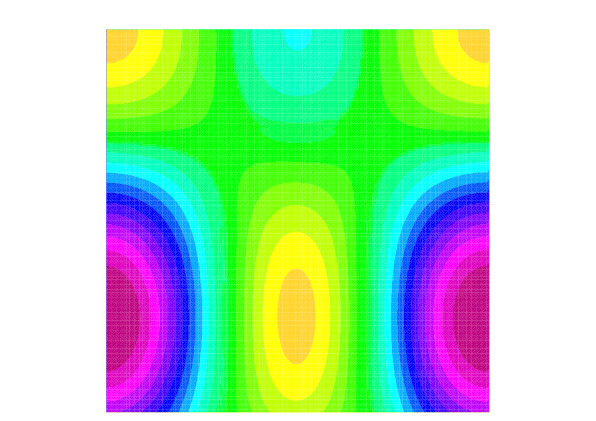}
      &
      \includegraphics[width=2.25cm]{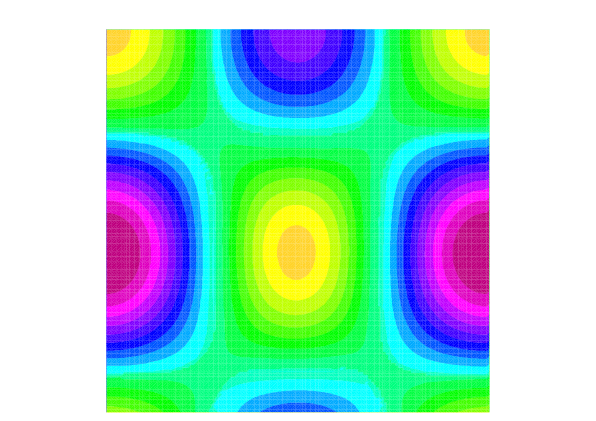}
      &
      \includegraphics[width=2.25cm]{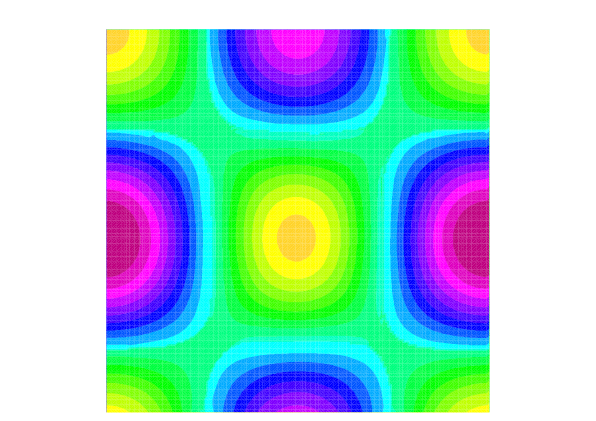}
      &
      \includegraphics[width=2.25cm]{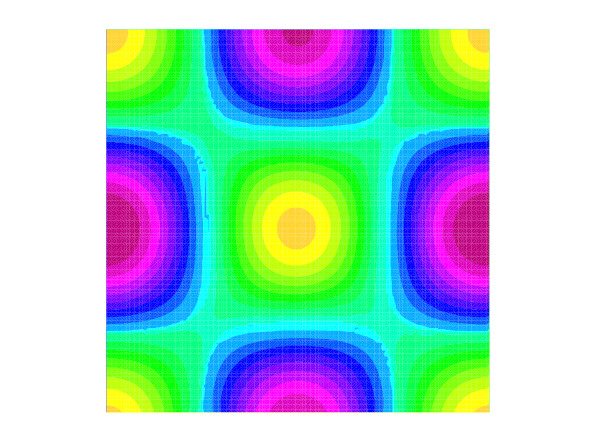}
      \cr
      &
      $\eta=0$
      &
      $\eta=0.2$
      &
      $\eta=0.5$
      &
      $\eta=1.0$
      &
      $\eta=5.0$
      \cr
    \end{tabular}
  \end{center}
  \caption{Maximizer for Problem \ref{prob-rho} on $\Omega_1$}
  (a) shows the maximizer $\rho_{\max}$ and 
  (b) shows $|\rho_{\max} \nabla u_{\max}|^2$ for the eigenfunction $u_{\max}$ of $\rho_1(\lambda_{\max})$ with various $\eta$
  The maximization criterion (\ref{inequality-max-rho}) is actually satisfied in each case. 
  The larger $\eta$ becomes, the closer $\rho_{\max}$ is to the maximizer with homogeneous Dirichlet boundary condition (Figure \ref{fig-rhoD-square}, c).
  \label{fig-rho-cont-max}
\end{figure}

Finally we consider Problem \ref{prob-rho} with the {\emx mixed} boundary conditions to discuss continuous dependence of optimizers on boundary conditions.
Let $\Omega_1 = (-\pi, \pi)\times (-\pi, \pi)$ and the boundary condition on $\partial \Omega_1$ be
\begin{equation}
  \label{robin-bc}
  u = 0\quad \text{ on }\quad x=\pm \pi \ \ \text{ and }\ \ y=\pi,\quad \eta u + \frac{\partial u}{\partial n} = 0 \quad \text{ on }\quad y=-\pi,
\end{equation}
where $\eta \geq 0$ is a non-negative constant.
In general, 
as in the case of the homogeneous Dirichlet and the homogeneous Neumann boundary conditions,
the unique existence of the boundary value problem of elliptic equation 
\begin{equation*}
  Lu = f \quad \text{ in }\Omega,\quad a(x)u + \frac{\partial u}{\partial n} = 0\quad \text{ on }\partial \Omega
\end{equation*}
is well-known under suitable assumptions,
where $L$ is an elliptic operator,
say, $Lu = -\Delta u + b(x)u$ for a given bounded function $b(x)$ and $a(x)$ is a piecewise continuous function on $\partial \Omega$.
Moreover, the unique solution $u$ depends continuously on $a(x)$ in $L^\infty$ (see e.g.\ \cite{CH}).
In previous subsections we observed that optimizers are dominated by eigenfunctions associated with corresponding eigenvalues.
It is then natural to consider that optimizers also depend continuously on boundary conditions as well as eigenfunctions as solutions of elliptic equations. 
\begin{observation}
  \label{observation-dependence-heat}
  Consider Problem \ref{prob-rho} with the boundary condition (\ref{robin-bc}). 
  Then the region $S_\ast$ given by the optimizer $\rho_\ast$ of $\lambda_1$ depends continuously on $\eta \geq 0$.
  Sample numerical results are shown in Figure \ref{fig-rho-cont-min} and  \ref{fig-rho-cont-max}.
\end{observation}
Although these are just simple examples,
combining this observation with continuous dependence of eigenfunctions on boundary conditions,
we may well expect that one can mathematically prove the continuous dependence of optimizers on boundary conditions in a suitable topology. 
\subsection{Comparative observations -- Problem \ref{prob-sigma}}
\label{section-prob-sigma}
As a comparison with Problem \ref{prob-rho},
we consider the optimization of the first eigenvalue for (\ref{ev-prob-sigma}).
There are many studies on this type of problems, for example \cite{OS}, 
while {\emx almost all such considerations are only with the homogeneous Dirichlet boundary condition on rectangular domain}.
In this section we set $S = \{x\in \Omega \mid \sigma(x) = c\}$ for $\sigma$ under consideration and fix $c = 2$.
\par
\begin{figure}[H]
  \begin{center}
    \begin{tabular}{llll}
      {\bfseries (a)} 
      &
      {\bfseries (b)} 
      & 
      {\bfseries (c)}
      & 
      {\bfseries (d)}
      \cr
      \includegraphics[width=3.0cm]{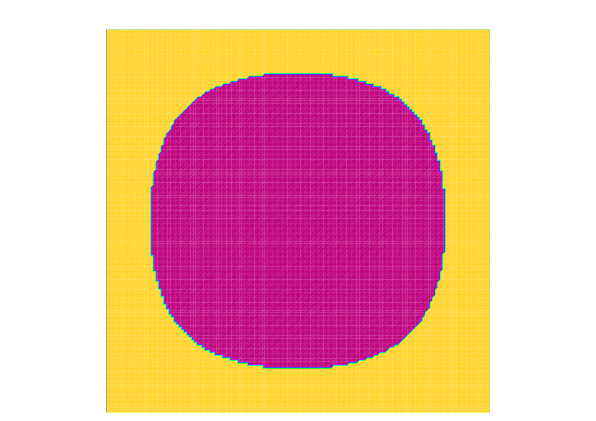}
      &
      \includegraphics[width=3.0cm]{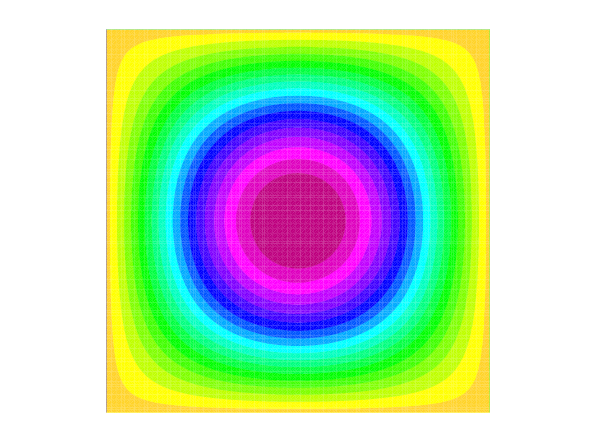}
      &
      \includegraphics[width=3.0cm]{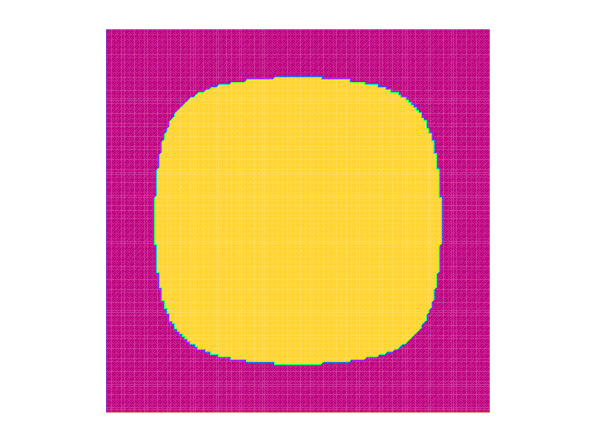}
      &
      \includegraphics[width=3.0cm]{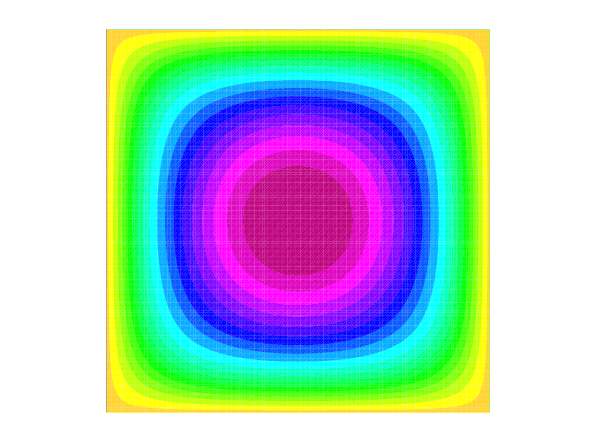}
      \cr
    \end{tabular}
  \end{center}
  \caption{Optimizers and eigenfunctions for Problem \ref{prob-sigma} on $\Omega_1$.}
  (a) minimizer $\sigma_{\min}$ 
  (b) $|u|^2$ of the associated eigenfunction of $\mu_1(\sigma_{\min})$.
  (c) maximizer $\sigma_{\max}$ and 
  (d) $|u|^2$ of the associated eigenfunction of $\mu_1(\sigma_{\max})$.
  \label{fig-sigmaD-square}
\end{figure}
\begin{figure}[H]
  \begin{center}
    \begin{tabular}{llll}
      {\bfseries (a)} 
      &
      {\bfseries (b)} 
      & 
      {\bfseries (c)}
      & 
      {\bfseries (d)}
      \cr
      \includegraphics[width=3.0cm]{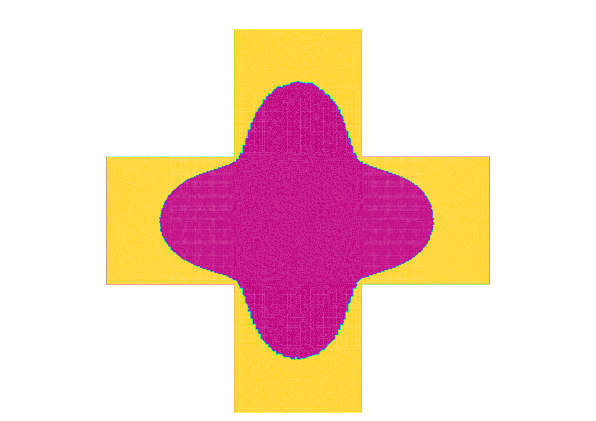}
      &
      \includegraphics[width=3.0cm]{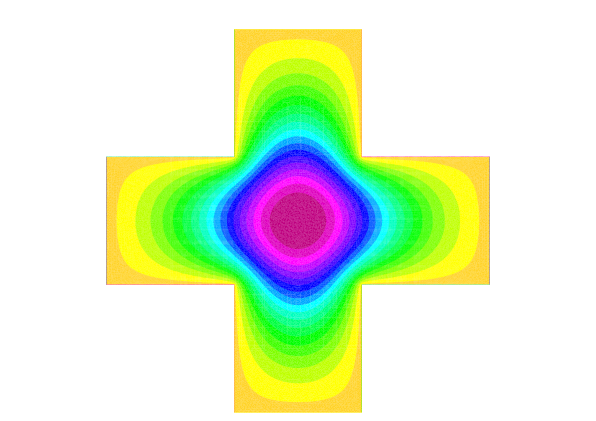}
      &
      \includegraphics[width=3.0cm]{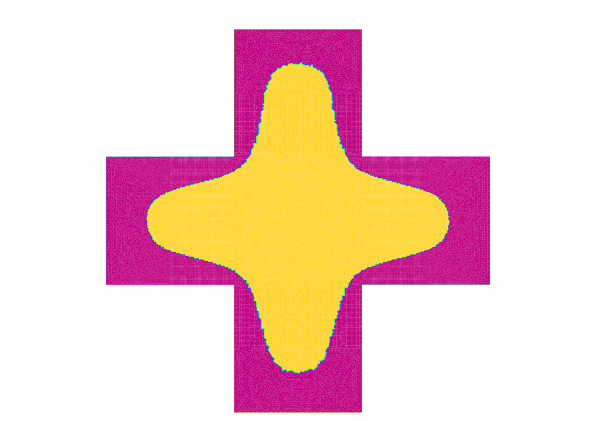}
      &
      \includegraphics[width=3.0cm]{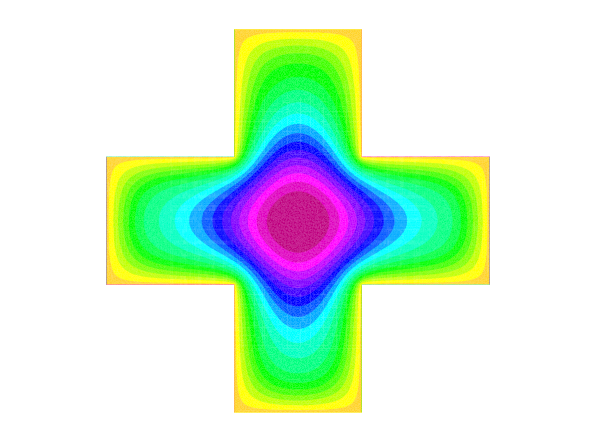}
      \cr
    \end{tabular}
  \end{center}
  \caption{Optimizers and eigenfunctions for Problem \ref{prob-sigma} on $\Omega_2$}
  (a) minimizer $\sigma_{\min}$ 
  (b) $|u|^2$ of the associated eigenfunction of $\mu_1(\sigma_{\min})$.
  (c) maximizer $\sigma_{\max}$ and 
  (d) $|u|^2$ of the associated eigenfunction of $\mu_1(\sigma_{\max})$.
  \label{fig-sigmaD-cross}
\end{figure}
\begin{figure}[H]
  \begin{center}
    \begin{tabular}{llll}
      {\bfseries (a)} 
      &
      {\bfseries (b)} 
      & 
      {\bfseries (c)}
      & 
      {\bfseries (d)}
      \cr
      \includegraphics[width=3.0cm]{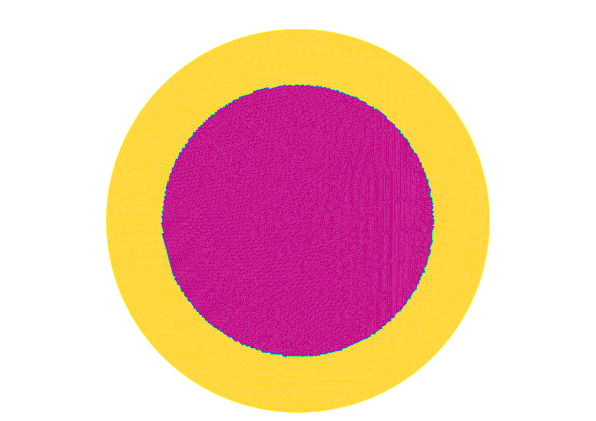}
      &
      \includegraphics[width=3.0cm]{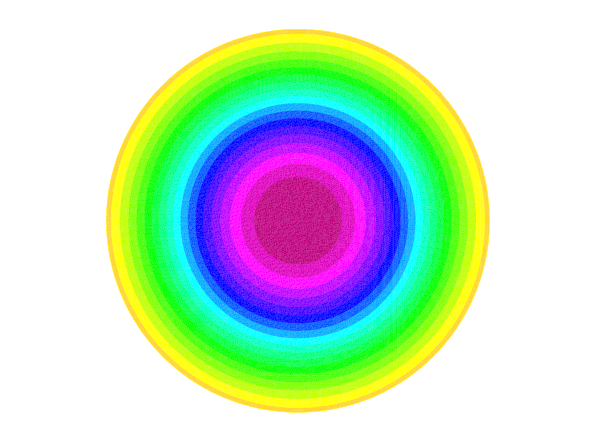}
      &
      \includegraphics[width=3.0cm]{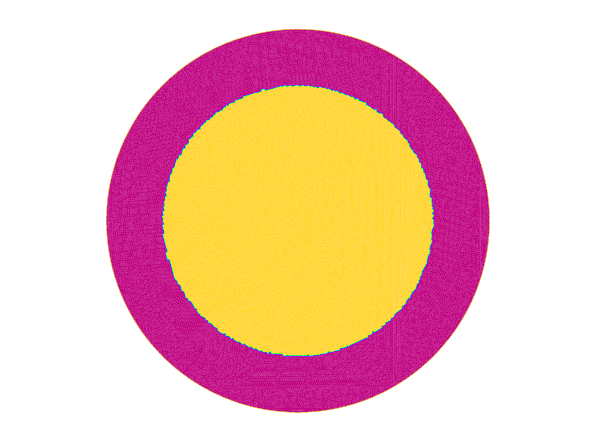}
      &
      \includegraphics[width=3.0cm]{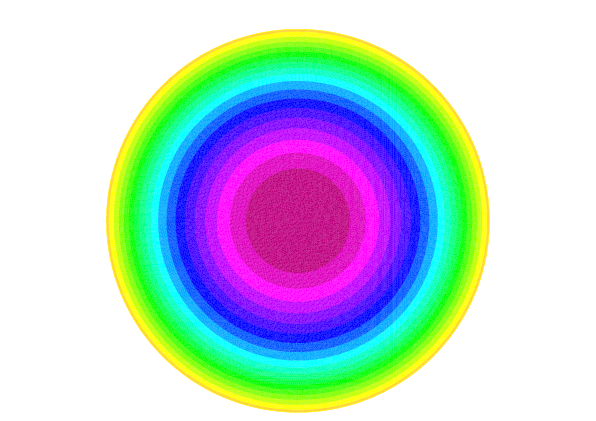}
      \cr
    \end{tabular}
  \end{center}
  \caption{Optimizers and eigenfunctions for Problem \ref{prob-sigma} on $\Omega_3$}
  (a) minimizer $\sigma_{\min}$ 
  (b) $|u|^2$ of the associated eigenfunction of $\mu_1(\sigma_{\min})$.
  (c) maximizer $\sigma_{\max}$ and 
  (d) $|u|^2$ of the associated eigenfunction of $\mu_1(\sigma_{\max})$.
  \label{fig-sigmaD-disk}
\end{figure}
\begin{figure}[H]
  \begin{center}
    \begin{tabular}{llll}
      {\bfseries (a)} 
      &
      {\bfseries (b)} 
      & 
      {\bfseries (c)}
      & 
      {\bfseries (d)}
      \cr
      \includegraphics[width=3.0cm]{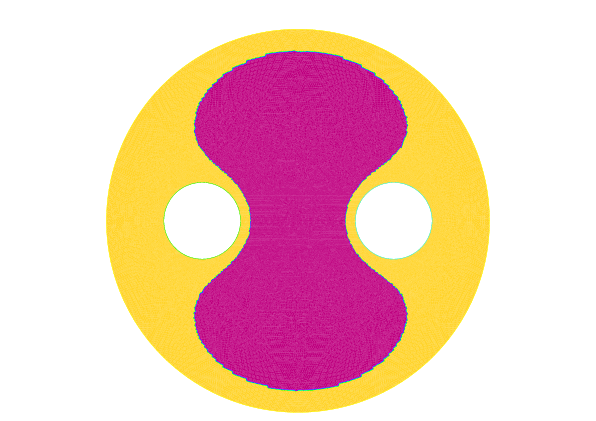}
      &
      \includegraphics[width=3.0cm]{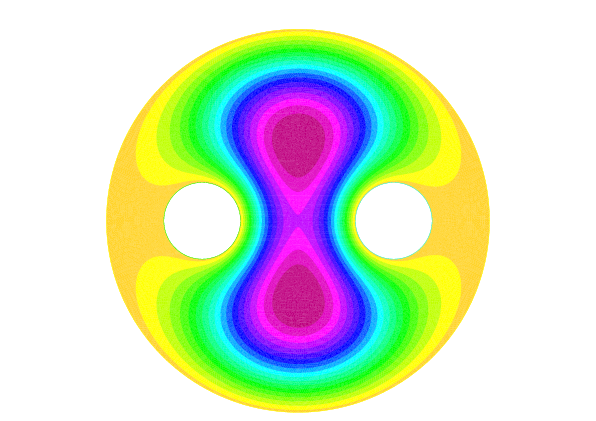}
      &
      \includegraphics[width=3.0cm]{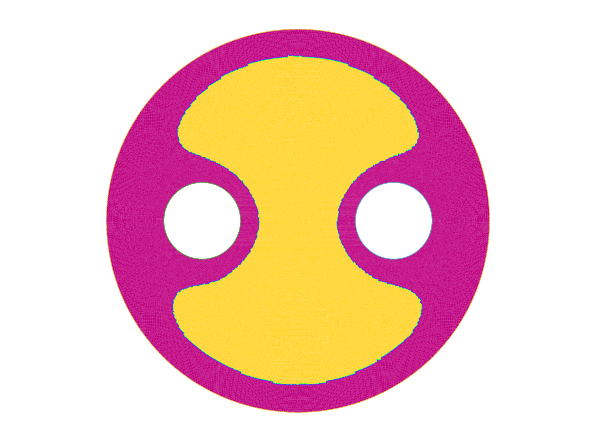}
      &
      \includegraphics[width=3.0cm]{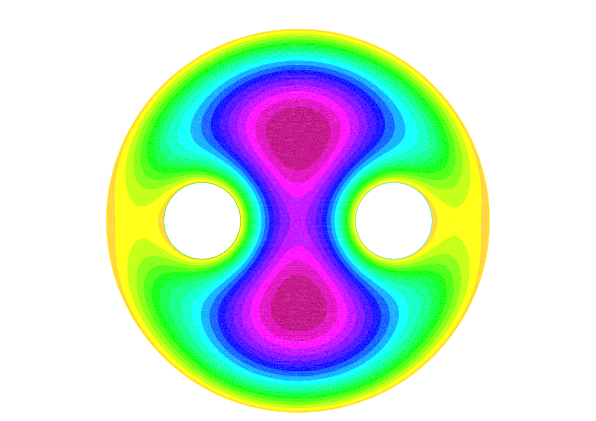}
      \cr
    \end{tabular}
  \end{center}
  \caption{Optimizers and eigenfunctions for Problem \ref{prob-sigma} on $\Omega_4$}
  (a) minimizer $\sigma_{\min}$ 
  (b) $|u|^2$ of the associated eigenfunction of $\mu_1(\sigma_{\min})$.
  (c) maximizer $\sigma_{\max}$ and 
  (d) $|u|^2$ of the associated eigenfunction of $\mu_1(\sigma_{\max})$.
  \label{fig-sigmaD-nonstar}
\end{figure}
\begin{figure}[H]
  \begin{center}
    \begin{tabular}{llll}
      {\bfseries (a)} 
      &
      {\bfseries (b)} 
      & 
      {\bfseries (c)}
      & 
      {\bfseries (d)}
      \cr
      \includegraphics[width=3.0cm]{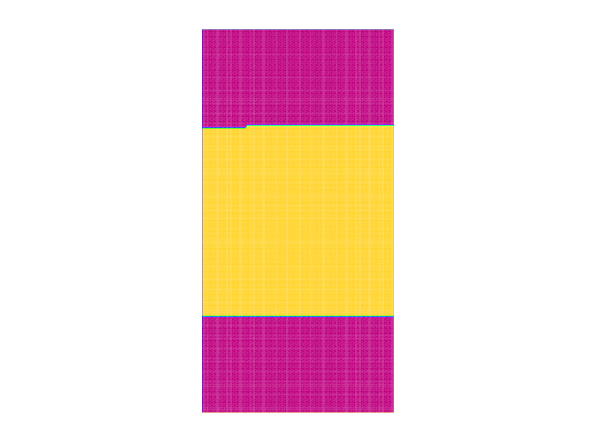}
      &
      \includegraphics[width=3.0cm]{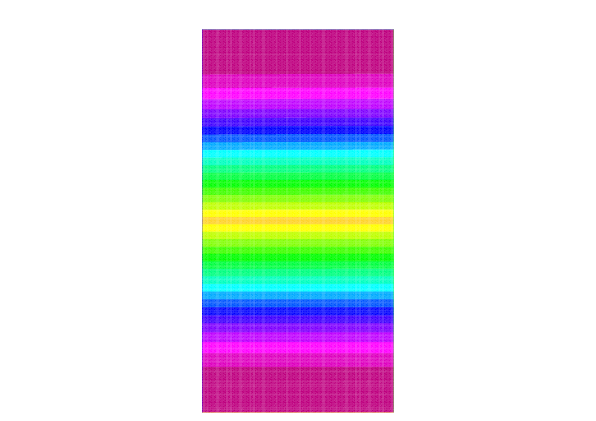}
      &
      \includegraphics[width=3.0cm]{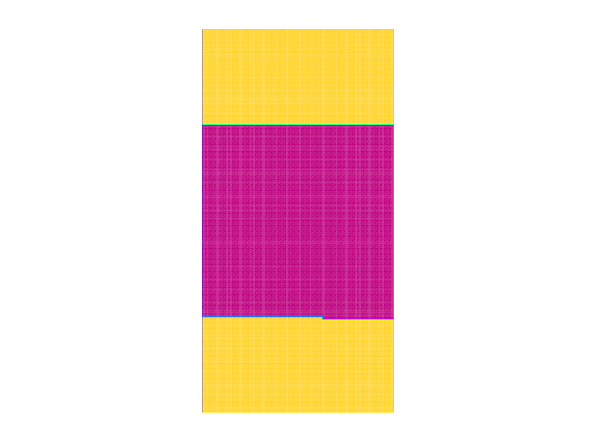}
      &
      \includegraphics[width=3.0cm]{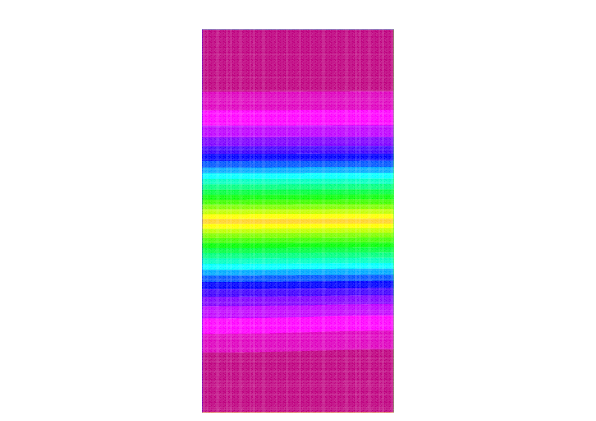}
      \cr
    \end{tabular}
  \end{center}
  \caption{Optimizers and eigenfunctions for Problem \ref{prob-sigma} on $\Omega_5$}
  (a) minimizer $\sigma_{\min}$ 
  (b) $|u|^2$ of the associated eigenfunction of $\mu_1(\sigma_{\min})$.
  (c) maximizer $\sigma_{\max}$ and 
  (d) $|u|^2$ of the associated eigenfunction of $\mu_1(\sigma_{\max})$.
  \label{fig-sigmaN-rect}
\end{figure}
\begin{figure}[H]
  \begin{center}
    \begin{tabular}{llll}
      {\bfseries (a)} 
      &
      {\bfseries (b)} 
      & 
      {\bfseries (c)}
      & 
      {\bfseries (d)}
      \cr
      \includegraphics[width=3.0cm]{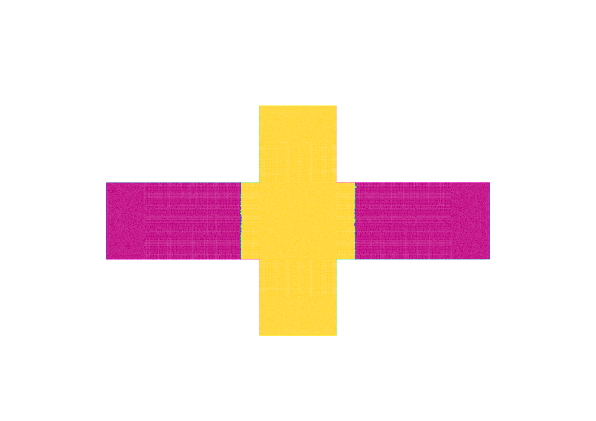}
      &
      \includegraphics[width=3.0cm]{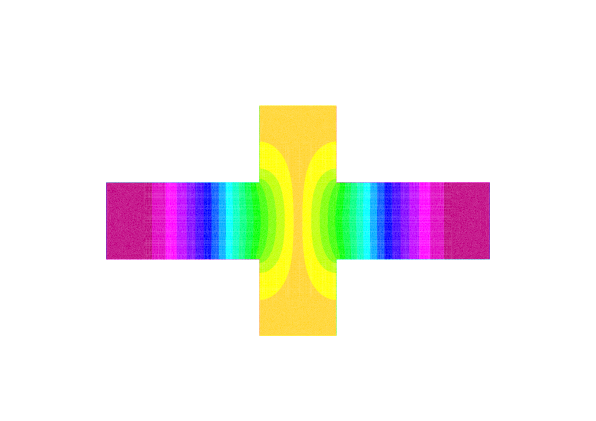}
      &
      \includegraphics[width=3.0cm]{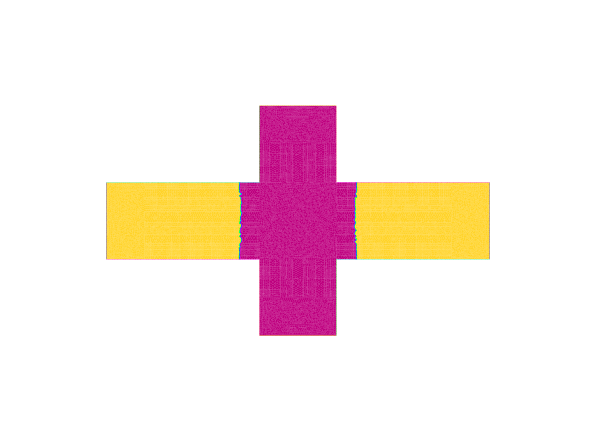}
      &
      \includegraphics[width=3.0cm]{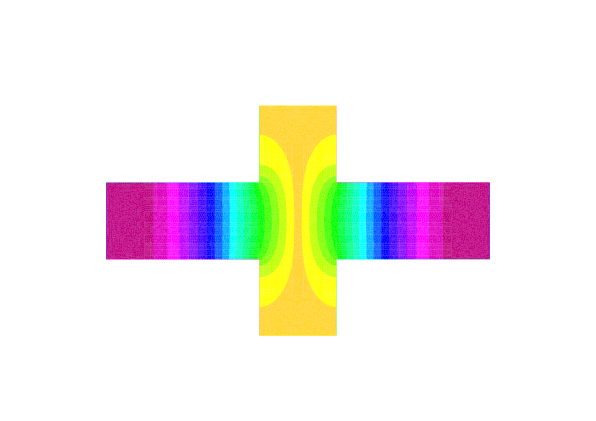}
      \cr
    \end{tabular}
  \end{center}
  \caption{Optimizers and eigenfunctions for Problem \ref{prob-sigma} on $\Omega_6$}
  (a) minimizer $\sigma_{\min}$ 
  (b) $|u|^2$ of the associated eigenfunction of $\mu_1(\sigma_{\min})$.
  (c) maximizer $\sigma_{\max}$ and 
  (d) $|u|^2$ of the associated eigenfunction of $\mu_1(\sigma_{\max})$.
  \label{fig-sigmaN-cross}
\end{figure}
\begin{figure}[H]
  \begin{center}
    \begin{tabular}{llll}
      {\bfseries (a)} 
      &
      {\bfseries (b)} 
      & 
      {\bfseries (c)}
      & 
      {\bfseries (d)}
      \cr
      \includegraphics[width=3.0cm]{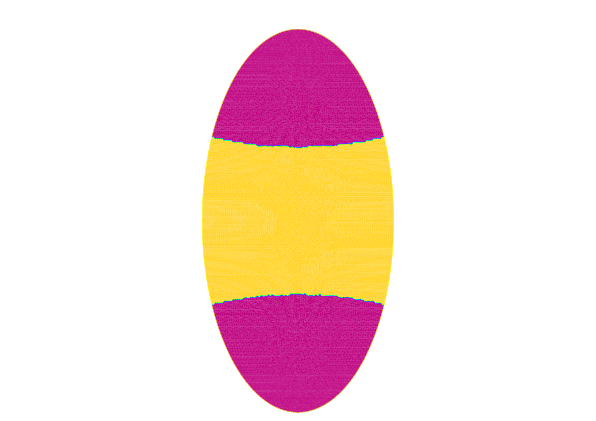}
      &
      \includegraphics[width=3.0cm]{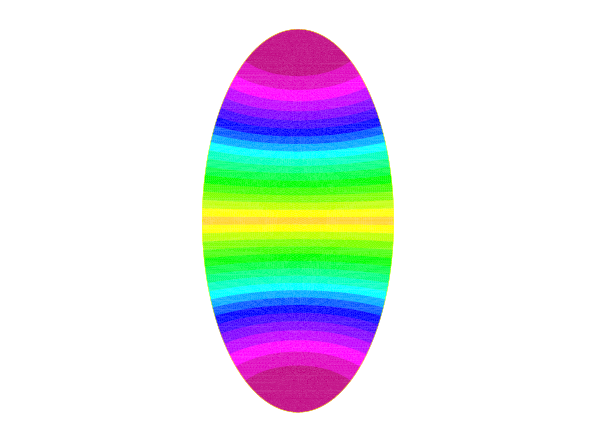}
      &
      \includegraphics[width=3.0cm]{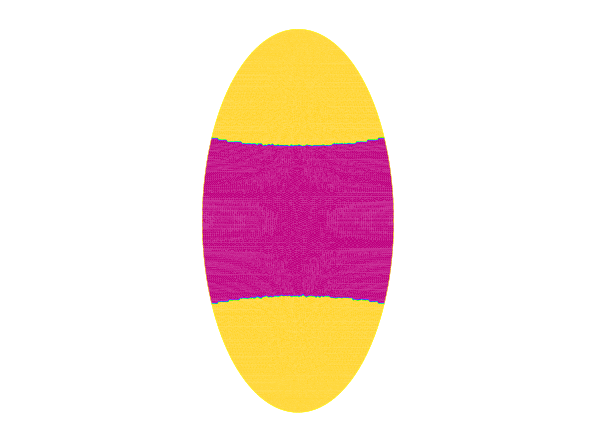}
      &
      \includegraphics[width=3.0cm]{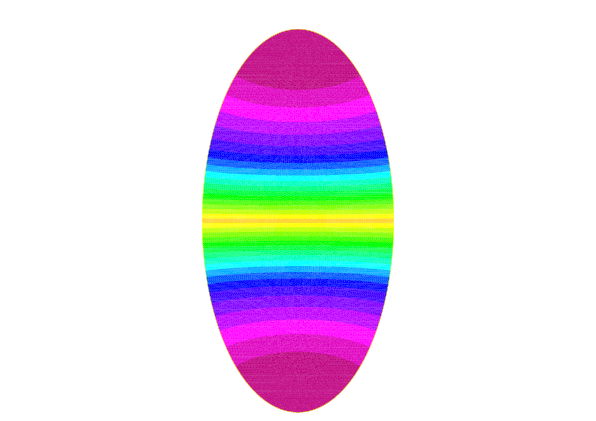}
      \cr
    \end{tabular}
  \end{center}
  \caption{Optimizers and eigenfunctions for Problem \ref{prob-sigma} on $\Omega_7$}
  (a) minimizer $\sigma_{\min}$ 
  (b) $|u|^2$ of the associated eigenfunction of $\mu_1(\sigma_{\min})$.
  (c) maximizer $\sigma_{\max}$ and 
  (d) $|u|^2$ of the associated eigenfunction of $\mu_1(\sigma_{\max})$.
  \label{fig-sigmaN-ellipse}
\end{figure}
\begin{figure}[H]
  \begin{center}
    \begin{tabular}{llll}
      {\bfseries (a)} 
      &
      {\bfseries (b)} 
      & 
      {\bfseries (c)}
      & 
      {\bfseries (d)}
      \cr
      \includegraphics[width=3.0cm]{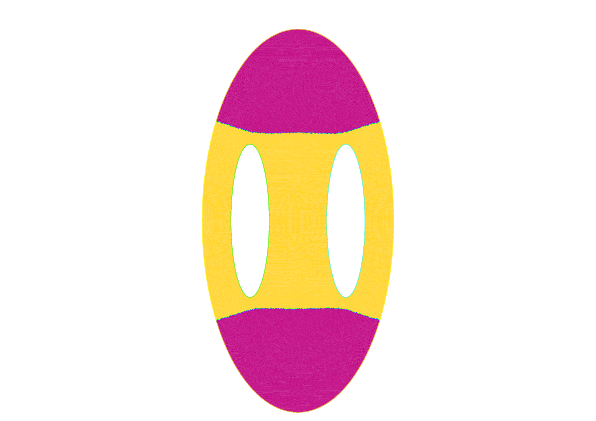}
      &
      \includegraphics[width=3.0cm]{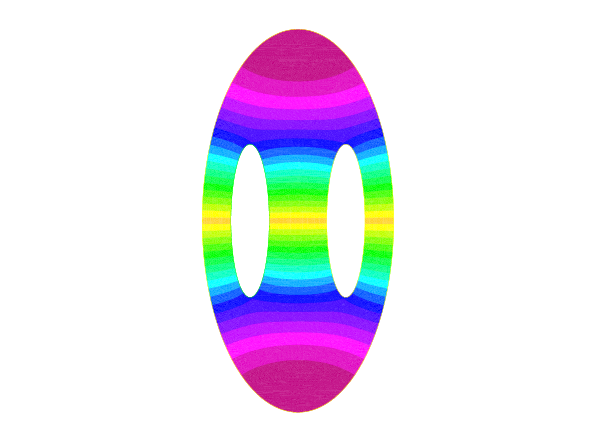}
      &
      \includegraphics[width=3.0cm]{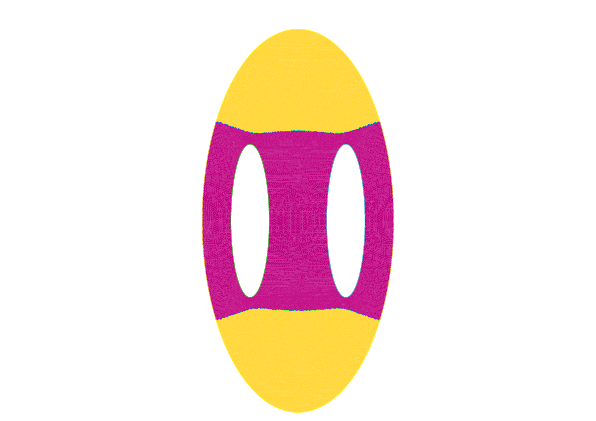}
      &
      \includegraphics[width=3.0cm]{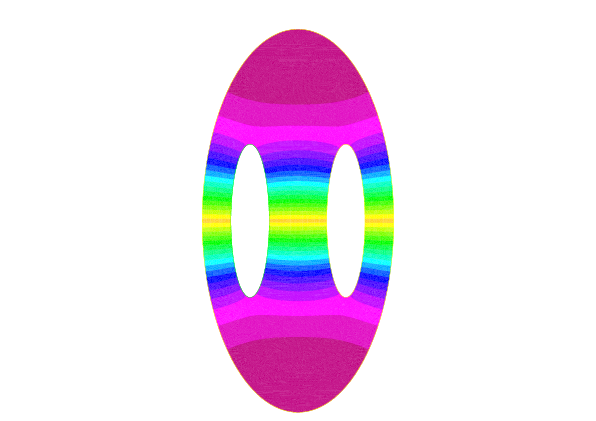}
      \cr
    \end{tabular}
  \end{center}
  \caption{Optimizers and eigenfunctions for Problem \ref{prob-sigma} on $\Omega_8$}
  (a) minimizer $\sigma_{\min}$ 
  (b) $|u|^2$ of the associated eigenfunction of $\mu_1(\sigma_{\min})$.
  (c) maximizer $\sigma_{\max}$ and 
  (d) $|u|^2$ of the associated eigenfunction of $\mu_1(\sigma_{\max})$.
  \label{fig-sigmaN-nonstar}
\end{figure}
\begin{figure}[H]
  \begin{center}
    \begin{tabular}{llll}
      {\bfseries (a)} 
      &
      {\bfseries (b)} 
      & 
      {\bfseries (c)}
      & 
      {\bfseries (d)}
      \cr
      \includegraphics[width=3.0cm]{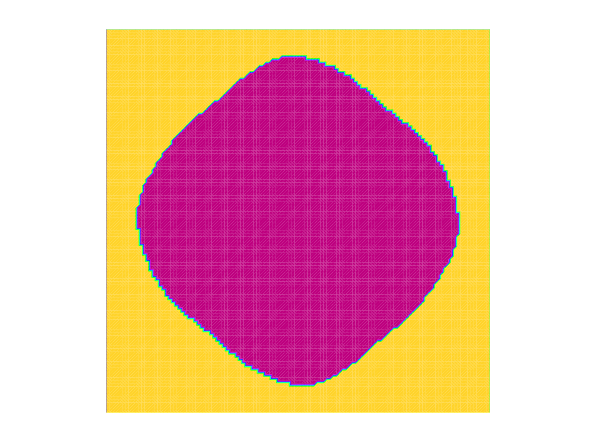}
      &
      \includegraphics[width=3.0cm]{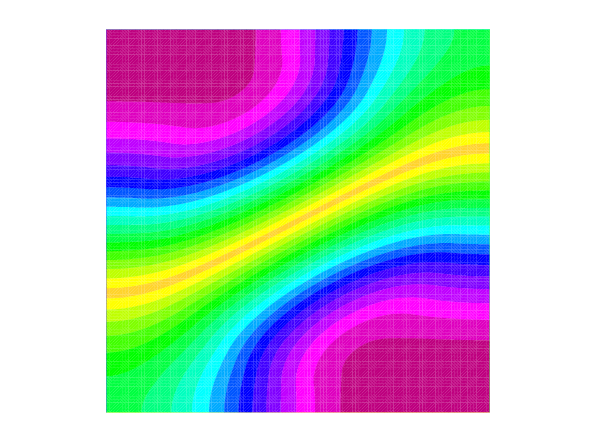}
      &
      \includegraphics[width=3.0cm]{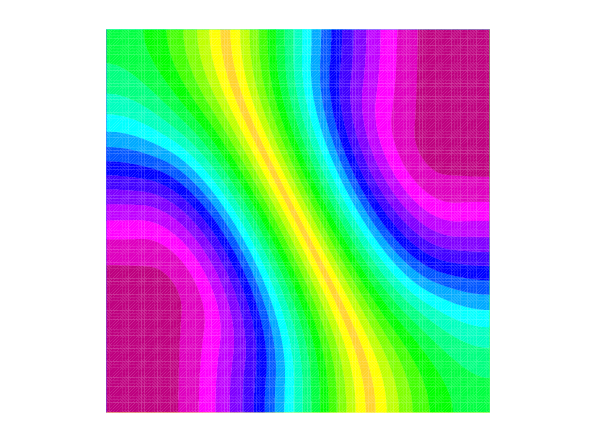}
      &
      \includegraphics[width=3.0cm]{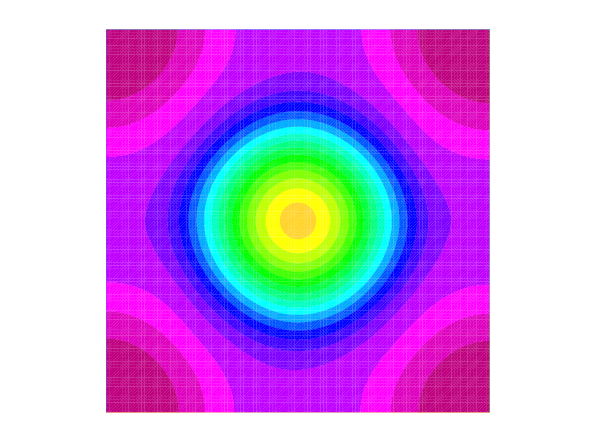}
      \cr
    \end{tabular}
  \end{center}
  \caption{Maximizers and eigenfunctions for Problem \ref{prob-sigma} on $\Omega_1$}
  Figure (a) shows the maximizer $\sigma_{\max}$ of $\mu_1(\sigma)$. 
  Figure (b) and (c) show corresponding $|u_{1,\max}|^2$ and $|u_{2,\max}|^2$, 
  where $u_{1,\max}$ and $u_{2,\max}$ are associated eigenfunctions of $\mu_1(\sigma_{\max})$ and $\mu_2(\sigma_{\max})$ (actually equal to $\mu_1(\sigma_{\max})$), respectively, 
  after the normalization so that $\int_\Omega \sigma_{\max} |u_{i,\max}|^2 = 1$ holds.  
  Figure (d) shows $|u_{1,\max}|^2+|u_{2,\max}|^2$ after normalizations. 
  \label{fig-sigmaN-square}
\end{figure}
\begin{figure}[H]
  \begin{center}
    \begin{tabular}{cccccc}
      (a) 
      &
      \includegraphics[width=2.25cm]{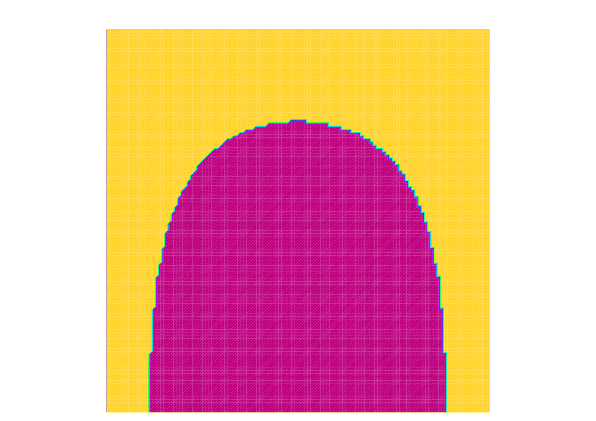}
      &
      \includegraphics[width=2.25cm]{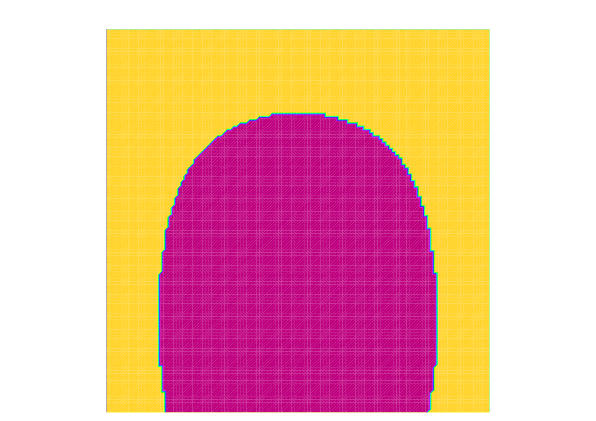}
      &
      \includegraphics[width=2.25cm]{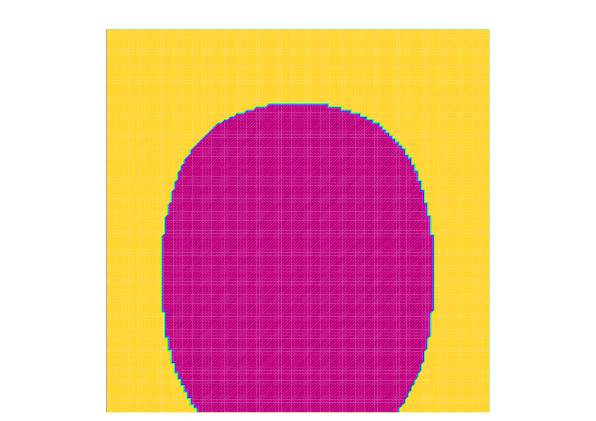}
      &
      \includegraphics[width=2.25cm]{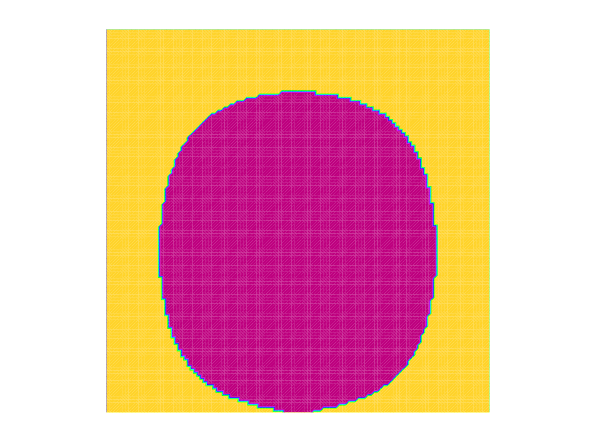}
      &
      \includegraphics[width=2.25cm]{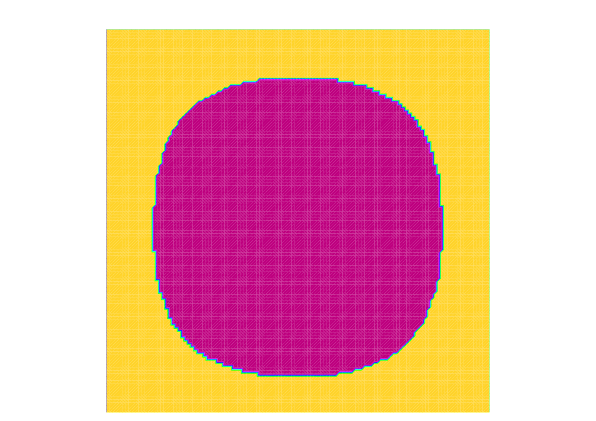}
      \cr
      (b) 
      &
      \includegraphics[width=2.25cm]{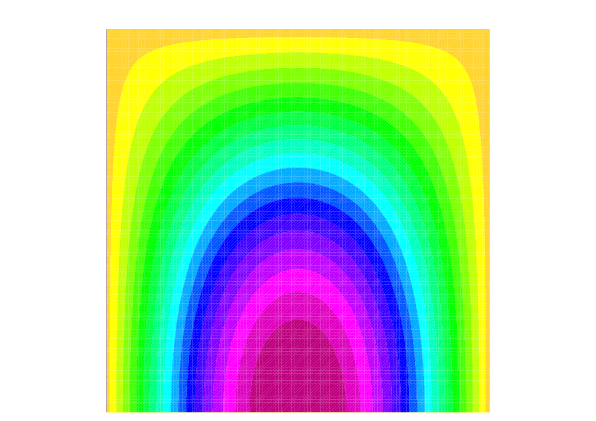}
      &
      \includegraphics[width=2.25cm]{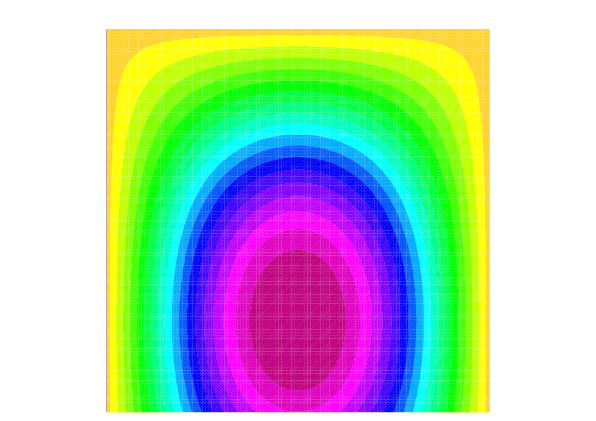}
      &
      \includegraphics[width=2.25cm]{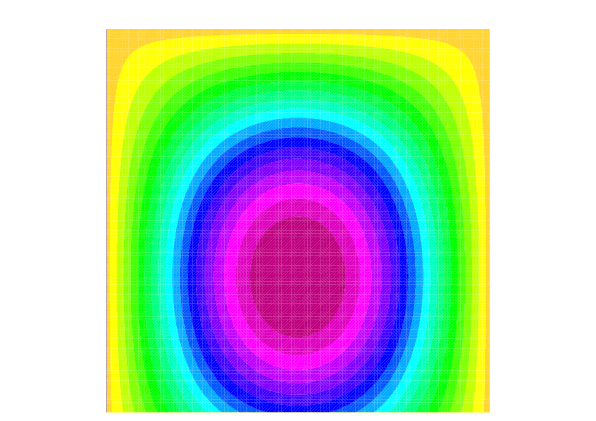}
      &
      \includegraphics[width=2.25cm]{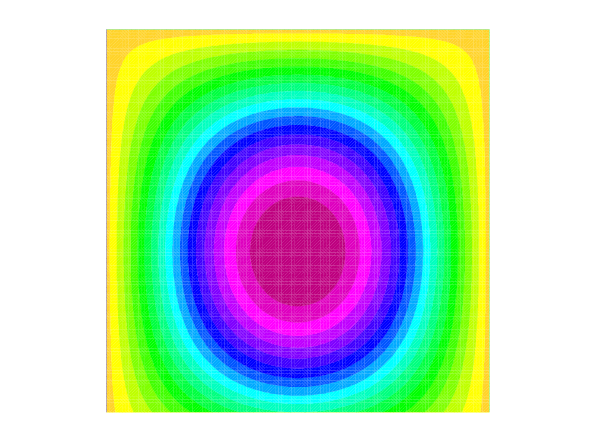}
      &
      \includegraphics[width=2.25cm]{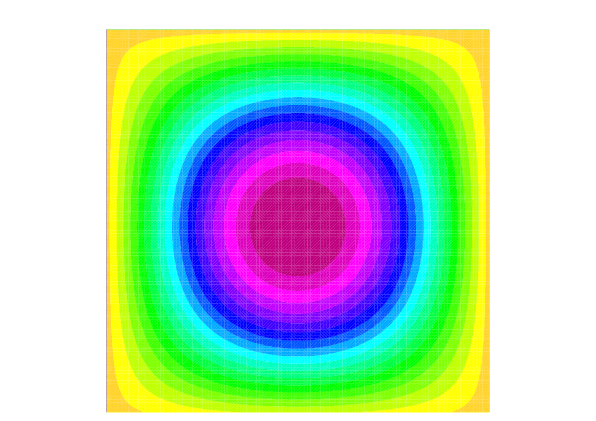}
      \cr
      &
      $\eta=0$
      &
      $\eta=0.2$
      &
      $\eta=0.5$
      &
      $\eta=1.0$
      &
      $\eta=5.0$
      \cr
    \end{tabular}
  \end{center}
  \caption{Minimizer for Problem \ref{prob-sigma} on $\Omega_1$}
  (a) shows the minimizer $\sigma_{\min}$ and 
  (b) shows $|u_{\min}|^2$ for the eigenfunction $u_{\min}$ of $\mu_1(\sigma_{\min})$ with various $\eta$
  The minimization criterion (\ref{inequality-min-sigma}) is actually satisfied in each case. 
  The larger $\eta$ becomes, the closer $\sigma_{\min}$ is to the minimizer with homogeneous Dirichlet boundary condition (Figure \ref{fig-sigmaD-square} a).
  \label{fig-sigma-cont-min}
\end{figure}
\begin{figure}[H]
  \begin{center}
    \begin{tabular}{cccccc}
      (a) 
      &
      \includegraphics[width=2.25cm]{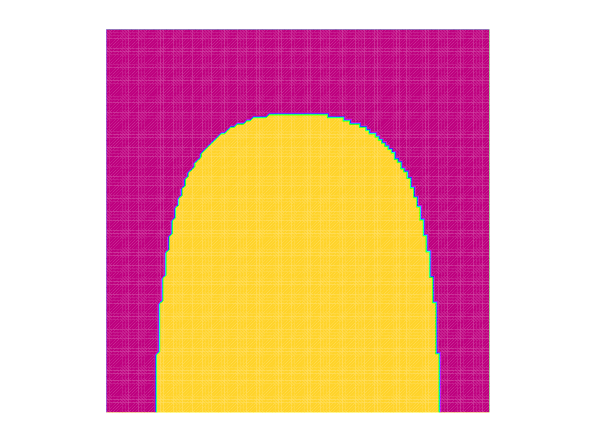}
      &
      \includegraphics[width=2.25cm]{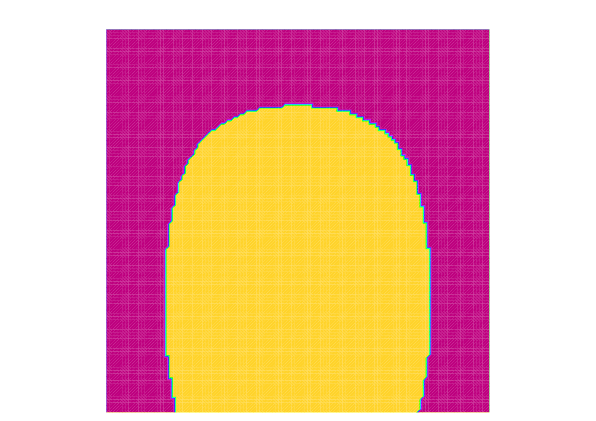}
      &
      \includegraphics[width=2.25cm]{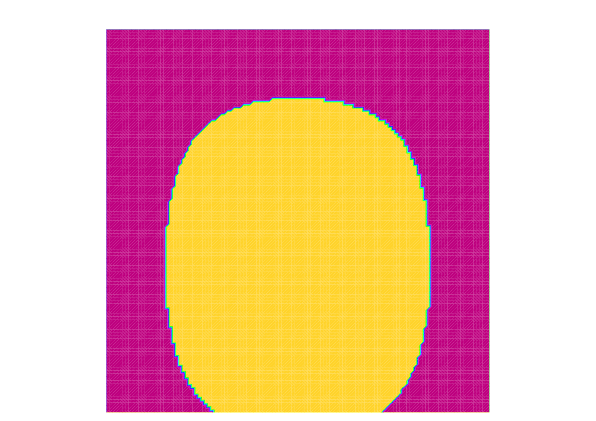}
      &
      \includegraphics[width=2.25cm]{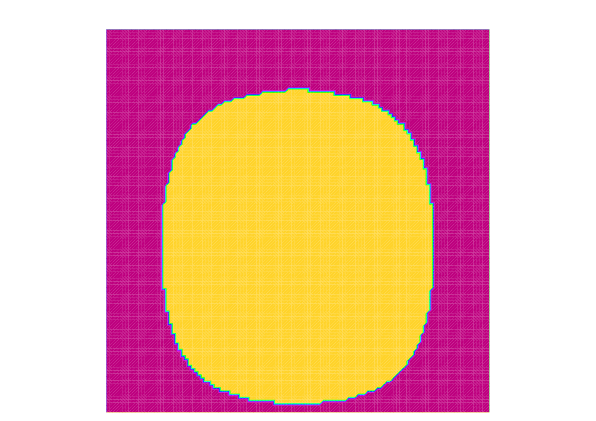}
      &
      \includegraphics[width=2.25cm]{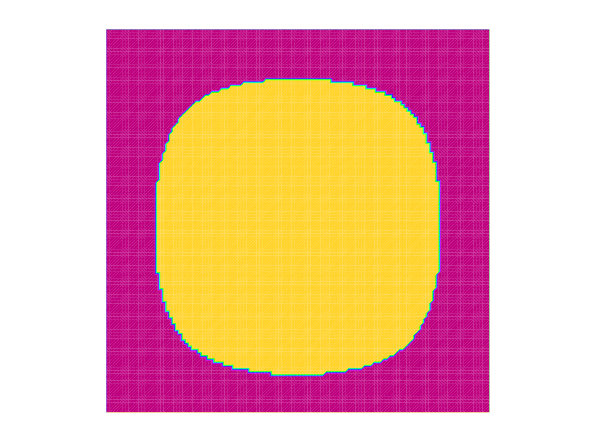}
      \cr
      (b) 
      &
      \includegraphics[width=2.25cm]{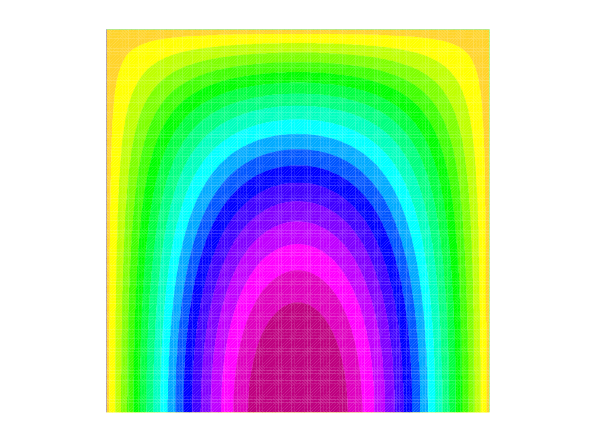}
      &
      \includegraphics[width=2.25cm]{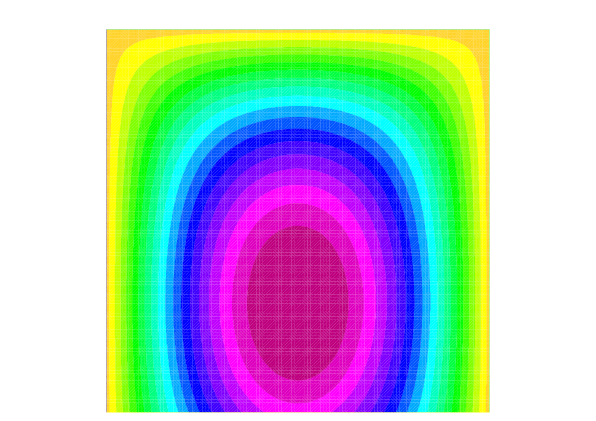}
      &
      \includegraphics[width=2.25cm]{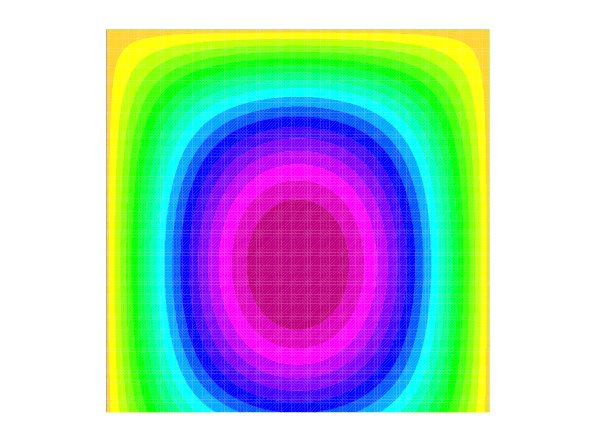}
      &
      \includegraphics[width=2.25cm]{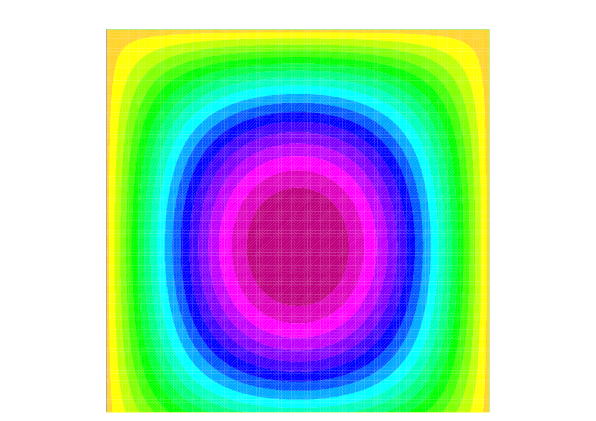}
      &
      \includegraphics[width=2.25cm]{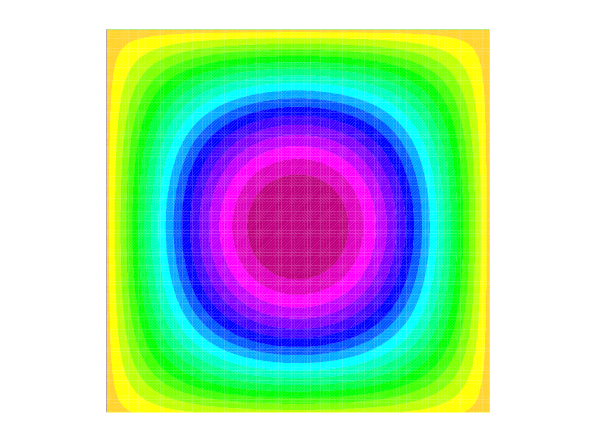}
      \cr
      &
      $\eta=0$
      &
      $\eta=0.2$
      &
      $\eta=0.5$
      &
      $\eta=1.0$
      &
      $\eta=5.0$
      \cr
    \end{tabular}
  \end{center}
  \caption{Maximizer for Problem \ref{prob-sigma} on $\Omega_1$}
  (a) shows the minimizer $\sigma_{\max}$ and 
  (b) shows $|u_{\max}|^2$ for the eigenfunction $u_{\min}$ of $\mu_1(\sigma_{\max})$ with various $\eta$
  The maximization criterion (\ref{inequality-max-sigma}) is actually satisfied in each case. 
  The larger $\eta$ becomes, the closer $\sigma_{\max}$ is to the maximizer with homogeneous Dirichlet boundary condition (Figure \ref{fig-sigmaD-square}, c).
  \label{fig-sigma-cont-max}
\end{figure}

Calculations similar to those in subsection \ref{section-prob-rho} with
\begin{equation*}
  v_0(x) = \frac{\mu_1(\phi)(c-1)}{\int_\Omega \sigma u_\phi^2dx }|u_\phi(x)|^2
\end{equation*}
in (\ref{v0}) yield numerical results listed in 
Figures \ref{fig-sigmaD-square}, \ref{fig-sigmaD-cross}, 
\ref{fig-sigmaD-disk}, \ref{fig-sigmaD-nonstar}, \ref{fig-sigmaN-rect}, 
\ref{fig-sigmaN-cross}, \ref{fig-sigmaN-ellipse}, 
\ref{fig-sigmaN-nonstar} and \ref{fig-sigmaN-square}.
These results imply if $\sigma \in \mathcal{K}$ minimizes or maximizes $\mu_1$ and 
if $\mu_1$ is {\emx simple} (see Table \ref{table-simple-sigma}) then the corresponding eigenfunction $u_1$ associated with $\mu_1$ satisfies
\begin{align}
  \label{inequality-min-sigma}
  \text{Minimization : }\quad &|u_1(x)| \leq |u_1(y)|\quad \text{for }x\in S\text{ and }y\in S^c,\\
  \label{inequality-max-sigma}
  \text{Maximization : }\quad &|u_1(y)| \leq |u_1(x)|\quad \text{for }x\in S\text{ and }y\in S^c.
\end{align}
Examples of mixed boundary value problem are also shown in Figures \ref{fig-sigma-cont-min} and \ref{fig-sigma-cont-max},
which give us the same observation as in Observation \ref{observation-dependence-heat}. 
Note that inequalities (\ref{inequality-min-sigma}) and (\ref{inequality-max-sigma}) are
exactly the same correspondences as the optimization criteria for the homogeneous Dirichlet boundary value problems:
Theorem \ref{preceding-sigma}. 
One of key considerations of these observations is that {\emx $\mu_1$ is assumed to be simple}.
Table \ref{table-simple-sigma} shows the ratio between $\mu_1(\sigma)$ and the second eigenvalue $\mu_2(\sigma)$.
By discussions similar to subsection \ref{section-prob-rho},
we obtain the following eigenvalue optimization criteria in Problem \ref{prob-sigma} including eigenvalues with the multiplicity two.
\begin{observation}
\label{observation-criteria-sigma}
  Consider Problem \ref{prob-sigma} in bounded domain $\Omega \subset \mathbb{R}^n$.
  If $\sigma_{\min}\in \mathcal{K}$ is the minimizer of $\mu_1$ in $\mathcal{K}$, then the eigenfunction $u_1$ associated with $\mu_1$ satisfies (\ref{inequality-min-sigma}).
  Similarly, if $\sigma_{\max}\in \mathcal{K}$ is the maximizer of $\mu_1$ in $\mathcal{K}$ and if $\mu_1(\sigma_{\max})$ is simple,
  then the eigenfunction $u_1$ associated with $\mu_1(\sigma_{\max})$ satisfies (\ref{inequality-max-sigma}).
  If $\sigma_{\max}\in \mathcal{K}$ is the maximizer of $\mu_1$ in $\mathcal{K}$ and if $\mu_1(\sigma_{\max})$ has multiplicity two,
  then eigenfunctions $u_1$ and $u_2$ associated with $\mu_1(\sigma_{\max})$ satisfies
  \begin{equation}
    \label{inequality-max-sigma-2}
    |u_1(y)|^2+ |u_2(y)|^2 \leq |u_1(x)|^2 + |u_2(x)|^2\quad \text{for }x\in S\text{ and }y\in S^c
  \end{equation}
  under the normalization $\int_\Omega \sigma |u_i|^2 dx = 1$ ($i=1$, $2$), 
  which can be confirmed in Figure \ref{fig-rhoN-square}.
\end{observation}
\par
Next we consider geometry of $S_\ast$.
In the case of the homogeneous Dirichlet boundary condition,
there is a mathematical result which describes the geometric property of $S_{\ast}$: Theorem \ref{preceding-geom}. 
In the case of the homogeneous Neumann boundary value problems in Problem \ref{prob-sigma},
unlike Dirichlet boundary value problems,
associated eigenfunctions do not have identical sign in $\Omega$.
We thus consider connected components of off-zero-level set of eigenfunctions, which is called {\em nodal domains}.
It is well-known that the eigenfunction associated with $\mu_1(\sigma)$ 
has exactly two nodal domains in the case of the homogeneous Neumann boundary value problems (see e.g.\ \cite{Ban}). 
\par
Assume that $u$ is the eigenfunction associated with $\mu_1(\sigma)$ under the homogeneous Neumann boundary condition.
Then the function $v_1$ and $v_2$ obtained by restricting $u$ on each nodal domain,
say,  $\Omega_1$ and $\Omega_2$, are eigenfunctions of the same equation with an identical sign in $\Omega_i$, respectively,
with
\begin{equation*}
\frac{\partial v_i}{\partial n} = 0\quad \text{ on }\partial \Omega_i\setminus \{u=0\},\quad v_i = 0\quad \text{ on } \{u=0\}
\end{equation*}
In this case, the same discussion as in the case of homogeneous Dirichlet boundary condition can be applied to
analyzing properties of $v_i$ (see e.g.\ \cite{T}) including our criteria (\ref{inequality-min-sigma}) and (\ref{inequality-max-sigma}). 
\par
\begin{figure}[H]
  \begin{center}
    \begin{tabular}{llll}
      {\bfseries (a)} 
      &
      {\bfseries (b)} 
      \cr
      \includegraphics[width=3.0cm]{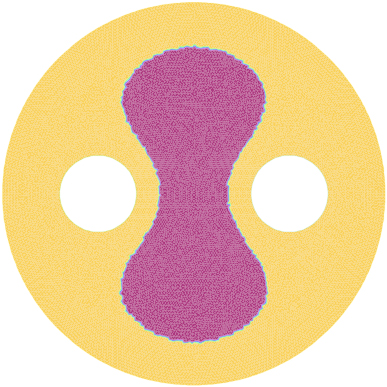}
      &
      \includegraphics[width=3.0cm]{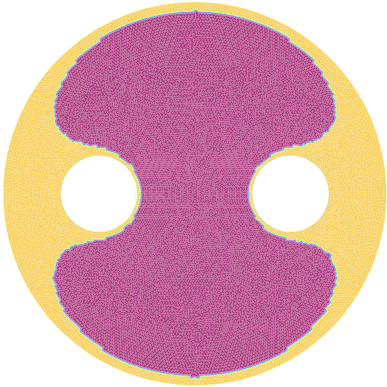}
      \cr
    \end{tabular}
  \end{center}
  \caption{{\bf Geometry of $\{x \mid \sigma_{\min}(x)=c\}$ for a non-star-shaped region $\Omega_4$}}
  Minimization in Problem \ref{prob-sigma} for the non-star-shaped region $\Omega_4$ (cf. Figure \ref{fig-sigmaD-nonstar}) with homogeneous Dirichlet boundary condition is considered.
  Figure (a) is the minimizer $\sigma_{\min}$ with the volume constraint ratio $m_0 = 0.3$ and 
  (b) is $\sigma_{\min}$ with the volume constraint ratio $m_0 = 0.7$.
  The super-level set $S_{\min} = \{x\in \Omega_4 \mid \sigma_{\min}(x)=c\}$ is star-shaped in case of (a), which is not the case of (b).
  \label{nonstar}
\end{figure}

On the other hand, in Figure \ref{fig-sigmaN-rect}, \ref{fig-sigmaN-cross}, \ref{fig-sigmaN-ellipse} and \ref{fig-sigmaN-nonstar},
$\Omega$ is symmetric with respect to the set $\{x\in \Omega \mid u(x) = 0\}$. 
In Figure \ref{fig-sigmaN-square}, $\Omega_1$ is rotationally symmetric with respect to the origin, namely,
$\{x\mid |u_1(x)|^2 + |u_2(x)|^2 = 0\}$.
In such cases, we can see that $S_\ast$ is also symmetric with respect to the null set of $u$ (or $|u_1|^2+|u_2|^2$). Consequently we observe the following: 
\begin{observation}
\label{observation-geometry-sigma}
  Consider Problem \ref{prob-sigma} with the homogeneous Neumann boundary condition.
  Assume that $\sigma_{\min} \in \mathcal{K}$ is the minimizer of $\mu_1$ in $\mathcal{K}$ and $S_{\min}  = \{x\in \Omega \mid \sigma_{\min}(x) = c\}$.
  If each nodal domain $\Omega_i$ ($i=1$, $2$) is convex and symmetric in $N$ orthogonal directions,
  then 
  $\{x \in \Omega_i \mid \sigma_{\min}(x)=c\}
  =
  S_{\min} \mid_{\Omega_i}$ 
  also has the same properties.
  \par
  Similarly, assume that $\sigma_{\max} \in \mathcal{K}$ is the maximizer of $\mu_1$ in $\mathcal{K}$ and $S_{\max}  = \{x\in \Omega \mid \sigma_{\max}(x) = c\}$.
  Then, under the same assumption as minimizers, the same statements hold for the set $\{x\in \Omega_i \mid \sigma_{\max}(x) = 1\} = S_{\max}^c\mid_{\Omega_i}$.
  \par
  Finally, if $\Omega$ is symmetric with respect to the null set of eigenfunctions 
  $\{x\in \Omega\mid u(x)=0\}$ (or $\{x\in \Omega \mid |u_1(x)|^2 + |u_2(x)|^2 = 0\}$ 
  in the case that $\mu_1(\sigma_\ast)$ is not simple), 
  then $S_\ast$ and $S_\ast^c$ are also symmetric.
  These observations are numerically confirmed in Figure \ref{fig-sigmaN-rect}, \ref{fig-sigmaN-cross}, 
  \ref{fig-sigmaN-ellipse}, \ref{fig-sigmaN-nonstar} and \ref{fig-sigmaN-square}. 
\end{observation}
We also discuss the geometry of $S_\ast$ in Problem \ref{prob-sigma} in the case that $\Omega$ is {\emx not even star-shaped}. 
Our numerical results partially answer the inheritance problem in this case, 
as shown in Figure \ref{fig-sigmaD-nonstar}, \ref{nonstar} and \ref{fig-sigmaN-nonstar}.
\par
Theorem \ref{preceding-geom} 
refers to the star-shapedness of $S_\ast$ {\emx only in the case that $\Omega$ is convex}.
Our numerical results newly suggest if $\Omega$ is star-shaped then $S_\ast$ is also star-shaped (see Figure \ref{fig-sigmaD-cross}).
However, if $\Omega$ is not star-shaped, neither $S_\ast$ nor $S_\ast^c$ are necessarily star-shaped.
In Figure \ref{nonstar}, minimization in Problem \ref{prob-sigma} for $\Omega_8$ with the homogeneous Dirichlet boundary condition is considered.
Only difference between two figures is the ratio of the volume constraint in (\ref{constraint}).
One of them is $m_0 = 0.3$ and the other is $m_0=0.7$. 
In the case of $m_0 = 0.3$, $S_{\min} = \{x\in \Omega \mid \sigma_{\min}(x) = c\}$ for the minimizer $\sigma_{\min}$ of $\mu_1(\sigma)$ is star-shaped.
On the other hand, in the case of $m_0 = 0.7$, neither $S_{\min}$ nor $S_{\min}^c$ is star-shaped. 
As a conclusion, 
if $\Omega$ is not star-shaped, in general, neither the optimized domain $S_\ast$ nor its complement $S_\ast^c$ is start-shaped.
\par
As for symmetry, if $\Omega$ is star-shaped and symmetric in a certain direction, then both $S_\ast$ and $S_\ast^c$ have the same symmetry both in the Dirichlet and the Neumann boundary problems. 
Summarizing our arguments we obtain the following observation:
\begin{observation}
\label{observation-symmetry-sigma}
  For Problem \ref{prob-sigma}, 
  let $\sigma_\ast$ be the optimizer of $\mu_1(\sigma)$,
  $u_\ast$ be the associated eigenfunction of $\mu_1(\sigma_\ast)$ and $S_\ast = \{x\in \Omega \mid \sigma_\ast(x) = c\}$.
  If $\Omega$ is not star-shaped, 
  the star-shapedness of $S_\ast$ and $S_\ast^c$
  generally depends on $\Omega$ and $m_0$.
  \par
  If $\Omega$ is star-shaped and
  is symmetric with respect to a certain direction or rotationally symmetric,
  then $S_\ast$ and $S_\ast^c$ also have the same symmetry.
\end{observation}
We finally remark that, even if $\Omega$ is symmetric and non-star-shaped, 
there are cases that $S_\ast$ and $S_\ast^c$ have symmetry which $\Omega$ possesses, as are shown in Figure \ref{fig-sigmaD-nonstar} and \ref{fig-sigmaN-nonstar}.
\subsection{Convergence rate}
\begin{figure}[H]
  \begin{center}
    \begin{tabular}{cccccc}
      (a)
      &
      \includegraphics[width=2.25cm]{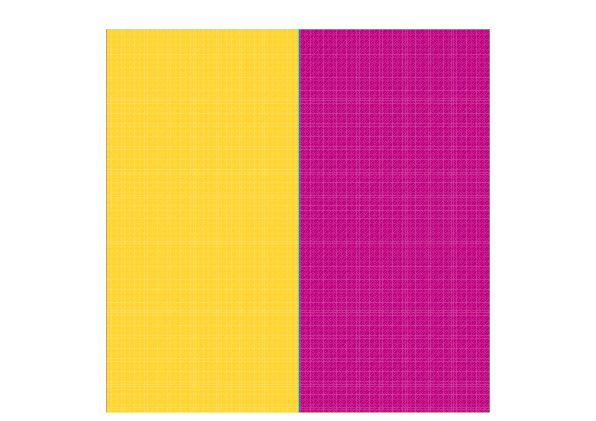}
      &
      \includegraphics[width=2.25cm]{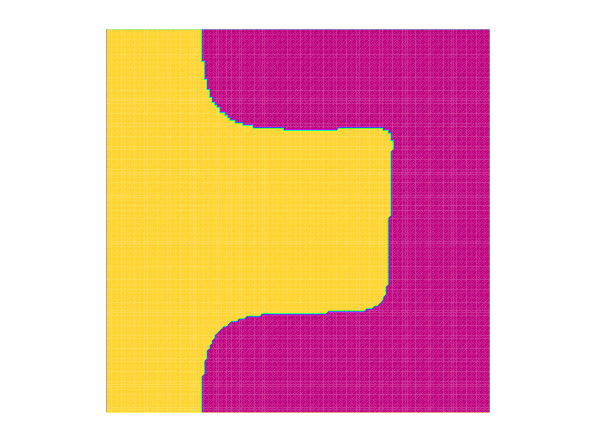}
      &
      \includegraphics[width=2.25cm]{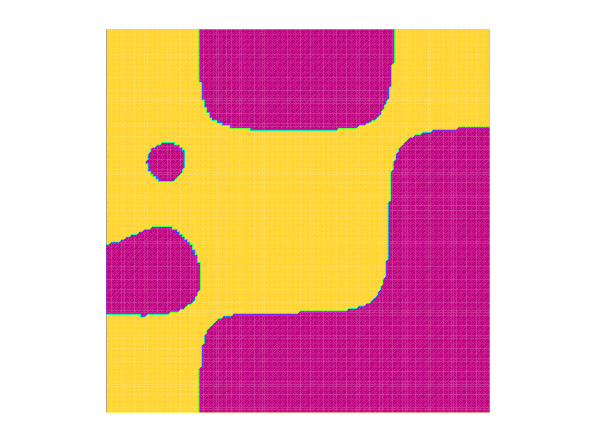}
      &
      \includegraphics[width=2.25cm]{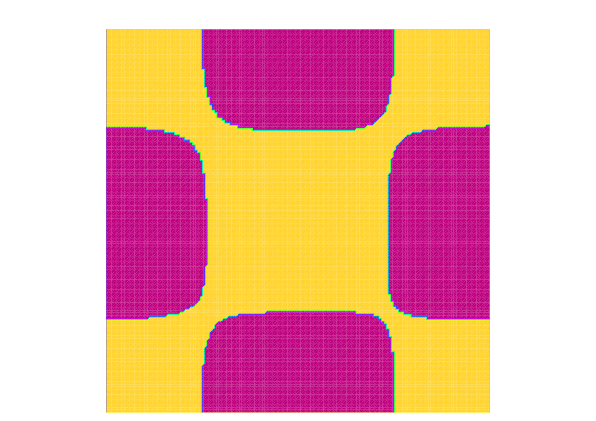}
      &
      \includegraphics[width=2.25cm]{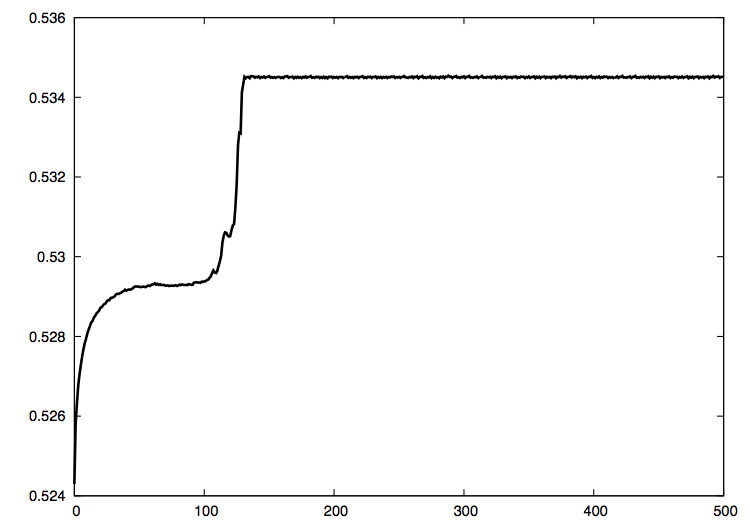}
      \cr
      &
      $t=0$ & $t=100$ & $t=125$ & $t=300$ \cr
      (b)
      &
      \includegraphics[width=2.25cm]{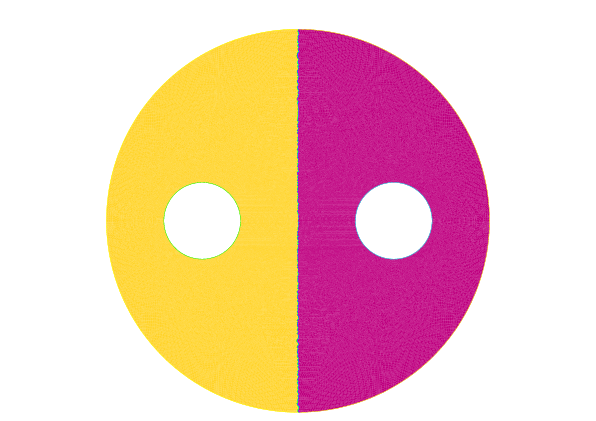}
      &
      \includegraphics[width=2.25cm]{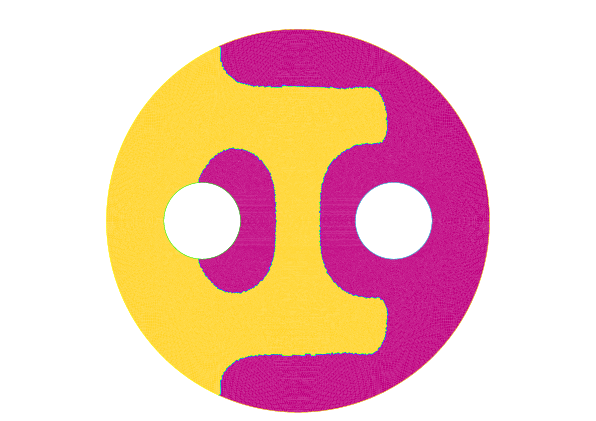}
      &
      \includegraphics[width=2.25cm]{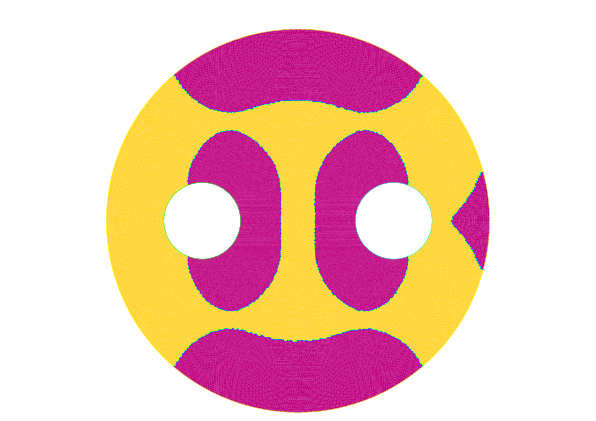}
      &
      \includegraphics[width=2.25cm]{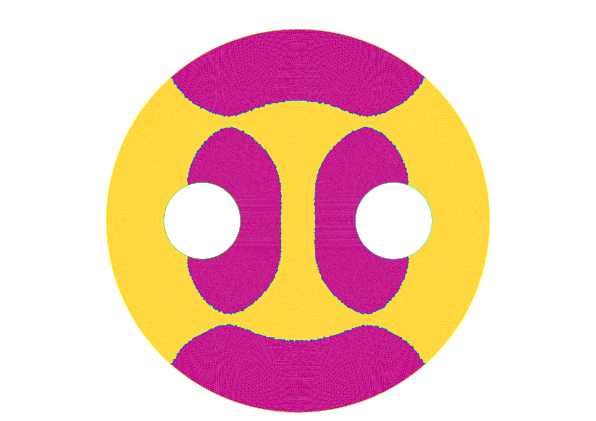}
      &
      \includegraphics[width=2.25cm]{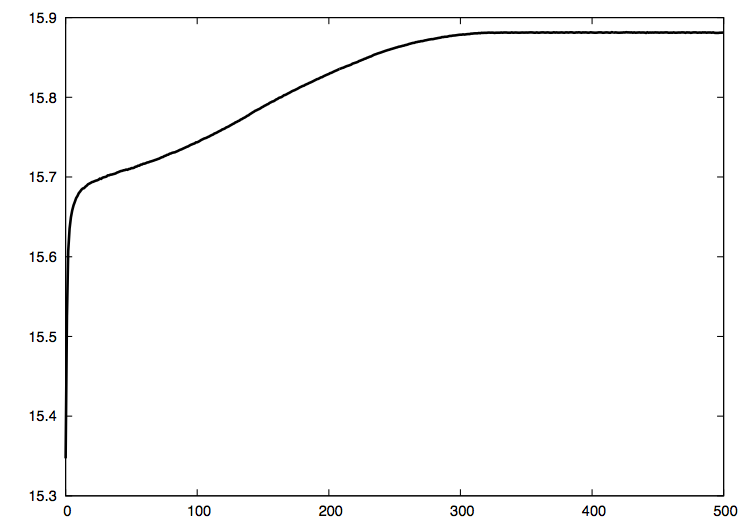}
      \cr
      &
      $t=0$ & $t=20$ & $t=250$ & $t=450$ \cr
    \end{tabular}
  \end{center}
  \caption{Problem \ref{prob-rho} with Dirichlet boundary}
  In all cases $c=1.1$ and $m_0 = 0.5$ are fixed. 
  The rightmost graph is the evolution of $\lambda_1(\rho)$ as we solve (\ref{viscos}).
  (a) : Maximization of $\lambda_1(\rho)$ on $\Omega_1$. 
  (b) : Maximization of $\lambda_1(\rho)$ on $\Omega_4$.\label{fig-rhoD-conv}
\end{figure}

\begin{figure}[H]
  \begin{center}
    \begin{tabular}{cccccc}
      (a)
      &
      \includegraphics[width=2.25cm]{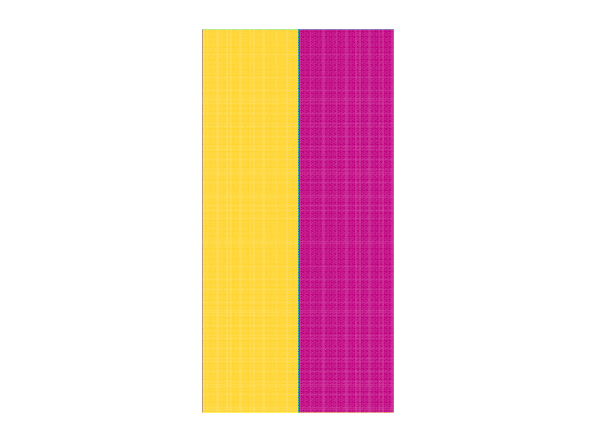}
      &
      \includegraphics[width=2.25cm]{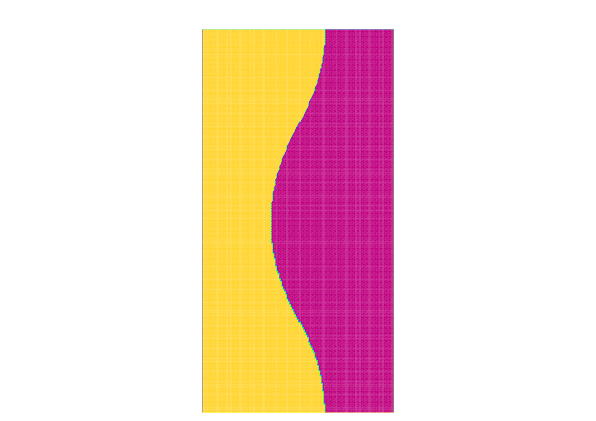}
      &
      \includegraphics[width=2.25cm]{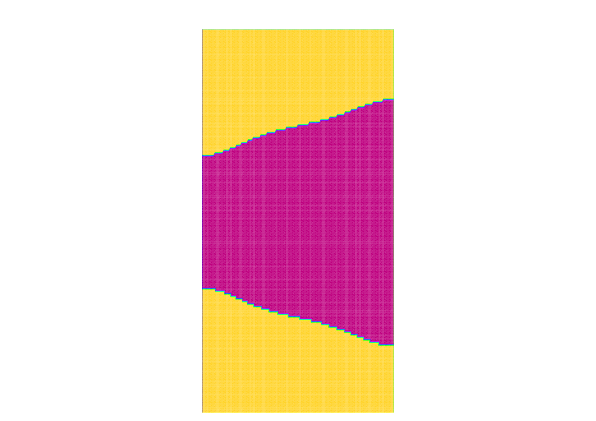}
      &
      \includegraphics[width=2.25cm]{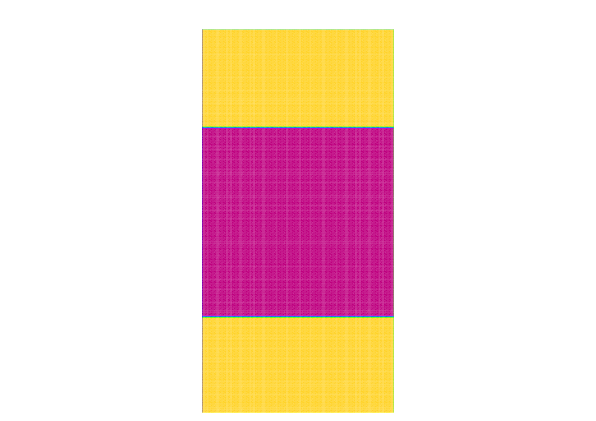}
      &
      \includegraphics[width=2.25cm]{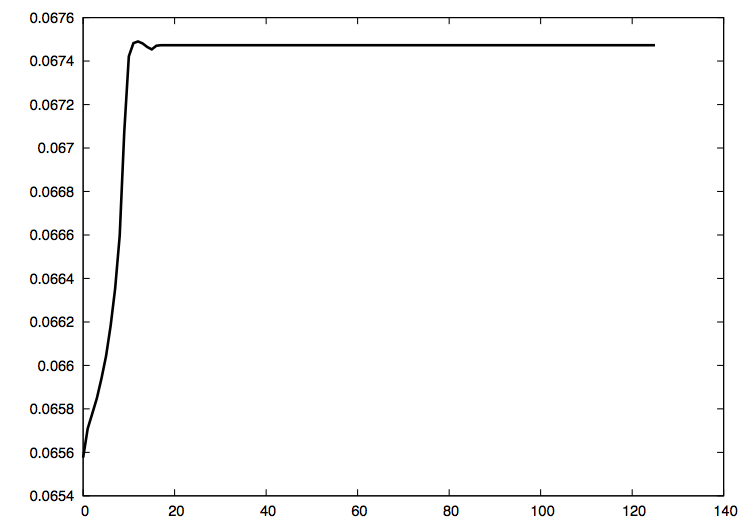}
      \cr
      &
      $t=0$ & $t=5$ & $t=10$ & $t=20$ \cr
      (b)
      &
      \includegraphics[width=2.25cm]{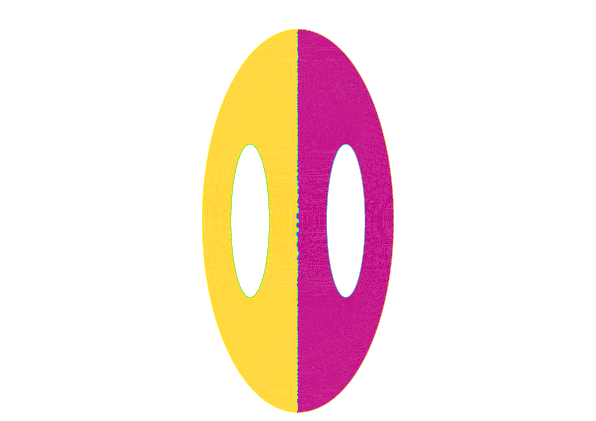}
      &
      \includegraphics[width=2.25cm]{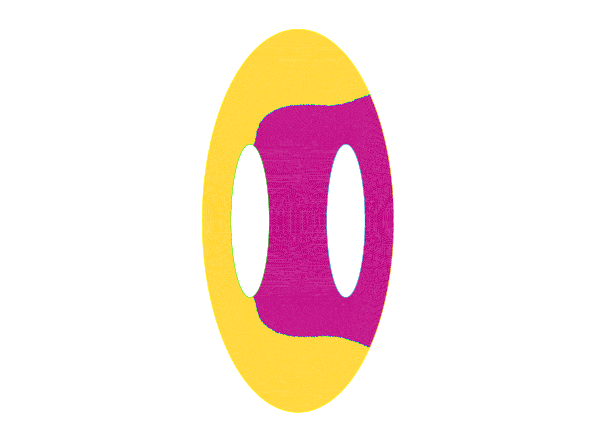}
      &
      \includegraphics[width=2.25cm]{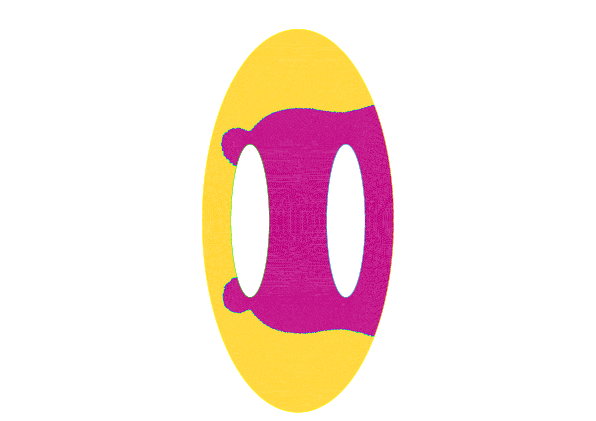}
      &
      \includegraphics[width=2.25cm]{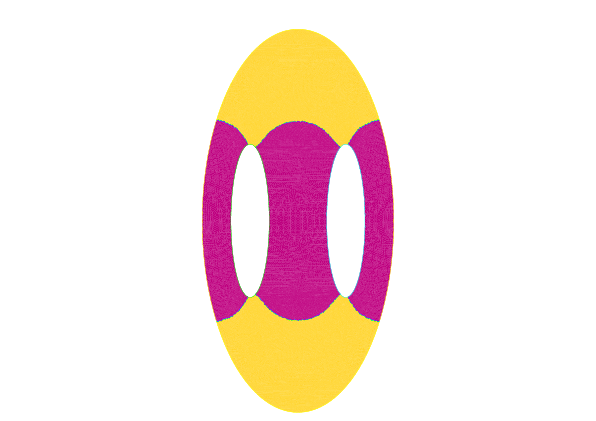}
      &
      \includegraphics[width=2.25cm]{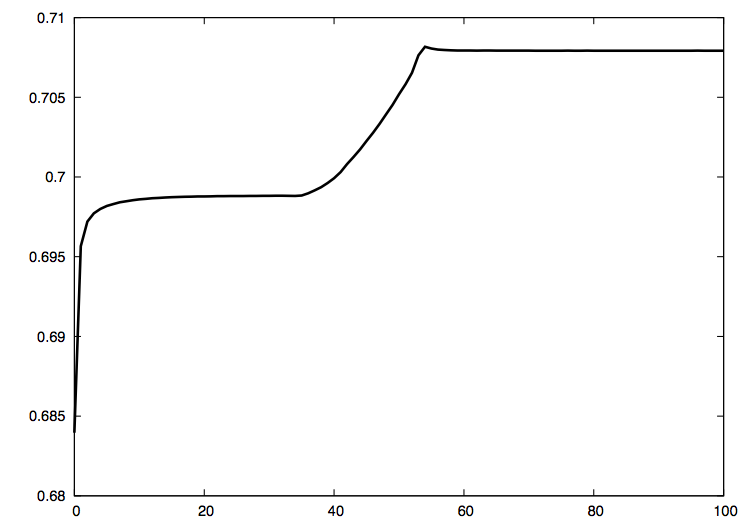}
      \cr
      &
      $t=0$ & $t=20$ & $t=40$ & $t=100$ \cr
      (c)
      &
      \includegraphics[width=2.25cm]{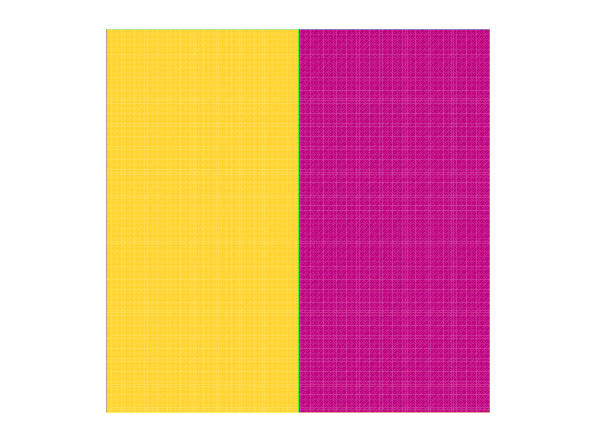}
      &
      \includegraphics[width=2.25cm]{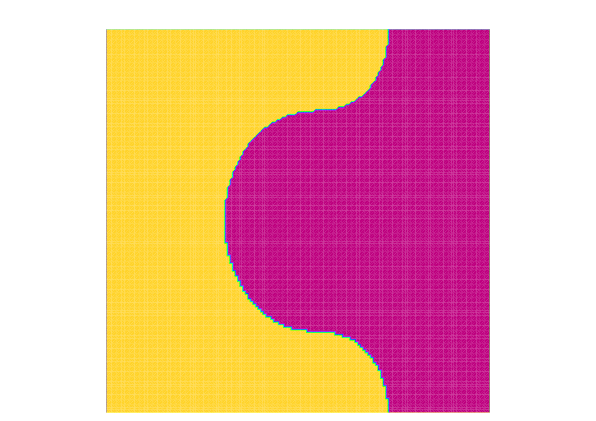}
      &
      \includegraphics[width=2.25cm]{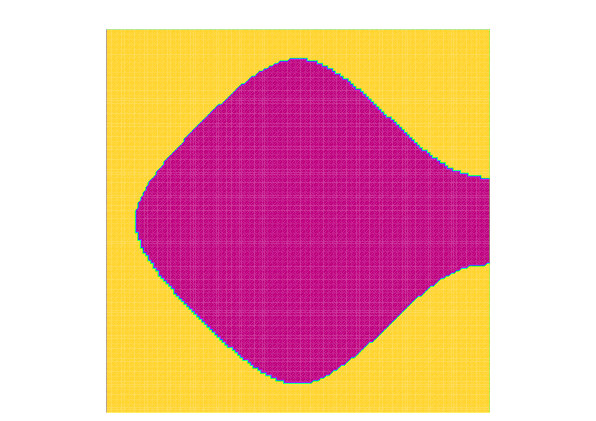}
      &
      \includegraphics[width=2.25cm]{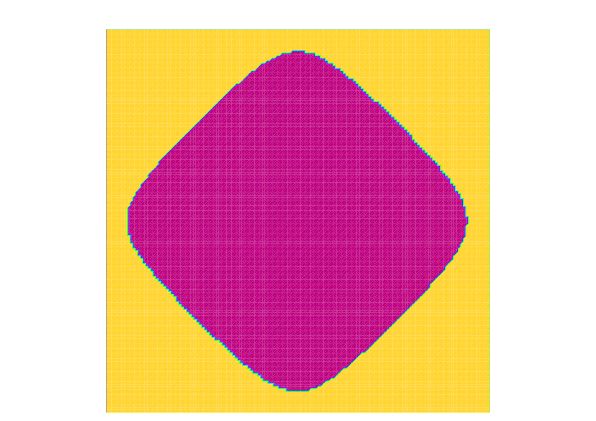}
      &
      \includegraphics[width=2.25cm]{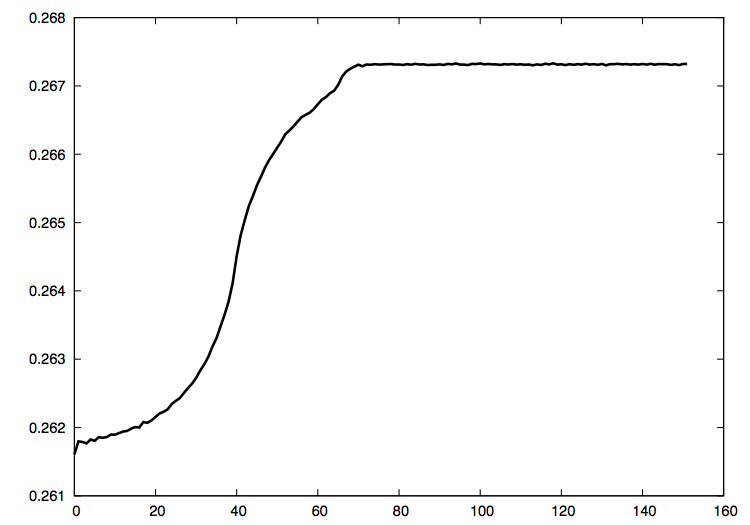}
      \cr
      &
      $t=0$ & $t=20$ & $t=60$ & $t=100$ \cr
    \end{tabular}
  \end{center}
  \caption{Problem \ref{prob-rho} with Neumann boundary}
  In all cases $c=1.1$ and $m_0 = 0.5$ are fixed.
  The rightmost graph is the evolution of $\lambda_1(\rho)$ as we solve (\ref{viscos}).
  (a) : Maximization of $\lambda_1(\rho)$ on $\Omega_5$.
  (b) : Maximization of $\lambda_1(\rho)$ on $\Omega_8$.
  (c) : Maximization of $\lambda_1(\rho)$ on $\Omega_1$.
  In this case the optimized eigenvalue $\lambda_1(\rho_{\max})$ has multiplicity two 
  (cf. Figure \ref{fig-rhoN-square} and Table \ref{table-simple-rho}) and hence the function $v_0$ in the level set evolution (\ref{SDF}) is 
  set $v_0(x) = -\{(c-1)|\nabla u_{1,\phi}(x)|^2 + (c-1)|\nabla u_{2,\phi}(x)|^2\}$ after the normalization $\int_{\Omega} |u_{i,\phi}(x)|^2 dx = 1$.
  \label{fig-rhoN-conv}
\end{figure}

\begin{figure}[H]
  \begin{center}
    \begin{tabular}{cccccc}
      (a)
      &
      \includegraphics[width=2.25cm]{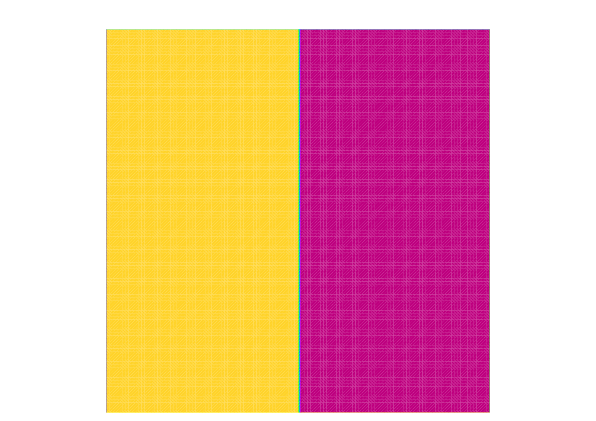}
      &
      \includegraphics[width=2.25cm]{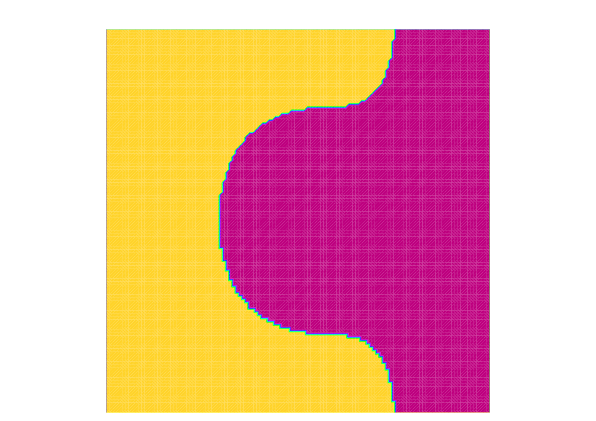}
      &
      \includegraphics[width=2.25cm]{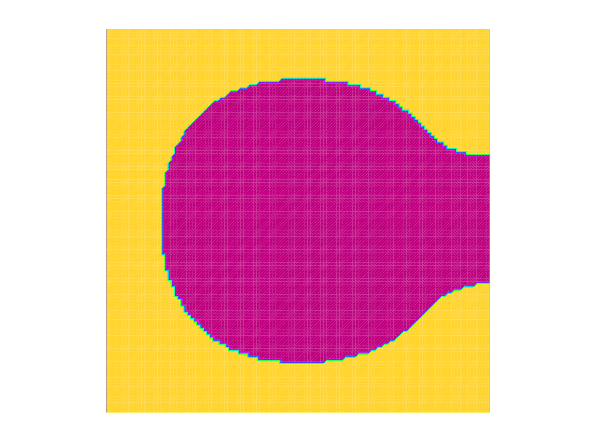}
      &
      \includegraphics[width=2.25cm]{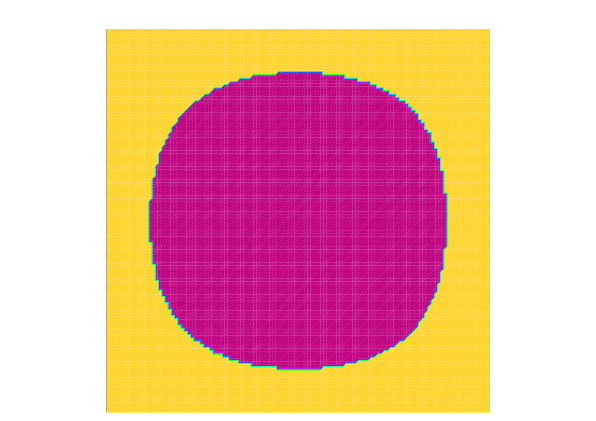}
      &
      \includegraphics[width=2.25cm]{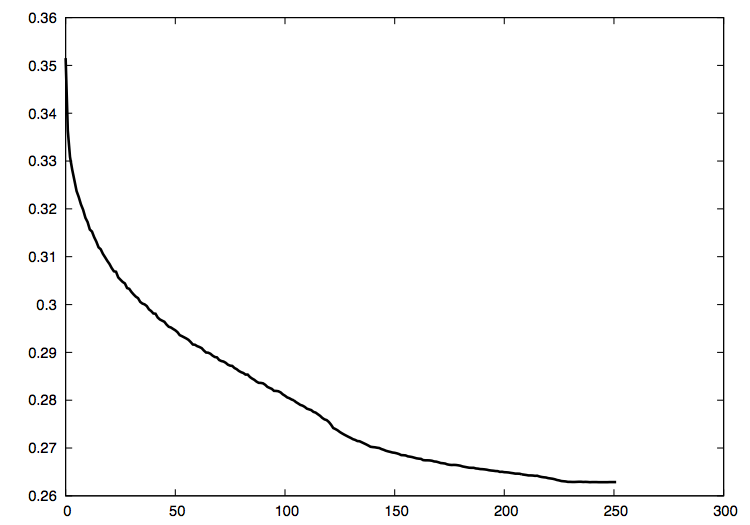}
      \cr
      &
      $t=0$ & $t=40$ & $t=180$ & $t=250$ \cr
      (b)
      &
      \includegraphics[width=2.25cm]{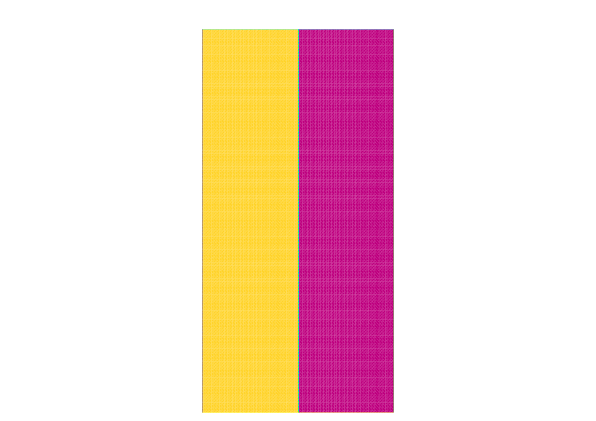}
      &
      \includegraphics[width=2.25cm]{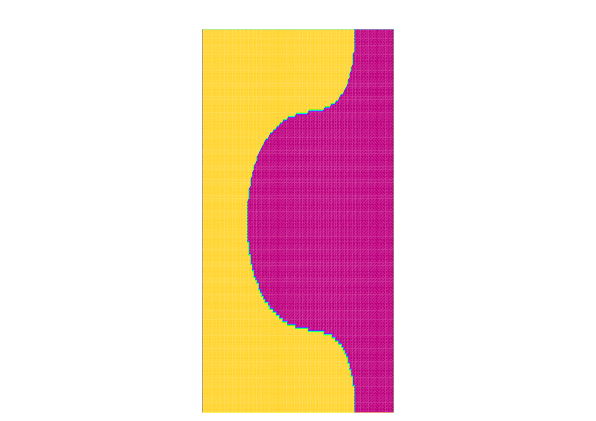}
      &
      \includegraphics[width=2.25cm]{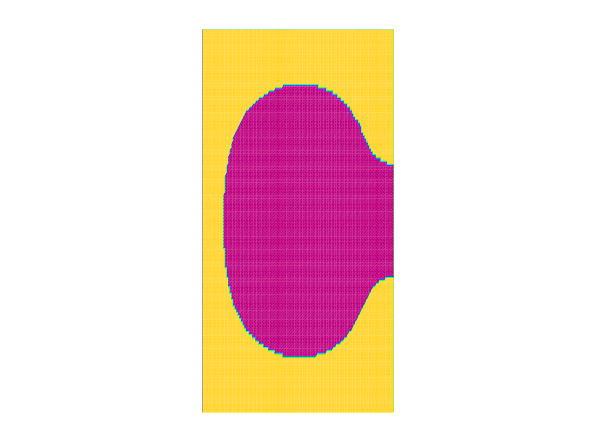}
      &
      \includegraphics[width=2.25cm]{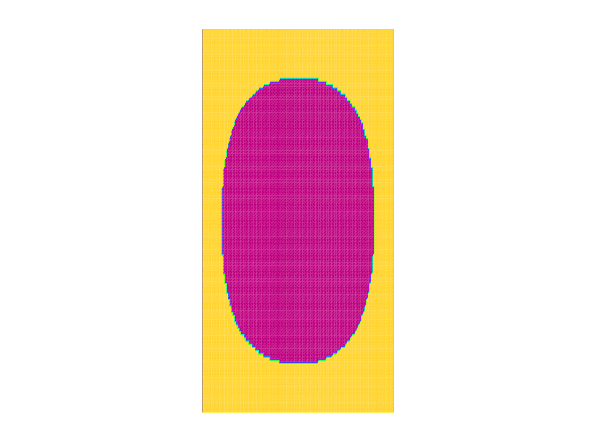}
      &
      \includegraphics[width=2.25cm]{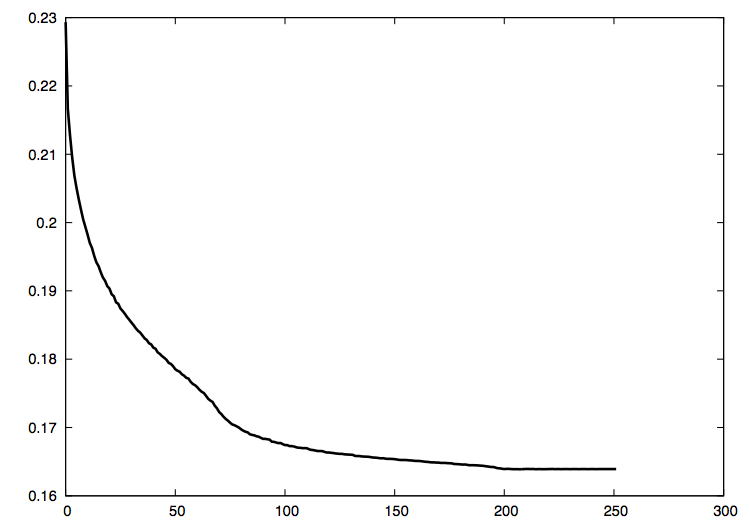}
      \cr
      &
      $t=0$ & $t=40$ & $t=150$ & $t=250$ \cr
      (c)
      &
      \includegraphics[width=2.25cm]{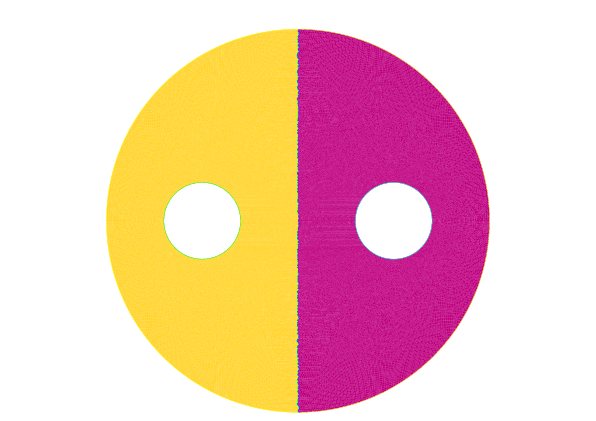}
      &
      \includegraphics[width=2.25cm]{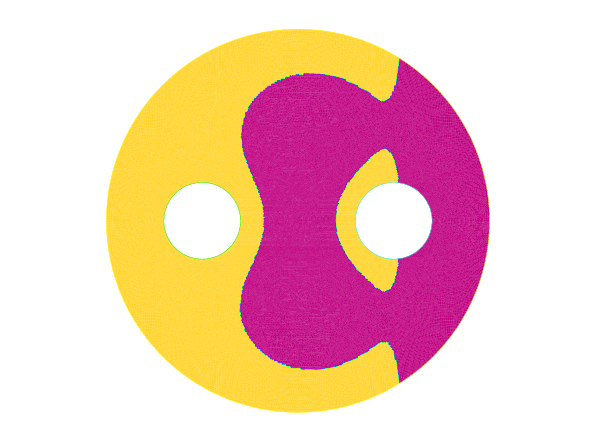}
      &
      \includegraphics[width=2.25cm]{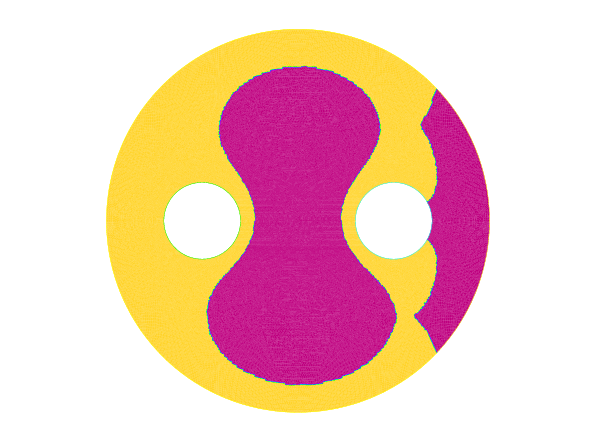}
      &
      \includegraphics[width=2.25cm]{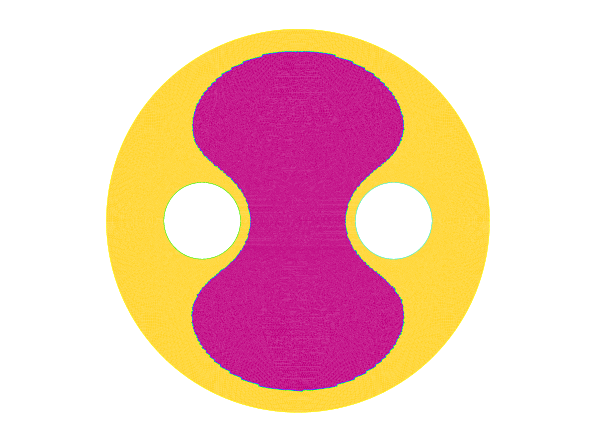}
      &
      \includegraphics[width=2.25cm]{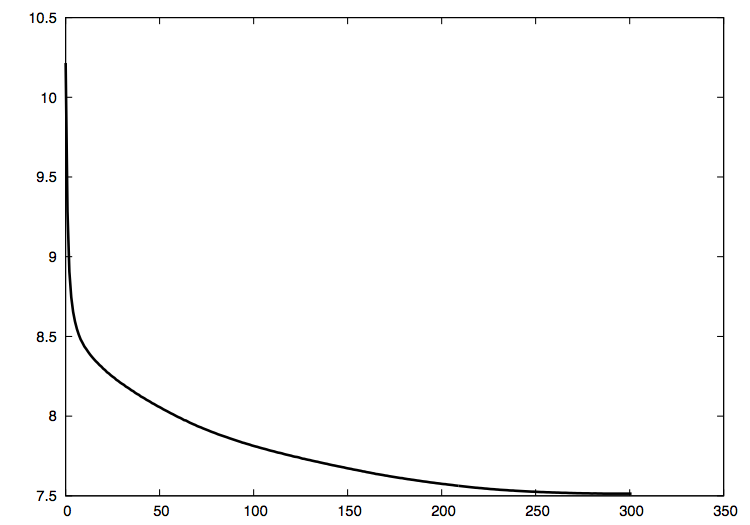}
      \cr
      &
      $t=0$ & $t=20$ & $t=100$ & $t=300$ \cr
    \end{tabular}
  \end{center}
  \caption{Problem \ref{prob-sigma} with Dirichlet boundary}
  \label{fig-sigmaD-conv}
  In all cases $c=2$ and $m_0 = 0.5$ are fixed.
  The rightmost graph is the evolution of $\mu_1(\sigma)$ as we solve (\ref{viscos}). 
  (a) : Minimization of $\mu_1(\sigma)$ on $\Omega_1$. 
  (b) : Minimization of $\mu_1(\sigma)$ on $\Omega_5$. 
  (c) : Minimization of $\mu_1(\sigma)$ on $\Omega_4$.
\end{figure}

\begin{figure}[H]
  \begin{center}
    \begin{tabular}{cccccc}
      (a)
      &
      \includegraphics[width=2.25cm]{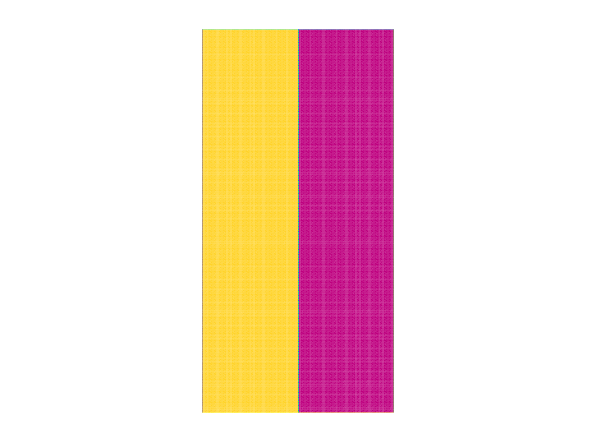}
      &
      \includegraphics[width=2.25cm]{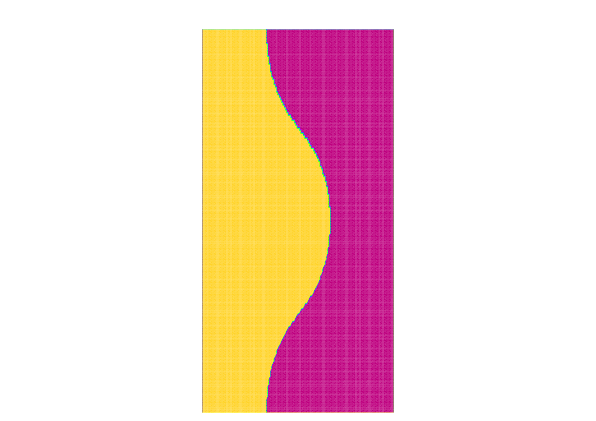}
      &
      \includegraphics[width=2.25cm]{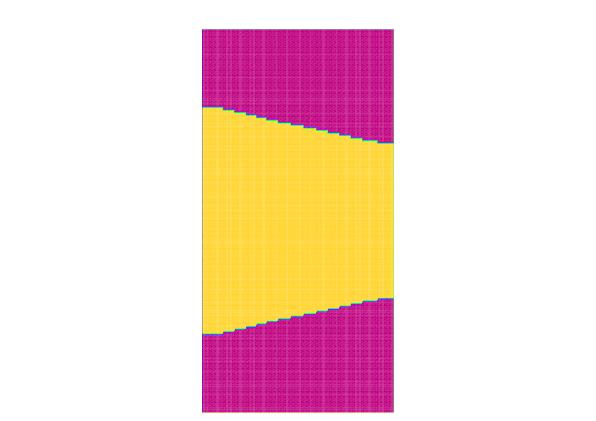}
      &
      \includegraphics[width=2.25cm]{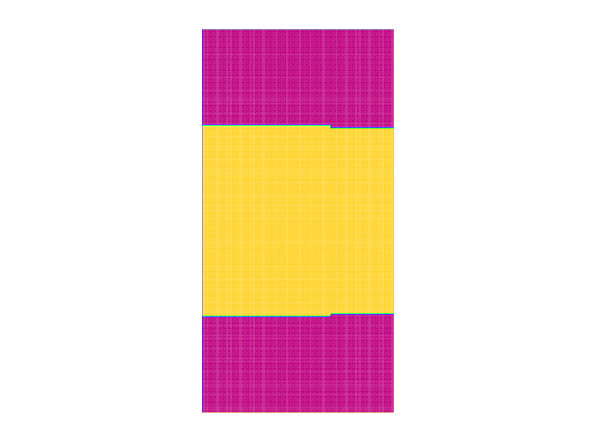}
      &
      \includegraphics[width=2.25cm]{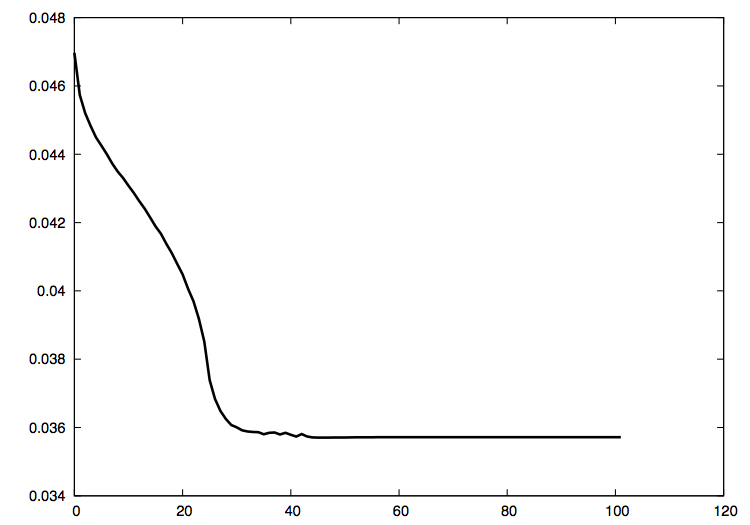}
      \cr
      &
      $t=0$ & $t=10$ & $t=30$ & $t=100$ \cr
      (b)
      &
      \includegraphics[width=2.25cm]{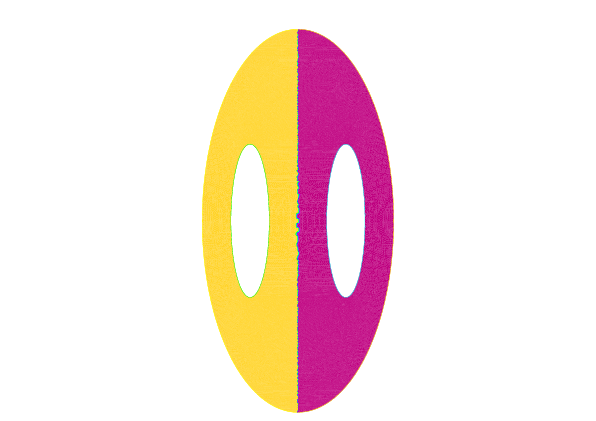}
      &
      \includegraphics[width=2.25cm]{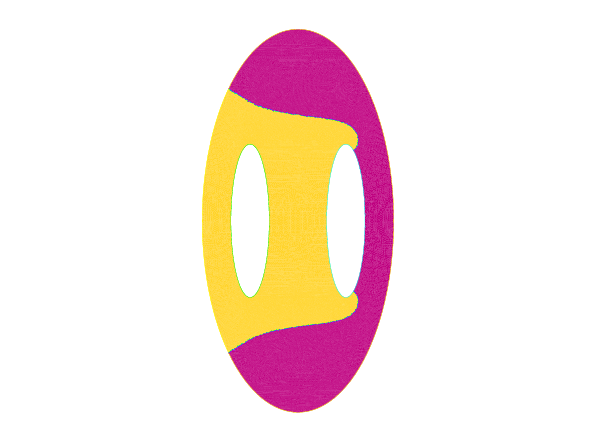}
      &
      \includegraphics[width=2.25cm]{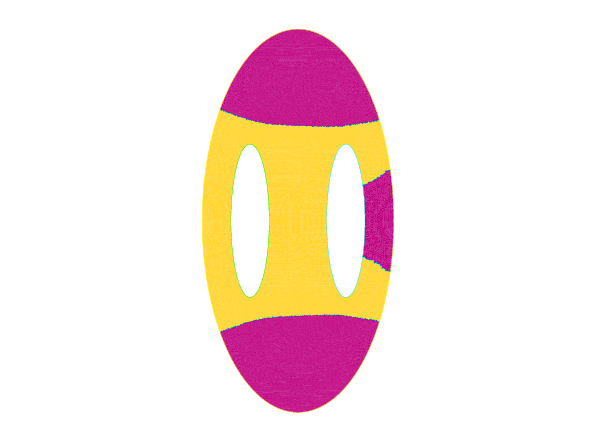}
      &
      \includegraphics[width=2.25cm]{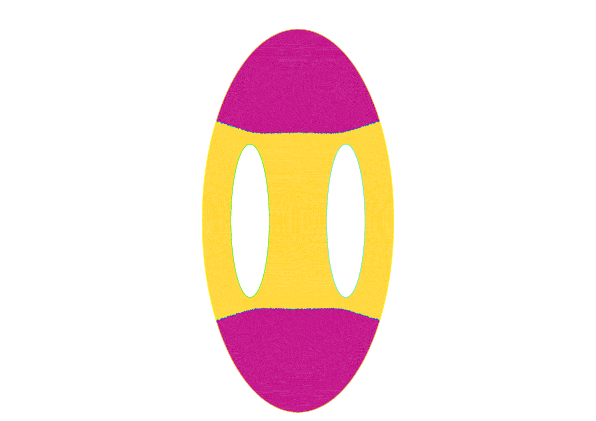}
      &
      \includegraphics[width=2.25cm]{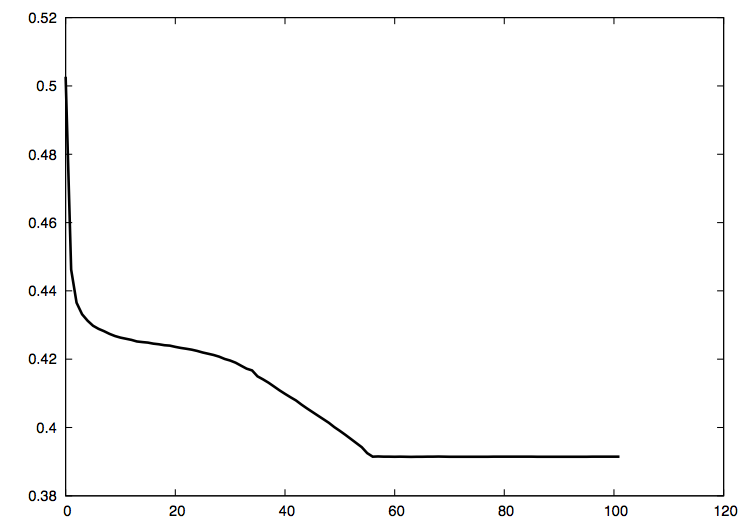}
      \cr
      &
      $t=0$ & $t=20$ & $t=40$ & $t=100$ \cr
      (c)
      &
      \includegraphics[width=2.25cm]{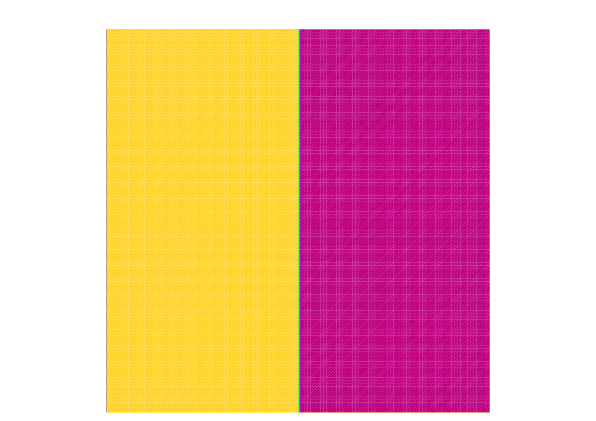}
      &
      \includegraphics[width=2.25cm]{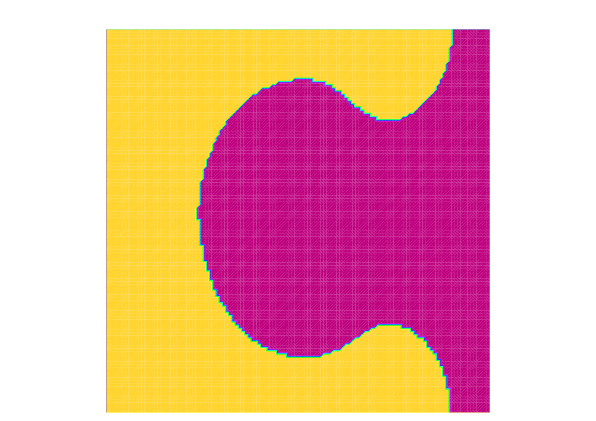}
      &
      \includegraphics[width=2.25cm]{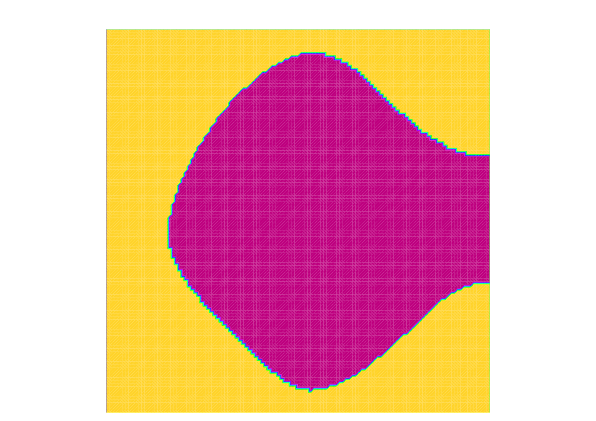}
      &
      \includegraphics[width=2.25cm]{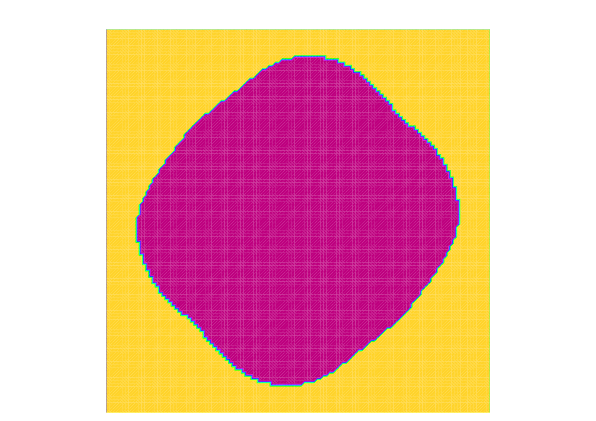}
      &
      \includegraphics[width=2.25cm]{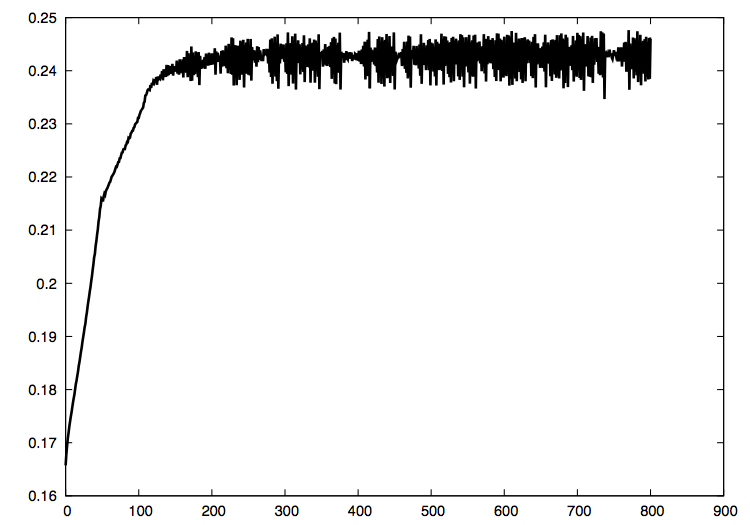}
      \cr
      &
      $t=0$ & $t=80$ & $t=200$ & $t=775$ \cr
    \end{tabular}
  \end{center}
  \caption{Problem \ref{prob-sigma} with Neumann boundary}
  \label{fig-sigmaN-conv}
  In all cases $c=2$ and $m_0 = 0.5$ are fixed.
  The rightmost graph is the evolution of $\mu_1(\sigma)$ as we solve (\ref{viscos}). 
  (a) : Minimization of $\mu_1(\sigma)$ on $\Omega_5$. 
  (b) : Minimization of $\mu_1(\sigma)$ on $\Omega_4$. 
  (c) : Maximization of $\mu_1(\sigma)$ on $\Omega_8$. 
  In this case the optimized eigenvalue $\mu_1(\sigma_{\max})$ has multiplicity two 
  (cf. Figure \ref{fig-sigmaN-square} and Table \ref{table-simple-sigma}) 
  and hence the function $v_0$ in the level set evolution (\ref{SDF}) is set 
  $v_0(x) = \{\mu_1(\phi)(c-1)|u_{1,\phi}(x)|^2 + \mu_2(\phi)(c-1)|u_{2,\phi}(x)|^2\}$ 
  after the normalization $\int_{\Omega} \sigma(x) |u_{i,\phi}(x)|^2 dx = 1$.
\end{figure}

Throughout numerical studies in this paper,
we numerically solved (\ref{viscos}) with
$\epsilon = 1.0\times 10^{-4}$ via the following implicit scheme
\begin{equation}
  \label{scheme}
  \int_\Omega \frac{\phi^n - \phi^{n-1}}{\Delta t^n} w_h dx = -\epsilon \int_\Omega \nabla \phi^n \cdot \nabla w_h dx 
  - \int_\Omega (v_0(\phi^{n-1})(x) + \nu)|\nabla \phi^{n-1}|w_h dx
\end{equation}
where $w_h$ is an arbitrary element of a finite element subspace of $H^1(\Omega)$ with suitable boundary conditions.
Fix the initial level set function by $\phi^0(x,y) = \phi(x,y) = x$.
Then we successively solve (\ref{scheme}) so that all $\phi^n\ (n\geq 0)$ keep the volume constraint $G(\phi^n)=0$ by following the discussion in \cite{OS}.
Here $c=1.1$ (Problem \ref{prob-rho}),
$c=2$ (Problem \ref{prob-sigma}) and $m_0 = 0.5$ are fixed in all cases,
and the step size $\Delta t^n$ is defined by $\Delta t^n = (\sup_{(x,y)\in \Omega}|\phi^n(x,y)|)^{-1}$ in each step.
Examples are shown in Figure  \ref{fig-rhoD-conv}-\ref{fig-sigmaN-conv}.
Although the initial shape of  $S=\sigma^{-1}(c)$ or $\rho^{-1}(c)$ matters,
we may observe that the convergence rate is independent of geometry of $\Omega$ in each problem with each boundary condition.

\section{Conclusion}
\label{section-discussion}
In this paper, we have studied the eigenvalue optimization of 
spatially inhomogeneous diffusion operator 
$A_{\rho}u = -\nabla \cdot (\rho(x) \nabla u)$
 with a given constraint,
which is motivated by the control of heat conductivity of spatially inhomogeneous media.
We applied the level set approach to characterize optimizers for Problem \ref{prob-rho}.
Collecting our numerical observations, one knows the following:
\begin{center}
  \begin{minipage}[t]{0.9\linewidth}
    The region $S_\ast =\{x\in \Omega \mid \rho_\ast(x) = c\}$ determined by the eigenvalue optimizer $\rho_\ast$ 
    is characterized by the super- or the sub-level set of $|\rho_\ast \nabla u_\ast|$ 
    even if $\lambda_1(\rho)$ has multiplicity greater than two, 
    where $u_\ast$ is the eigenfunction associated with $\lambda_1(\rho_\ast)$.
    This characterization is independent of topologies of $\Omega$ and boundary conditions on $\partial \Omega$.
    Moreover, 
    if $\Omega$ is star-shaped and has symmetry in a certain direction, 
    $S_\ast$ also possesses the same symmetry.
  \end{minipage}
\end{center}
One of key considerations in our numerical studies here is that 
{\emx eigenvalue optimizers can be characterized by associating eigenfunctions}
including symmetry. 
Once such characterizations are mathematically confirmed, 
various properties of optimizers, 
such as symmetry and continuous dependence on boundary conditions in a suitable topology,
will follow from those of corresponding eigenfunctions,
as is true in the case of Problem \ref{prob-sigma}.

\section*{Acknowledgements}
This research was partially supported by JST, CREST : 
A Mathematical Challenge to a New Phase of Material Science, Based on Discrete Geometric Analysis.
KM was partially supported by Coop with Math Program, a commissioned project by MEXT.
HN was partially supported by Grants-in-Aid for Scientific Research (C) (No. 26400067).
We thank Prof. Hideyuki Azegami for providing us with very meaningful advice for our computational study.
We also thank Prof. Motoko Kotani 
for introducing us to this study and giving us a lot of advice for writing this paper.


\appendix
\section{Tables}
In Problem \ref{prob-rho} (Problem \ref{prob-sigma}) with the homogeneous Neumann boundary condition, 
the simpleness of the first eigenvalue $\lambda_1(\rho)$ ($\mu_1(\sigma)$)  is nontrivial. 
$\lambda_{\min}$ and $\lambda_{\max}$ ($\sigma_{\min}$ and $\sigma_{\max}$) denote the minimizer and the maximizer of $\lambda_1(\rho)$ ($\mu_1(\sigma)$), respectively. 
As for the homogeneous Neumann boundary value problems for $\Omega = \Omega_5, \Omega_6, \Omega_7$ or $\Omega_8$, 
$\lambda_1(\rho)/\lambda_2(\rho) < 1$  ($\mu_1(\sigma)/\mu_2(\sigma) < 1$) holds for minimizations and maximizations, 
which implies that $\lambda_1(\rho)$ ($\mu_1(\sigma)$) is simple in those cases. 
\par
In the case of the maximization problem for $\Omega = \Omega_1$,
the ratio 
$\lambda_1(\rho_{\max})/\lambda_2(\rho_{\max})$ 
($\mu_1(\sigma_{\max})/\mu_2(\sigma_{\max})$) 
is close to $1$ and hence 
we may not consider that the first eigenvalue is simple.
Indeed, Figure \ref{fig-rhoN-square} and Figure \ref{fig-sigmaN-square} imply that (\ref{inequality-max-rho}) and (\ref{inequality-max-sigma}) do not hold, respectively. 
On the other hand,
$\lambda_2(\rho_{\max})/\lambda_3(\rho_{\max}) < 1$
and 
$\mu_2(\sigma_{\max})/\mu_3(\sigma_{\max}) < 1$
hold, 
which imply that the first eigenvalue for $\Omega_1$ has multiplicity two.
\begin{center}
  \begin{table}[h]
    \begin{center}
      \begin{tabular}{|c|c|c|c|c|} \hline
	{} & $\lambda_1(\rho_{\min}) / \lambda_2(\rho_{\min})$  & $\lambda_1(\rho_{\max})/\lambda_2(\rho_{\max})$ & $\lambda_2(\rho_{\max})/\lambda_3(\rho_{\max})$ & \\ 
        \hline
	$\Omega_1$ & $0.968758$ & $0.999951$ & $0.510619$ & Figure \ref{fig-rhoN-square}  \\
	$\Omega_5$ & $0.243209$ & $0.258043$ & ---- & Figure \ref{fig-rhoN-rect} \\
	$\Omega_6$ & $0.491202$ & $0.500293$ & ---- & Figure \ref{fig-rhoN-cross} \\
	$\Omega_7$ & $0.289396$ & $0.313473$ & ---- & Figure \ref{fig-rhoN-ellipse} \\
	$\Omega_8$ & $0.396238$ & $0.417894$ & ---- & Figure \ref{fig-rhoN-nonstar} \\ 
        \hline
      \end{tabular}
    \end{center}
    \caption{Multiplicity of $\lambda_1(\rho)$ ($c=1.1$, $m_0=0.5$).}
    \label{table-simple-rho}
  \end{table}
\end{center}
\begin{center}
  \begin{table}[h]
    \begin{center}
      \begin{tabular}{|c|c|c|c|c|} \hline
	{} & $\mu_1(\sigma_{\min})/\mu_2(\sigma_{\min})$  & $\mu_1(\sigma_{\max})/\mu_2(\sigma_{\max})$ & $\mu_2(\sigma_{\max})/\mu_3(\sigma_{\max})$ 
        & \\ 
        \hline 
	$\Omega_1$ & $0.605469$ & $0.999452$ & $0.424123$ & Figure \ref{fig-sigmaN-square}\\
	$\Omega_5$ & $0.221712$  & $0.338271$ & ---- & Figure \ref{fig-sigmaN-rect} \\
	$\Omega_6$ & $0.353771$ & $0.649946$ & ---- & Figure \ref{fig-sigmaN-cross}\\
	$\Omega_7$ & $0.236784$  & $0.440369$ & ---- & Figure \ref{fig-sigmaN-ellipse} \\
	$\Omega_8$ & $0.298349$  & $0.746472$ & ---- & Figure \ref{fig-sigmaN-nonstar}\\ 
        \hline
      \end{tabular}
    \end{center}
    \caption{Multiplicity of $\mu_1(\sigma)$ ($c=2$, $m_0=0.5$).}
    \label{table-simple-sigma}
  \end{table}
\end{center}


\section{Figures}
\subsection{1-dimensional case}
Optimization of $\lambda_1(\rho)$ in Problem \ref{prob-rho} on $\Omega_0 = (0,1)$. 
The optimal region $S_\ast$ and the graph of associated eigenfunction $u$ are drawn. 

\subsection{Dirichlet boundary condition}
For Problem \ref{prob-rho}, figures of eigenfunction stand for $|\rho \nabla u|^2$,
and 
for Problem \ref{prob-sigma}, figures of eigenfunction stand for $|u|^2$. Each figure shows optimizer of the first eigenvalue and associated eigenfunction.
The red region in figures of optimizer stands for 
$\{x \in \Omega \mid \rho(x) = c\}$ or $\{x \in \Omega \mid \sigma(x) = c\}$.
In the following figures, we calculate the optimal configuration for $c=1.1$ and $m_0 = 0.5$ unless otherwise noted.

\subsection{Neumann boundary condition}
For Problem \ref{prob-rho}, figures of eigenfunction stand for $|\rho \nabla u|^2$
and 
for Problem \ref{prob-sigma}, figures of eigenfunction stand for $|u|^2$. Each figure shows optimizer of the first eigenvalue and associated eigenfunction.
The red region in figures of optimizer stands for 
$\{x \in \Omega \mid \rho(x) = c\}$ or $\{x \in \Omega \mid \sigma(x) = c\}$.
In the following figures, we calculate the optimal configuration for $c=1.1$ and $m_0 = 0.5$ unless otherwise noted.

\subsection{Continuous dependency on boundary condition}
We calculate the dependency on boundary conditions.
The boundary condition is given by (\ref{robin-bc}).

\subsection{Convergence to optimizer}
We calculate the convergence to the optimizer.
Each figure shows the density function $\lambda$ and $\sigma$ after $t$ steps as we solve (\ref{viscos}) towards the optimizer 
in Problem \ref{prob-rho} and Problem \ref{prob-sigma}, respectively.

\end{document}